\documentclass[11pt]{article}
\usepackage{graphicx}
\usepackage{epstopdf}
\usepackage{amsmath,amssymb,amsthm,amsfonts}
\usepackage{color}
\usepackage{cite}
\usepackage{verbatim}
\usepackage{appendix}

\newtheorem{thm}{Theorem}[section]

\newtheorem{lem}[thm]{Lemma}
\newtheorem{prop}[thm]{Proposition}

\theoremstyle{definition}
\newtheorem{defn}{Definition}[section]
\theoremstyle{remark}
\newtheorem{rem}{Remark}[section]

\numberwithin{equation}{section}

\newcommand{\R}{{\mathbb{R}}}

\newcommand{\A}{{\mathcal{A}}}

\def\H{\mathcal{H}}
\def\A{\mathcal{A}}
\def\fh{\f{1}{2}}
\def\p{\partial}

\def\f{\frac}
\def\R{\mathbb R}
\def\C{\mathbb C}

\def\al{\alpha}

\def\Om{\Omega}

\def\va{\varepsilon}
\def\vaa{\varepsilon_a}
\def\ap{\alpha_a}

\def\l{\left}
\def\r{\right}

\def\inte{\int_{\mathbb{R}^N}}
\def\ovl{\overline}

\begin{document}
\title{Qualitative Properties of Normalized Solutions for Magnetic Gross-Pitaevskii Equations with Anharmonic Potentials}
\author{Yujin Guo\thanks{Email: yguo@ccnu.edu.cn. }, Yan Li\thanks{Email: yanlimath@mails.ccnu.edu.cn.},
Yong Luo\thanks{Email: yluo@ccnu.edu.cn. },   and Shuangjie Peng\thanks{Email: sjpeng@ccnu.edu.cn. }\\
\small \it	School of Mathematics and Statistics,\\
\small \it  Key Laboratory of Nonlinear Analysis {\em \&} Applications (Ministry of Education),\\
\small \it Central China Normal University, Wuhan 430079, P. R. China\\}
%\footnotetext{This paper is supported by National Key R\&D Program of China (2022YFA1006900).}

\date{\today}
\smallbreak  \maketitle

\begin{abstract}
This paper is concerned with normalized solutions of the magnetic focusing Gross-Pitaevskii equations with   anharmonic potentials in $\R^N$, where $N=2$ or $3$. We construct axially symmetric normalized concentrating solutions as
the parameter \(a>0\) approaches \(a_*(N)\), where \(a_*(N)\geq0\) is a
critical constant depending only on \(N\). We further prove that up to a constant phase (and a
rotational transformation for $N=2$), normalized concentrating solutions are unique and
axially symmetric as \(a\to a_*(N)\).
When \(N=3\), we also prove that  the corresponding unique normalized concentrating solution is free of vortices  as \(a\to a_*(3)\), even if the anharmonic potential is non-radially symmetric.

\end{abstract}


%\noindent {\it MSC(2010):} ***, ***, ***
\tableofcontents

\bigskip

\section{Introduction}
In this paper, we study the following   magnetic focusing Gross-Pitaevskii (GP) equation with a trapping potential in $\R^N$ ($N=2, 3$):
\begin{equation}\label{GP-1}
\left\{
\begin{aligned}
&-\Delta u+V(x)u+i\Om(x^\perp\cdot\nabla) u-a|u|^2u=\lambda _a u\quad\hbox{in}\,\,\ \R^N,\ \ a>0,\ \ \lambda_a\in\R, \\
&\inte |u|^2dx=1,\ \ u\in\mathcal{H},
\end{aligned}
\right.
\end{equation}
where $i$ denotes the imaginary unit, $\H$ is defined by
\begin{equation}\label{N1.2}
\H:=\l\{u\in H^1(\R^N,\C): \,\ \inte \Big[|\nabla u|^2+\big(1+V(x)\big)|u|^2\Big]<\infty\r\},
\end{equation}
and $V(x): \R^N\to\R$ denotes a trapping potential. Here we denote for $x=(x_{1},\cdots,x_{N})\in \R^N$,
\begin{equation}\label{1.1M}
x^\perp:=\left\{
\begin{aligned}
&(-x_2,x_1),&\   \hbox{if}\ \ N=2;\\
&(-x_2,x_1,0),&    \hbox{if}\ \ N=3;
\end{aligned}
\right.
\end{equation}
where the fixed constant $0\le \Om<\infty $ represents the rotational speed of the magnetic trap $V(x)$.
The equation (\ref{GP-1}) is used (cf. \cite{bec1,bec2,Brad1,Brad2,Coo}) to describe the trapped attractive Bose-Einstein condensates (BECs) under rotation, where the parameter $a>0$  characterizes the strength of attractive interactions among ultracold atoms,  and $\lambda _a:=\lambda_a(a,\Omega)\in\R$ is a suitable Lagrange multiplier such that the $L^2$-constraint condition is satisfied.

For the non-rotational  case $\Om=0$,  concentrating solutions of \eqref{GP-1} were well studied recently, where the equation \eqref{GP-1} is essentially mass critical for $N=2$, while it is mass supercritical for $N\ge 3$. Particularly, when $N=2$, Guo et al. \cite{GZZ} studied the existence and nonexistence of ground state solutions for \eqref{GP-1}  with $\Om=0$ by the variational methods, where they also analyzed the concentrating behavior and  local uniqueness of  ground state solutions. Moreover, if $V(x)$ admits at least one isolated critical point $x_0$, Pellacci et al. \cite{Pistoia}  employed the reduction method to obtain the existence of normalized solutions, which  concentrate at $x_0$, for \eqref{GP-1} with $\Om=0$.
More recently, suppose $x_0$ is a non-isolated critical point of $V(x)$ degenerating   at $N-1$ directions,
%\begin{equation}\label{N1.3}
%\f{\partial^2V(x_0)}{\partial \nu^2}\neq 0\quad\hbox{and}\quad\det\l(\f{\partial^2\Delta V(x_0)}%{\partial \tau_l\partial \tau_j}\r)_{1\leq l,j\leq N-1}\neq 0,
%\end{equation}
%where $\nu$ is the  outward unit normal and $\tau_j$ is the $j$-th principal tangential unit vector,
Luo et al. \cite[Theorems 1.3 and 1.4]{LPWY}  obtained the existence and uniqueness of normalized solutions,  which concentrate  at $x_0$, for \eqref{GP-1} with $\Om=0$. We also refer to  \cite{Jean,Bart2,CXZ} and the references therein for analyzing the saddle solutions of \eqref{GP-1} with $\Om=0$.

%By using the rearrangement symmetry method, \cite{Lieb2} proved ground state  solution (i.e., minimizer) must be radially symmetrical. As for  positive solution, one can using the classical moving plane method developed by \cite{GNN} to prove that any positive solution of \eqref{GP-1} is radially symmetrical in the case where $\mu<0$. Based on these works mentioned above \cite{Lieb2,GNN}, \cite{Bye,Kwon} using some ODE methods proved the uniqueness of positive solution equation \eqref{GP-1} without $L^2$-constrain.

As for the rotational case $\Om>0$,  complex-valued solutions of  \eqref{GP-1} were studied recently in \cite{BoN,BarT1,Bao,Este,GLY,GLP2,Lewin} and the references therein. However, most existing results of  \eqref{GP-1} are mainly restricted to the following harmonic potential
\begin{equation}\label{hampot}
 V(x)=\sum_{j=1}^Nk_j^2x_j^2,\ \,\ k_j>0, \ \ x=(x_{1},\cdots,x_{N})\in \R^N.
\end{equation}
Under the assumption \eqref{hampot}, one can check that \eqref{GP-1}  admits a finite critical rotational speed $$\Omega^*:=2\min\big\{k_1,k_2,\cdots,k_N\big\},$$
which can be defined as
\begin{equation}\label{Omega}
\Omega^*=\Omega^*\big(V(x)\big):=\sup\big\{\Omega>0:\  V(x)-\f{\Omega^2}{4}{|x|^2}\to\infty\ \ \hbox{as}\,\ |x|\to\infty\big\}.
\end{equation}
Suppose $0<\Omega<\Omega^*$ is fixed,
Lewin-Nam-Rougerie \cite{Lewin} and Guo-Luo-Yang \cite{GLY} studied the existence and concentrating behavior of ground state solutions for \eqref{GP-1} with $N=2$. Moreover, Guo and his collaborators  \cite{GLY,GLP2} also proved the uniqueness and nonexistence of vortices for ground state solutions of  \eqref{GP-1} for the case where $V(x)=|x|^2$ and  $N=2$.
While for $ \Omega>\Omega^*$, there is no ground state solution of \eqref{GP-1}, see \cite{GLY}. As for $\Omega=\Omega^*$, the existence and nonexistence of  normalized solutions  for \eqref{GP-1} are more complicated, see \cite{GLP2}, which depend strongly on the shape of $V(x)$ and the value of $a>0$.
It should be pointed out that  the existing works mentioned above were mainly studied through the variational approaches.

%Luo-Yang \cite{L-Y} proved the existence of  solutions of \eqref{GP-1} via the mountain pass lemma for $N=3$.

In order to study the equation \eqref{GP-1} at arbitrarily large rotational speed $\Om>0$, physicists (cf. \cite{Fett}) usually replace the harmonic potential \eqref{hampot} by the following stiffer anharmonic potential
\begin{equation}\label{ahm}
V(x)=|x|^4+\sum_{j=1}^N\kappa_jx_j^2,\ \ \kappa_j\in\R, \ \ x=(x_{1},\cdots,x_{N})\in \R^N.
\end{equation}
Anharmonic trapping potentials of the type \eqref{ahm} were achieved experimentally in \cite{Bret,Bret2}, where
the quadratic trapping potential was modified (cf. \cite{Fett2,R1,Lun2,Af1}) by superimposing a blue detuned laser beam to
the magnetic trap holding the atoms.

%use the reduction method to

Motivated by above mentioned works \cite{BoN,GLY,LPWY,Fett}, in this paper we shall mainly study the normalized solutions of  \eqref{GP-1} under the  following anharmonic trap
\begin{equation}\label{pot}
V(x)=\left\{
\begin{aligned}
&|x|^4+\big(\f{\Om^2}{4}-2\big)|x|^2,&  \ \hbox{if}\ \ N=2;\\
&|x|^4+\big(\f{\Om^2}{4}-2\big)(x_1^2+x_2^2)-2x_3^2,&  \, \hbox{if}\ \ N=3;
\end{aligned}
\right.
\end{equation}
where $0\leq\Om<\infty$ is fixed. We remark that the critical rotational speed $\Om^*$ defined in (\ref{Omega}) satisfies  $\Om^*=\infty$ under the assumption (\ref{pot}). By choosing the anharmonic trap \eqref{pot}, the equation \eqref{GP-1} can be equivalently rewritten as the following form
\begin{equation}\label{GP-2}
\left\{
\begin{aligned}
&\left(-\Delta_{x^{\Omega}}\right) u+(|x|^2-1)^2u-a|u|^2u=\mu_a u\ \ \hbox{in}\, \ \R^N,\ \ a>0, \ \ N=2,3, \\
&\inte |u|^2dx=1,\ \ u\in\mathcal{H},
\end{aligned}
\right.
\end{equation}
where $0\leq\Om<\infty$ is fixed, the suitable Lagrange multiplier $\mu _a=\lambda _a+1\in \R$ depends on $a>0$ and $\Omega\ge 0$, so that the $L^2$-constraint condition is satisfied, and the magnetic Laplacian operator $-\Delta_{x^{\Omega}}$ is defined by
\begin{equation*}
-\Delta_{x^{\Omega}}:=-(\nabla-ix^{\Omega})^2=-\Delta+2ix^{\Omega}\cdot\nabla+|x^{\Omega}|^2,
\end{equation*}
together with the fact that for $x=(x_{1},\cdots,x_{N})\in \R^N$,
\begin{equation}\label{1.1}
 x^\Omega:=\left\{
\begin{aligned}
&\frac{\Om}{2}(-x_2,x_1),&\   \hbox{if}\ \ N=2;\\
&\frac{\Om}{2}(-x_2,x_1,0),&    \hbox{if}\ \ N=3.
\end{aligned}
\right.
\end{equation}
The main features of the ring-shaped potential $(|x|^2-1)^2$ for (\ref{GP-2}) are as follows: it admits  non-isolated critical points on the sphere $|x|=1$, and
its directional derivative   along the tangent direction
is always zero at each critical point lying on  the sphere $|x|=1$. Thus, the critical points of the ring-shaped potential $(|x|^2-1)^2$ satisfy neither  the non-degenerate conditions of \cite{Pistoia}, nor the assumptions of \cite[Theorem 1.3]{LPWY}. Therefore,  as far as we know, the existing approaches such as \cite{LPWY,Pistoia} cannot be directly applied to the analysis of (\ref{GP-2}).

We now consider the following equation
\begin{equation}\label{scal equ}
-\Delta u+u-|u|^2u=0\quad\hbox{in}\ \ \R^N,\ \  u\in H^1(\R^N,\C).
\end{equation}
It follows from  Kurata \cite{Kurata} that the
least energy solution of \eqref{scal equ} has the form
$Q_{x_0}(x)=Q(x-x_0)e^{i\sigma}$,
where $\sigma\in\R$, $x_0\in\R^N$, and  $Q=Q(|x|)>0$ is (cf. \cite{Kwon}) the unique positive  solution of the following nonlinear scalar field equation
\begin{equation}\label{A1}
-\Delta u+u- u^3=0\quad\hbox{in}\ \ \R^N,\ \ u\in H^1(\R^N,\R).
\end{equation}
%It follows from   \cite{Kwon} that the above equation admits a unique positive solution $Q=Q(|x|)$.
Note from \cite{Cing} that $Q(x)$ is non-degenerate, in the sense that its linearized operator
\begin{equation}\label{lin op}
L:=-\Delta+1-Q^2-2\operatorname{Re}(Q\cdot)Q \ \,\ \hbox{in}\, \ L^2(\R^N,\C)
\end{equation}
satisfies
\begin{equation}\label{nond}
ker L=span\Big\{iQ,\,\frac{\partial Q}{\partial x_1},\,\cdots,\,\frac{\partial Q}{\partial x_N}\Big\}.
\end{equation}
Furthermore, we get  from \cite[Lemma 8.1.2] {CAZ} and \cite[Proposition 4.1]{GNN}  that $Q(x)$ satisfies
\begin{equation}\label{N1:1}
\inte Q^2dx=\f{4-N}{4}\inte Q^4dx=\f{4-N}{N}\inte |\nabla Q|^2dx,
\end{equation}
and
\begin{equation} \label{1:3}
Q(x) \, , \ |\nabla Q(x)| = O(|x|^{-\frac{N-1}{2}}e^{-|x|})
\ \ \text{as} \  \ |x|\to \infty.
\end{equation}
Recall also from  \cite{LL} the following diamagnetic inequality
\begin{equation}\label{1:4}
\big|(\nabla-i\A)u\big|^2\geq\big|\nabla |u|\big|^2
\ \ \hbox{a.e. in}\ \ \R^N, \  \ u\in H^1(\R^N,\C),
\end{equation}
where $\A\in L^2_{loc}(\R^N,\R^N)$ is any given vector function.

%Some adaptations need to be made on those existing reduction approaches to obtain the qualitative properties of concentration  solutions for \eqref{GP-2}.

%Note that the extra non-degenerate assumptions on the critical point $P$ of the potentials  are  essential  when using reduction method.

\subsection{Main results}

In this subsection, we shall introduce the main results of the present paper. Throughout the whole paper, we always denote
 \begin{equation}\label{A:1.16}
\arraycolsep=1.5pt
 a^*:=\inte Q^2dx,\ \ a_{*}(N):=\left\{\begin{array}{lll} &a^*,  \ \ &\mbox{if}\ \ N=2;\\[3mm]
&0, \,\ &\mbox{if}\ \ N=3;\\[2mm]
\end{array}\right.
\end{equation}
and
\begin{equation}\label{A:1.17}
\arraycolsep=1.5pt
a\to a_{*}(N)\ \ \Longleftrightarrow\ \ \left\{\begin{array}{lll}
&a\nearrow a^*,  \  \ &\mbox{if}\ \ N=2;\\[3mm]
&a\searrow0, \,\ &\mbox{if}\ \ N=3.\\[2mm]
\end{array}\right.
\end{equation}
The main purpose of the present paper is to address qualitative properties, including the existence, axial symmetry and vortices, of normalized  concentrating solutions $w_a$ of (\ref{GP-2}) as $a\to a_*(N)$, where $0\leq \Omega<\infty$ is always fixed, and $\mu_a=\mu_a(a,\Omega)<0$ depends on $a>0$ and $\Omega \ge 0$. Towards this purpose, we define

\begin{defn} {\em  $w_a\in\H $ is called an axially symmetric solution of (\ref{GP-2}), if $w_a$ satisfies $$w_a(x_1, \cdots, x_{N-1}, x_N)\equiv  w_a(r', x_N)\ \  \mbox{in}\, \ \R^N ,\ \  r'=\sqrt{x_1^2+ \cdots +x_{N-1}^2}.$$  }
\end{defn}

%For  $a>0$ and $\mu_a<0$, denote
%\begin{equation}\label{A.6}
%\vaa:=\f{1}{\sqrt{-\mu_a}}>0,\ \  u_a:=\sqrt a\vaa w_a(\vaa x+x_a)e^{-i\sigma_a-i\vaa x\cdot x_a^\Om},
%\end{equation}
%where $x_a\in\R^N$ is the global maximum point of $|w_a|$, and $\sigma_a\in[0,2\pi)$ is a constant phase.
%%$x_{a}=(0,\cdots,0,p_{a})\in\R^N$ is the unique global maximum point of $|w_{a}(x)|$, and
%%$\sigma_{a}\in[0,2\pi)$ is chosen properly such that
%%\begin{equation}\label{1:23}
%%Re\Big(\inte u_{a}(iQ)\Big)=0.
%%\end{equation}
%We derive  from \eqref{GP-2} that $u_a(x)$ satisfies
%\begin{equation}\label{1NGP-5}
%\left\{
%\begin{split}
%& (-\Delta_{\vaa^2x^{\Omega}}) u_{a}+\Big[1+\vaa^2\big(|\vaa x+ x_a|^2-1\big)^2\Big]u_{a}-|u_{a}|^2u_{a}=0\ \ \hbox{in}\ \ \R^N,\\
%&\inte |u_a|^2dx=a\vaa^{2-N}.
%\end{split}
%\right.
%\end{equation}
%We treat \eqref{1NGP-5} a perturbation of \eqref{scal equ}, and consider solutions approximated by a suitable scaling of $Q(x)$, with $\vaa>0$ sufficiently small. In view of above facts and inspired by
%%One can note formally that the equation of (\ref{1NGP-5}) is a perturbation of \eqref{scal equ} as $a\to a_*(N)$. By choosing a suitable constant phase $\sigma_a\in[0,2\pi)$, it is thus natural to assume that
%\begin{equation}\label{A.5}
%u_a(x)=Q(x)+v_a(x),
%\end{equation}

\begin{defn}\label{def2} {\em  $w_a\in\H $ is called a normalized  concentrating solution of (\ref{GP-2}) concentrating at a point $P\in\R^N$ as $a\to a_*(N)$, if there exists $\mu_a<0$ such that $w_a$ satisfies (\ref{GP-2}) and
\begin{equation}\label{sigpek}
\sqrt {\f{a}{-\mu_a}} w_a\Big(\f{x}{\sqrt{-\mu_a}} +x_a\Big)e^{-i\sigma_a-\f{i}{\sqrt{-\mu_a}} x\cdot x_a^\Om}=Q(x)+v_a(x)\ \ \hbox{as}\ \ a\to a_*(N),
\end{equation}
where
\begin{equation*}
  x_a\to P\ \ \hbox{and}\ \ \mu_a\to-\infty\ \ \hbox{as}\ \  a\to a_*(N),
\end{equation*}
and $v_a(x)$ is a lower order term in some suitable sense.
Here  $x_a\in\R^N$ is a global maximum point of $|w_a|$, and
 $\sigma_a\in[0,2\pi)$ is a constant phase.}
\end{defn}

%Throughout the whole paper, we always consider the limit $a\to a_{*}(N)$ in the sense of (\ref{A:1.17}).
Employing the  reduction method, in this paper we first establish the following existence of normalized concentrating solutions for \eqref{GP-2}.

%
%For convenience, we define the set $\mathcal{F}$ of axially symmetric functions as follows:
%\begin{equation}\label{wk space}
%\mathcal{F}:=\Big\{%u:\,\
%u\in \H:\ u(x_1,x_2,\cdots,x_N)=u(r^\prime,x_N),\ \ r^\prime=\big(\sum^{N-1}_{i=1}x_i^2\big)^\frac{1}{2}\Big\}\subset\H,
%\end{equation}
%where $\mathcal{H}$ is defined by \eqref{N1.2}.

\begin{thm}\label{thm1}
Let $0\leq\Omega<\infty$ be fixed, and suppose $a>0$ satisfies
\begin{equation}\label{thm1:1}
 \ \, a\in (a^*-\delta,a^*), \ \  if \ \ N=2;  \ \ a\in(0,\delta),   \ \ if \ \ N=3,
\end{equation}
where $\delta=\delta (\Om)>0$ is a sufficiently small constant. Then \eqref{GP-2} admits at least one normalized concentrating solution concentrating at $P_0=(0,\cdots,0,1)$, which is also  axially symmetric.
\end{thm}

To the best of our knowledge, this is the first work on the existence and axial symmetry of normalized concentrating solutions  for magnetic GP equations.
%Since the origin is another critical point of the trapping potential $(|x|^2-1)^2$ for \eqref{GP-2}, we expect that the same argument of  Theorem \ref{thm1} can also yield the  existence of radially symmetric solutions of \eqref{GP-2} concentrating at the origin, and its proof is however left to the interested reader.
 %On the other hand,
To prove Theorem \ref{thm1}, we shall first analyze the existence of axially symmetric solutions $w_a$ for a magnetic GP equation without constraint conditions, based on which we then address the refined mass estimates of $\sqrt{a}w_a$ as $a\to a_*(N)$ by establishing Pohozaev identities. Since the non-isolated critical points of   $(|x|^2-1)^2$ do not satisfy the non-degenerate conditions of \cite[Theorem 1.3]{LPWY}, the method of \cite[Theorem 1.3]{LPWY} does not directly yield the
existence of axially symmetric solutions \(w_a\).  To overcome this difficulty, we shall make full use of the radial symmetry of the  potential $(|x|^2-1)^2$ for \eqref{GP-2}. Moreover,  comparing with the existing works \cite{LPWY,Pistoia}, which handled the GP equation without the magnetic term, it is challenging to tackle the magnetic term of  \eqref{GP-2}.  For this reason, we need to construct in \eqref{2:13} a Pohozaev identity, corresponding to the magnetic translation, to further estimate the global maximum point of $|w_a|$ as $a\to a_*(N)$.

%{\color{red}We remark that the idea of constructing symmetric solutions was also widely used  in studying  other equations. More precisely, the non-conical and non-planar traveling fronts with axial symmetry  for the balanced Allen--Cahn equation were  proved to exist in \cite{CGH} by  constructing a sequence of approximate waves.    Besides, Gui et al. \cite{Gui3} constructed  blow-up solutions to a mean field equation, which  possess the corresponding symmetries of the configuration. Based on these existence result, the symmetric properties of solutions can be further investigated. We refer the reader to \cite{Gui1} for the  even symmetry of monotone traveling wave solutions to the balanced Allen--Cahn equation in the entire plane, and \cite{Gui2} for the axially symmetric solutions of Allen--Cahn equation with finite Morse index. }

%Specially,  the above mentioned Pohozaev identity, instead of the usual Pohozaev identity,  is fully used to
%  %a Pohozaev identity corresponding to magnetic translation, instead of the usual Pohozaev identity, is constructed to obtain
% derive the estimates of the global maximum point of $|w_a|$.

Following Theorem \ref{thm1} on the existence of axially symmetric normalized concentrating
solutions, it is natural to wonder whether \eqref{GP-2} may admit
non-axially symmetric solutions in the class of normalized concentrating
solutions. The following theorem gives a negative answer
for this issue.

%For \(N=3\), the same conclusion is still true under the
%additional assumption that the maximum point lies on the \(x_3\)-axis, which is natural in view of Remark \ref{remA.1} below.

\begin{thm}\label{cor3}
Assume that \(0\leq\Omega<\infty\) is fixed, and let
\(w_a\in\H\) be a normalized concentrating solution of \eqref{GP-2}
concentrating at \(P_0=(0,\cdots,0,1)\) as \(a\to a_*(N)\), where $N=2,\,3$.
Suppose additionally that for $N=3,$ the maximum point \(x_a\)
of \(|w_a|\) satisfies
\begin{equation}\label{thmK:1}
x_a=(0,0,p_a),\ \ p_a\to1
	\ \ \hbox{as }\ a\to a_*(3).
\end{equation}
Then up to a constant phase (and a rotational transformation for  $N=2$), \(w_a\in\H\) is unique and axially symmetric as
\(a\to a_*(N)\).
\end{thm}

%\begin{thm}\label{cor3}
%Assume $0\leq\Omega<\infty$ is fixed, and suppose $w_a(x)\in\H$ is a normalized concentrating solution of \eqref{GP-2} concentrating at $P_0=(0,\cdots,0,1)$ as $a\to a_*(N)$, where   the maximum point \(x_a=(0,0,p_a)\) of \(|w_a|\) satisfies  for \(N=3\),
%\textcolor{red}{
%\begin{equation}\label{thmK:1}
% p_a\to1
%\quad\hbox{as }\ a\to a_*(3).
%\end{equation}}
%Then up to a constant phase (and a rotational transformation for  $N=2$), \(w_a\in\H\) must be unique and axially symmetric as
%\(a\to a_*(N)\).
%\end{thm}

\begin{rem}\label{remA.1}
Theorems~\ref{thm1} and \ref{cor3} give the existence and uniqueness of normalized concentrating solutions for \eqref{GP-2}, which satisfy \eqref{thmK:1} for  $N=3$. We also remark that for \(N=3\), the additional assumption $(\ref{thmK:1})$ is natural and necessary,
because the magnetic potential
\[
x^\Omega=\frac{\Omega}{2}(-x_2,x_1,0)
\]
is only invariant under rotations around the \(x_3\)-axis, while it is not invariant under
arbitrary rotations in \(SO(3)\). In other words, if a normalized concentrating solution  of \eqref{GP-2} with \(N=3\) has a maximum point away from the \(x_3\)-axis, then it cannot coincide with the axially symmetric solution, up to a constant phase and a rotation around the \(x_3\)-axis. Therefore, the uniqueness of Theorem \ref{cor3} with $N=3$ holds only
in the class of normalized concentrating solutions, whose maximum points lie on the \(x_3\)-axis.
\end{rem}

%\begin{rem}
%Theorem \ref{thm1} shows the existence of axially symmetric solutions
%for \eqref{GP-2} as \(a\to a_*(N)\), and Theorem \ref{cor3} further
%implies that, up to the constant phase and rotational transformation,
%normalized concentrating solutions of \eqref{GP-2} concentrating at
%\(P_0=(0,\cdots,0,1)\) must be unique and axially symmetric as
%\(a\to a_*(N)\). In the case \(N=3\), this uniqueness conclusion is
%understood under the additional assumption that the maximum point lies
%on the \(x_3\)-axis.
%\end{rem}

\begin{rem}
By deriving suitable Pohozaev identities, one can obtain that any normalized
concentrating solution of \eqref{GP-2} must concentrate at a critical
point, including the origin, of the trapping potential
\((|x|^2-1)^2\) as \(a\to a_*(N)\). Hence, if the  point $P_0=(0,\cdots,0,1)$ is replaced by any other critical point of \((|x|^2-1)^2\), then Theorem \ref{cor3}  still holds true under suitable normalization and symmetry assumptions.
\end{rem}

%\begin{rem}
%By deriving Pohozaev identities, one can obtain that any normalized concentrating solution of \eqref{GP-2} must concentrate at a critical point, including the origin, of the trapping potential $(|x|^2-1)^2$ as $a\to a_*(N)$. Moreover, the existence of Theorem \ref{thm1} and the local uniqueness of Theorem \ref{cor3} can actually be extended to the case where the normalized concentrating solution of \eqref{GP-2} concentrates at any critical point of $(|x|^2-1)^2$  as $a\to a_*(N)$. Therefore, if the particular point $P_0$ of Theorems \ref{thm1} and \ref{cor3} is replaced by any critical point of $(|x|^2-1)^2$, then both Theorems \ref{thm1} and \ref{cor3} still hold true. This further yields that up to the constant phase and rotational transformation, all normalized concentrating solutions of \eqref{GP-2} concentrating at the origin must be unique and radially symmetric as $a\to a_*(N)$.
%\end{rem}
%

%More precisely, by deriving  the detailed asymptotic formula of the level curve and using the moving plane method, Gui et al. \cite{Gui1} proved the even symmetry of monotone traveling wave solutions to the balanced Allen--Cahn equation in the entire plane. Besides, the axially symmetric solutions of Allen-Cahn
%equation with finite Morse index were also studied  in \cite{Gui2} by exploring the decaying properties of the curvature.
Since the proof of Theorem \ref{cor3} is involved, for the reader's convenience we next sketch its proof strategy by the following four steps.

\vspace{0.1cm}

{\noindent \em   Strategy of Proving Theorem \ref{cor3}.} Suppose $w_a(x)\in\H$ is a normalized concentrating solution of \eqref{GP-2} concentrating at $P_0=(0,\cdots,0,1)$ as $a\to a_*(N)$.
%The  proof of Theorem \ref{cor3} can be sketched by the following four steps.

In the first step, we shall derive the refined spike profiles of $w_a(x)$ in terms of $\mu_a<0$ as $a\to a_*(N)$, where $\mu_a\to -\infty$ as $a\to a_*(N)$ in view of Definition \ref{def2}.
%More precisely, the purpose of
%this step is to prove that $u_a(x)$ defined by \eqref{A.6} and \eqref{A.5} satisfies
% \begin{equation}\label{A3:8}
%	\begin{aligned}
%	u_{a}(x) =Q(x)+\vaa^4\psi_{1}(x)+\vaa^5\psi_{2}(x)+\vaa^6\psi_{3}(x)+o(\vaa^6)\ \ \hbox{in}\ \ \R^N\ \ \hbox{as}\ \ a\to a_*(N),
%	\end{aligned}
%	\end{equation}
%where $\psi_{j}(x)\in C^2(\R^N,\C)\cap L^\infty(\R^N,\C)$ are uniquely  defined by \eqref{3:9}.
We first establish in Lemma \ref{lem3.1} the more refined estimates of $v_a(x)$ defined in \eqref{sigpek}, based on which the refined estimates for the maximum point $x_a$ of $|w_a|$ are further obtained  by constructing a Pohozaev identity involved with the magnetic Laplacian operator. The refined spike profiles of $w_a(x)$ as $a\to a_*(N)$ are finally established in Lemma \ref{lem3.3}.

 As the second step, we shall follow the first step to establish the refined spike profiles of $w_a(x)$ in terms of $\ap>0$ as $a\to a_*(N)$, where $\ap$ is defined by
 \begin{equation*}
  \ap:=\left\{
\begin{aligned}
&\Big(\f{a^*-a}{B_{1}}\Big)^{\f{1}{4}},&\ \ \hbox{if}\ \ N=2;\\
&\f{a}{a^*},&\ \ \hbox{if}\ \ N=3;
\end{aligned}
\right.
\end{equation*}
and  $B_1>0$ is as in Lemma \ref{Nlem3.5}.
By constructing another type of Pohozaev identities, we shall first derive the refined expansions of  $\mu_a<0$ in terms of $\ap$. It deserves to comment that we shall follow the equation of $v_a(x)$ to investigate the more refined  estimates of $v_a(x)$, whose procedure improves largely the  process of \cite[Lemma 4.1]{Guo}. Applying the first step and the refined expansions of
$\mu_a<0$ in terms of $\ap$, the refined spike profiles of $w_a(x)$ in terms of $\ap$ are finally obtained in Lemmas \ref{lem3.6} and \ref{lem3.5}.

In the third step, we shall employ the results of the second step to prove the local
uniqueness of Propositions \ref{thm2} and \ref{prop:N3-axial-maximum-unique}.
%More precisely, for
%\(N=2\), the uniqueness holds up to the constant phase and rotational
%transformation. For \(N=3\), the uniqueness is proved under the
%additional assumption that the maximum point lies on the \(x_3\)-axis.
The main idea of our argument is to  construct various Pohozaev identities and apply the non-degeneracy of the limiting equation \eqref{scal equ}. In spite of this fact, comparing with \cite{CLL,LPWY,GZZ,Grossi}, one cannot obtain  the useful information from  Pohozaev identities corresponding to the tangent direction of the sphere $|x| = 1$, due to the degeneracy of the potential $(|x|^2-1)^2$.  In order to overcome this difficulty, we shall set up a new
transformation \eqref{4:1} for \(N=2\). When \(N=3\), the
argument is modified by using the additional
assumption $(\ref{thmK:1})$ and   comparing two maximum points in the \(x_3\)-direction.

As the fourth step, Theorem \ref{cor3} is finally proved by combining the existence result of Theorem
\ref{thm1} and the local uniqueness result of the third step. We refer to Sections 3 and 4 for the detailed proof of Theorem \ref{cor3}.

%As the last step, we are ready to complete the proof of Theorem \ref{cor3}. Recall from Theorem \ref{thm1} that there exists at least one normalized concentrating solution of \eqref{GP-2} concentrating  at the critical point $P_0=(0,\cdots,0,1)$ as $a\to a_*(N)$, which is also axially symmetric. Following the local uniqueness of the third step, one can then obtain that Theorem \ref{cor3} holds true.  For the detailed proof of Theorem \ref{cor3}, we refer the reader to Sections 3 and 4.

We comment that the axially symmetric solutions of Allen--Cahn equations were recently investigated in \cite{Gui1,Gui2} and the references therein.
As a byproduct of Theorem \ref{cor3}, we are able to establish the following nonexistence of vortices for normalized concentrating solutions for the case \(N=3\).

\begin{prop}\label{cor4}
Assume \(N=3\) and \(0\leq\Omega<\infty\) is fixed. Suppose
\(w_a(x)\in\H\) is a normalized concentrating solution of \eqref{GP-2}
concentrating at \(P_0=(0,0,1)\) as \(a\searrow0\), where  the maximum point \(x_a\) of \(|w_a|\) satisfies
$(\ref{thmK:1})$.
Then for all sufficiently small \(a>0\), up to a constant phase,
\(w_a\) is real-valued and free of vortices.
\end{prop}

The proof of Proposition~\ref{cor4} relies strongly on the uniqueness
result of Theorem~\ref{cor3}, which also gives that, up to a
constant phase, the unique normalized concentrating solution \(w_a\) of \eqref{GP-2} is axially symmetric for all sufficiently small
\(a>0\). Hence, the rotational term  of \eqref{GP-2} vanishes, and the vortex-free property
follows from a standard linearized argument around \(Q\), if $a>0$ is sufficiently small.
To the best of our knowledge, the existing results (such as \cite{Af2,GLY}) on the nonexistence of vortices in the whole space are mainly restricted to the special case where the trapping potential must be radially symmetric.
Proposition \ref{cor4} successfully proves that normalized concentrating solutions of \eqref{GP-2} satisfying the natural assumption $(\ref{thmK:1})$ are free of vortices in \(\R^3\), even if  the trapping potential $V(x)$ is non-radially symmetric. We also remark that for $N=2$, it remains open whether the nonexistence of vortices for
Proposition \ref{cor4} holds true for normalized concentrating
solutions concentrating on the unit circle, for instance at
\(P_0=(0,1)\).
We leave it to the interested reader.

\medskip
This paper is organized as follows. In Section 2, we shall use the reduction method to prove Theorem \ref{thm1} on the existence and axial symmetry of normalized concentrating solutions for \eqref{GP-2} as $a\to a_*(N)$. Section 3 is then devoted to the refined spike profiles of the normalized concentrating solution $w_a(x)$ for \eqref{GP-2} as $a\to a_*(N)$. Following the previous sections, in Section 4 we shall complete the proofs of Theorem \ref{cor3} and Proposition \ref{cor4}. In the Appendix, we shall prove two lemmas used in this paper.

\section{Existence of Axially Symmetric Solutions}
This section is devoted to the proof of Theorem \ref{thm1} on  the existence of  axially symmetric   solutions  for  the equation \eqref{GP-2}, which concentrate at $P_0=(0,\cdots,0,1)$ as $a\to a_*(N)$. Towards this purpose, we shall first use the reduction method to derive the existence of concentrating solutions, which are axially symmetric, for a GP equation without the $L^2$-constraint condition. We then complete in Subsection 2.2 the proof of Theorem \ref{thm1}.

\subsection{Analysis of a  GP equation without constraint}
In this subsection, we mainly consider the following equation
\begin{equation}\label{GP-4}
(-\Delta_{\va^2x^{\Omega}}) u+\Big[1+\va^2\big(|\va x+x_\va|^2-1\big)^2\Big]u-|u|^2u=0\ \ \hbox{in}\, \ \R^N,
\end{equation}
where $0\leq \Om<\infty$ is fixed, $\va>0$ is small, $x_\va\in\R^N$  satisfies
\begin{equation}\label{2.1}
x_\va=(0,\cdots,0,p_\va)\to P_0=(0,\cdots,0,1)\,\ \hbox{as}\,\ \va\to 0,
\end{equation}
and
\begin{equation}\label{magn Lap2}
-\Delta_{\va^2x^{\Omega}}:=-(\nabla-i\va^2x^{\Omega})^2=-\Delta+2i\va^2x^{\Omega}\cdot\nabla+\va^4|x^{\Omega}|^2.
\end{equation}
Note that  after rescaling, \eqref{GP-4} is similar to \eqref{GP-2} without the $L^2$-constraint  condition.

We shall prove that for $\va>0$ small enough, \eqref{GP-4} has a  solution $u_\va(x)\in H_\va$ satisfying
\begin{equation}\label{2.82}
u_\va(x)=Q(x)+v_\va(x),\ \ v_\va(x)\in F_\va,
\end{equation}
where $v_\va(x)$ satisfies
\begin{equation}\label{02.14}
\|v_\va\|_\va\to 0\ \ \hbox{as\ \  $\va\to 0$}.
\end{equation}
%Following  Lemma \ref{lemA.1} and De Giorgi-Nash-Moser theory \cite[Theorem 4.1]{HL}, one can also derive that $\|v_\va\|_{L^\infty(\R^N)}=o(1)$ as $\va\to0$.
%and $\hat{v}_\va(x)$ satisfy
%\begin{equation}\label{lower 2}
%\inte \l\{{\va}|\nabla \hat{v}_\va(x)-ix^\Om\hat{v}_\va|^2+\l(\f{1}{\va}+\va V_\Om(x)\r)|\hat{v}_\va(x)|^2\r\}\to 0\quad\hbox{as\,\  $\va\to 0$}.
%\end{equation}
Here  the Hilbert space $H_\va$ is  defined by
\begin{equation}\label{2.83}
H_\va:=\l\{u:\,\  u(x)=u(r^\prime,x_N)\,\ \hbox{for $r^\prime=\sqrt{x_1^2+\cdots+x_{N-1}^2}$,}\,\ \|u\|^2_\va<\infty\r\},
\end{equation}
and the norm $\|u\|_\va$
 is induced by the inner product
\begin{equation}\label{2.7}
\begin{split}
 \langle u,v\rangle_\va:=\operatorname{Re}\inte\Big\{&(\nabla u-i\va^2x^\Om u)(\ovl{\nabla v-i\va^2x^\Om v})\\
&+\Big[1+\va^2\big(|\va x+x_\va|^2-1\big)^2\Big]u\ovl{v}\Big\}dx,\quad u,\ v\in H_{\va}.
\end{split}
\end{equation}
The subspace $F_\va\subset H_\va$ is defined as follows:
\begin{equation}\label{AA2.10}
F_\va=\Big\{u\in H_\va :\ \ \langle u,iQ\rangle_\va=0\,\ \hbox{and}\,\ \big\langle u,\frac{\partial Q}{\partial x_N}\big\rangle_\va=0\Big\}.
\end{equation}

%We next show that there exists $v_{\va}(x)\in F_\va$ such that $u_\va(x)$ defined in \eqref{u eps} solves the equation \eqref{GP-4}.
%has a concentration solution with the form
%\begin{equation}\label{2.2}
%u_{\va,x_\va}(x):=Q(x)+v_{\va,x_\va}(x),
%\end{equation}
%where $v_{\va,x_\va}(x)\in F_\va$ satisfies
%\begin{equation}\label{2.3}
%\inte (|\nabla v_{\va,x_\va}(x)|^2+|v_{\va,x_\va}(x)|^2)=o(\va^2)\quad\hbox{as\,\  $\va\to 0$}.
%\end{equation}

Note that \eqref{GP-4} admits a solution $u_\va(x)\in H_\va$ satisfying  \eqref{2.82} and (\ref{02.14}), if and only if $v_{\va}(x)$ satisfies
\begin{equation}\label{NL-1}
\begin{split}
&\quad L_\va v_\va=l_\va(x)+R_\va(v_{\va})\quad \hbox{in}\ \ \R^N,
\end{split}
\end{equation}
where $l_\va(x)\in H_\va$ and $R_\va(v_{\va})\in H_\va$ satisfy
\begin{equation}\label{2.4}
\langle l_\va(x),\phi\rangle_\va=\operatorname{Re}\Big\{\inte -\Big[\va^4|x^\Om|^2+\va^2\big(|\va x+x_\va|^2-1\big)^2\Big]Q\ovl{\phi}\Big\},\ \ \forall\phi\in H_\va,
\end{equation}
\begin{equation}\label{2.5}
\langle R_\va(v_{\va}),\phi\rangle_\va=\operatorname{Re}\Big\{\inte \Big[(2\operatorname{Re}(v_{\va})v_{\va}+|v_{\va}|^2)Q+|v_{\va}|^2v_{\va}\Big]\ovl{\phi}\Big\},\ \ \forall\phi\in H_\va,
\end{equation}
and $L_\va$ is a bounded linear operator in $H_\va$, defined by
%and
%the operator $L_\va$ is defined by the following bi-linear functional in $H_\va$:
%\begin{equation}\label{2.6N}
%\begin{split}
%&<L_\va\psi,\phi>_{\va}\\
%:=&Re\inte\Big\{(\nabla\psi-i\va^2x^\Om \psi)\overline{(\nabla\phi-i\va^2x^\Om \phi)}+\Big[1+\va^2\big(|\va x+x_\va|^2-1\big)^2\Big]\psi\ovl{\phi}\\
%&\quad-Q^2\psi\ovl{\phi}-2Re(Q\psi)Q\ovl{\phi}\Big\}dx,\quad \forall\psi,\ \phi\in H_\va,
%\end{split}
%\end{equation}
%
%We define the following bi-linear functional
\begin{equation}
\begin{split}
&\langle L_\va\psi,\phi\rangle_{\va}\\
=&\operatorname{Re}\inte\Big\{(\nabla\psi-i\va^2x^\Om \psi)\overline{(\nabla\phi-i\va^2x^\Om \phi)}+\Big[1+\va^2\big(|\va x+x_\va|^2-1\big)^2\Big]\psi\ovl{\phi}\\ &\quad-Q^2\psi\ovl{\phi}-2\operatorname{Re}(Q\psi)Q\ovl{\phi}\Big\}dx,\ \ \forall\psi,\ \phi\in H_\va.\label{2.6}
\end{split}
\end{equation}
%Since the above bi-linear functional \eqref{2.6} is bounded in $H_\va$, it is obvious that $L_\va$ is a bounded linear operator in $H_\va$.
%\begin{lem}\label{lem2.3}
%$L_\va$ is a bounded linear operator from $\hat{H}_\va$ onto $\hat{H}_\va$.
%\end{lem}
%Next we prove $T_\va L_\va$ is bijective from $F_\va$ to $F_\va$. Following \eqref{2.8}, we obtain that $T_\va L_\va$ is injective in $F_\va$. We next prove that $T_\va L_\va F_\va=F_\va$.
%By a density argument, it suffices to prove $T_\va L_\va$ is bijective from $C_c^\infty(\R^3,\C)\cap F_\va$ to $C_c^\infty(\R^3,\C)\cap F_\va$.
%For any $u\in F_\va$, we first prove that $T_\va L_\va u\in F_\va$. By the definition of $F_\va$, we only need to check that
%$(T_\va L_\va u)(x_1,x_2,x_3)=(T_\va L_\va u)(r\prime,x_3)$ for $(r^\prime)^2=x_1^2+x_2^2.$
%For any $\phi\in \hat{H}_\va$, we denote the following rotation operator
%\begin{equation}\label{2.92}
%R(\phi(x)):=\phi(R(x_1,x_2),x_3),
%\end{equation}
%where $R(x_1,x_2)$ denotes a rotation about the $x_1$ and $x_2$ compoment.
Define the projection $T_\va$ from $H_\va$ to $F_\va$ as follows:
\begin{equation*}%\label{2.9}
T_\va u=u-c_0(iQ)-c_N\frac{\partial Q}{\partial x_N},\ \ u\in H_\va,
\end{equation*}
where $c_0\in\R$ and $c_N\in\R$ are chosen such that $$\langle T_\va u,\,iQ\rangle_\va=0\ \ \hbox{and}\ \ \langle T_\va u,\,\frac{\partial Q}{\partial x_N}\rangle_\va=0.$$ We then have the following estimate.

\begin{lem}\label{lem2.2}
Assume $\va>0$ is small enough, then the bounded linear operator $T_\va L_\va$ satisfies
\begin{equation}\label{2.8}
\|T_\va L_\va u\|_\va\geq\rho\|u\|_\va\quad\hbox{for all}\,\ u\in F_\va,
\end{equation}
where the constant $\rho>0$ is independent of $\va>0$, and $F_\va$ is as in \eqref{AA2.10}.
Moreover, $T_\va L_\va$ is bijective from $F_\va$  to $F_\va$.
\end{lem}

{\noindent \bf Proof.} By contradiction, suppose \eqref{2.8} is false. We then get that up to a subsequence if necessary, $w_\va\in F_\va$ satisfies $\|w_\va\|_\va= 1$ and
\begin{equation}\label{2.85}
\|T_\va L_\va w_\va\|_\va\to 0\quad\hbox{as}\ \ \va\to 0.
\end{equation}
This further yields that for any  $\phi\in {F}_\va$,
\begin{equation}\label{2.87}
\begin{split}
\langle L_\va w_\va,\phi\rangle_\va=\langle T_\va L_\va w_\va,\phi\rangle_\va\to 0\quad\hbox{as}\ \ \va\to 0.
\end{split}
\end{equation}

Since $\|w_\va\|_\va=1$, we deduce from Lemma \ref{lemA.1}  that
\[
\inte (|\nabla w_\va|^2+|w_\va|^2)\leq C\|w_\va\|^2_\va\leq C,
\]
where the constant $C>0$ is independent of $\va>0$.
Therefore,  there exists a function $w_0\in H^1(\R^N,\mathbb C)$ such that up to a subsequence if necessary,
\begin{equation}\label{02-1}
w_\va\rightharpoonup w_0 \quad\hbox{weakly in}\,\ H^1(\R^N,\mathbb C)\ \ \hbox{as}\ \ \va\to0,
\end{equation}
and
\begin{equation}\label{02-2}
w_\va\to w_0 \quad\hbox{strongly in}\,\ L^2_{loc}(\R^N,\mathbb C)\ \ \hbox{as}\ \ \va\to0.
\end{equation}
We then derive from \eqref{2.6} and \eqref{2.87} that for any $\phi\in C^\infty_c(\R^N,\mathbb C)\cap F_\va$, \begin{equation*}%\label{2.88}
\operatorname{Re}\inte \Big\{\nabla w_0\nabla \ovl{\phi}+w_0\ovl{\phi}-Q^2w_0\ovl{\phi}-2\operatorname{Re}(Qw_0)Q\ovl{\phi}\Big\}dx=0.
%\quad\hbox{for any $\phi\in C^\infty_c(\R^N,\C)\cap {F}_\va$.}
\end{equation*}
Moreover, we deduce from \eqref{nond} that the above identity also holds for $\phi=iQ$ and $\phi=\f{\p Q}{\p x_j}$, where $j=1,\cdots,N$. This implies that
\begin{equation*}%\label{2.89}
\operatorname{Re}\inte \Big\{\nabla w_0\nabla \ovl{\phi}+w_0\ovl{\phi}-Q^2w_0\ovl{\phi}-2\operatorname{Re}(Qw_0)Q\ovl{\phi}\Big\}dx=0,\ \ \forall\phi\in C^\infty_c(\R^N,\C)\cap H_\va.
\end{equation*}
% Following Lemma \ref{lemA.1}, we also have $\phi\in H^1(\R^N,\C)$.
Therefore, we obtain from \eqref{nond} that there exist some constants $c_j\in\R$, where $j=0,1,\cdots,N$, such that
\begin{equation*}%\label{2.90}
w_0=c_0(iQ)+\sum_{j=1}^Nc_j\f{\p Q}{\p x_j}\ \ \hbox{in}\ \ \R^N.
\end{equation*}

Next, we prove that $c_j=0$ for $j=0,1,\cdots,N$, and thus $w_0(x)\equiv 0$ in $\R^N$. Indeed, by the definition of $H_\va$, we obtain that $w_{\va}$  is even for $x_1,\cdots,x_{N-1}$.   It then follows from  \eqref{02-1}  that
\[
\operatorname{Re}\Big(\inte  w_0\f{\p Q}{\p x_j}\Big)=0,\ \ j=1,\cdots,N-1,
\]
which implies that $c_1=\cdots=c_{N-1}=0$.
Moreover, since $w_\va\in F_\va$, we  obtain from \eqref{02-1} and \eqref{02-2} that
\[
\operatorname{Re}\inte \Big\{\nabla w_0\cdot\nabla (iQ)+w_0(iQ)\Big\}dx=0,
\]
and
\[
\operatorname{Re}\inte \Big\{\nabla w_0\cdot\nabla \Big(\f{\p Q}{\p x_N}\Big)+w_0\Big(\f{\p Q}{\p x_N}\Big)\Big\}dx=0,
\]
which yield that $c_0=0$ and $c_N=0$.
We hence conclude from above that $w_0(x)\equiv0$ in $\R^N$, and
\begin{equation}\label{3.91}
w_\va\rightharpoonup 0 \ \ \hbox{weakly in}\,\ H^1(\R^N,\mathbb C)\ \ \hbox{and}\ \ w_\va\to 0 \ \ \hbox{strongly in}\,\ L^2_{loc}(\R^N,\mathbb C)\, \ \hbox{as}\ \ \va\to0.
\end{equation}

On the other hand, taking $R>0$ large enough, we derive from \eqref{1:3} and \eqref{3.91} that
\[
\begin{aligned}
&\quad \inte \Big(Q^2|w_\va|^2+2\big(Q\operatorname{Re}(w_\va)\big)^2\Big)\\
&=\int_{B_R(0)}\l(Q^2|w_\va|^2+2\big(Q\operatorname{Re}(w_\va)\big)^2\r)+\int_{\R^N\backslash B_R(0)}\l(Q^2|w_\va|^2+2\big(Q\operatorname{Re}(w_\va)\big)^2\r) \\
&\leq  3\|Q\|_{L^\infty(\R^N)}^2 \int_{B_R(0)}|w_\va |^2 +\f14 \int_{\R^N\backslash B_R(0)}|w_\va|^2\leq \f 12, \ \ \hbox{if}\ \ \va>0\ \ \hbox{is small}.
\end{aligned}
\]
Therefore, since $w_\va\in F_\va$ satisfies $\|w_{\va}\|_{\va}=1$, setting $\phi=w_\va$ into \eqref{2.87} yields  that
\[
\begin{aligned}
0&=\lim_{\va\to 0}\inte\Big\{ \Big|(\nabla-i\va^2 x^\Om)w_\va\Big|^2+\Big[1+\va^2\big(|\va x+x_\va|^2-1\big)^2\Big]|w_\va|^2\\
&\qquad \quad \ \qquad -Q^2|w_\va|^2-2\big(Q\operatorname{Re}(w_\va)\big)^2\Big\}\\
&=\lim_{\va\to 0}\|w_\va\|^2_\va-\lim_{\va\to 0}\inte \Big(Q^2|w_\va|^2+2\big(Q\operatorname{Re}(w_\va)\big)^2\Big)\geq \f 12,
\end{aligned}
\]
a contradiction. We thus conclude that  \eqref{2.8} holds true.
%
%where
%\[
%I_{1\va}:=\int_{B_R(0)} |(\nabla-i\va^2 x^\Om)w_\va|^2+(1+\va^2V_\Om(\va x+x_\va))|w_\va|^2-Q^2|w_\va|^2-2(QRe(w_\va))^2
%\]
%and
%\[
%I_{2\va}:=\int_{\R^3\backslash B_R(0)} |(\nabla-i\va^2 x^\Om)w_\va|^2+(1+\va^2V_\Om(\va x+x_\va))|w_\va|^2-Q^2|w_\va|^2-2(QRe(w_\va))^2
%\]
%with $R>0$ is a sufficiently large constant. The convergence in \eqref{3.91} implies that
%$$\lim_{\va\to 0}I_{1\va}=0.$$
%Thus we have
%\begin{equation}\label{2.91}
%\begin{split}
%0&=\lim_{\va\to 0}I_{2\va}\\
% &\geq \lim_{\va\to 0}\int_{\R^3\backslash B_R(0)} |(\nabla-i\va^2 x^\Om)w_\va|^2+(1+\va^2V_\Om(\va x+x_\va))|w_\va|^2-\f 12|w_\va|^2\\
% &\geq \lim_{\va\to 0}\f 12\int_{\R^3\backslash B_R(0)} |(\nabla-i\va^2 x^\Om)w_\va|^2+(1+\va^2V_\Om(\va x+x_\va))|w_\va|^2.
%\end{split}
%\end{equation}
%Combining \eqref{3.91} and \eqref{2.91} yields that
%\[
%\|w_\va\|_\va\to 0\quad\hbox{as}\,\ \va\to 0,
%\]
%which however contradicts $\|w_\va\|_\va=1$. This proves \eqref{2.8}.

Once \eqref{2.8} is true, it is standard to prove that $T_{\va}L_\va$ is  bijective  from $F_\va$  to  $F_\va$. This completes the proof of Lemma \ref{lem2.2}. \qed

Applying Lemma \ref{lem2.2}, one can establish the following existence.

\begin{lem}\label{lem2.4}
	Suppose $x_\va=(0,\cdots,0,p_\va)\in\R^N$ satisfies $x_\va\to P_0=(0,\cdots,0,1)$ as $\va\to 0$. Then for sufficiently small $\va>0$, there exists a unique $v_\va(x) \in F_\va$ defined by \eqref{AA2.10}, which  satisfies $\|v_\va\|_\va\to 0$ as $\va\to 0$, such that
	$u_{\va }(x)=Q(x)+v_{\va}(x)$ solves
	\begin{equation}\label{eq redu1}
	\begin{split}
	&\left(-\Delta_{\va^2x^{\Omega}}\right) u_{\va}+\Big[1+\va^2\big(|\va x+x_\va|^2-1\big)^2\Big]u_{\va}-|u_{\va}|^2u_{\va}\\
	=&a_{\va,x_\va}(iQ)+b_{\va,x_\va}\frac{\partial Q}{\partial x_N}\quad \hbox{in}\ \ \R^N,
	\end{split}
	\end{equation}
	where $a_{\va,x_\va}$, $b_{\va,x_\va}\in\R$  are some constants   depending only on $\va>0$ and $x_\va$.
\end{lem}

Since the proof of Lemma \ref{lem2.4} is similar to that of \cite[Proposition 2.2.6]{CPY}, we omit the details for simplicity.
To obtain  the existence of solutions for \eqref{GP-4}, we also need to establish the following lemma.

 %Following above lemmas, we are now ready to establish the following existence of solutions for \eqref{GP-4} satisfying  \eqref{2.82} and \eqref{02.14}.

\begin{lem}\label{lem2.5}
	Suppose $\va> 0$ is small enough, then there exists $x_\va=(0,\cdots,0,p_\va)\in\R^N$ satisfying $x_\va\to P_0=(0,\cdots,0,1)$ as $\va\to0$,  such that the constants of \eqref{eq redu1} satisfy $a_{\va,x_\va}=b_{\va,x_\va}=0$ for sufficiently small $\va>0$.
%Moreover, there exists a solution $u_\va(x)=Q(x)+v_\va(x)$ of \eqref{GP-4}, where $v_\va\in F_\va$ satisfies $\|v_\va\|_\va\to0$ as $\va\to 0$.	
\end{lem}

\noindent{\bf Proof.}
Multiplying \eqref{eq redu1} by $\ovl{iu}_{\va}$, integrating over $\R^N$ and then taking its real part, we have
\begin{equation}\label{2.93}
a_{\va,x_{\va}}\operatorname{Re} \Big(\inte (iQ)\ovl{iu}_{\va}\Big)+b_{\va,x_{\va}} \operatorname{Re}\Big(\inte \f{\p Q}{\p x_N}\ovl{iu}_\va\Big)=0.
\end{equation}	
On the other hand, multiplying \eqref{eq redu1} by  $\frac{\partial \ovl{u}_{\va}}{\partial x_N}+\f{-i\va^2\Om\delta_{2N}x_{1}\ovl u_{\va}}{2}$, integrating over $\R^N$ and then taking its real part, we have
\begin{equation}\label{2.10}
\begin{aligned}
&\operatorname{Re}\Big\{\inte\Big((-\Delta_{\va^2x^\Om})u_{\va}+\Big[1+\va^2\big(|\va x+x_\va|^2-1\big)^2\Big]u_{\va}-|u_{\va}|^2u_{\va}\Big)\\
&\qquad\qquad\Big(\frac{\partial \ovl{u}_{\va}}{\partial x_N}+\f{-i\va^2\Om\delta_{2N}x_{1}\ovl u_{\va}}{2}\Big )\Big\}\\
=&\operatorname{Re}\Big\{\inte \Big(a_{\va,x_\va}(iQ)+b_{\va,x_\va}\frac{\partial Q}{\partial x_N}\Big)\Big(\frac{\partial \ovl{u}_{\va}}{\partial x_N}+\f{-i\va^2\Om\delta_{2N}x_{1}\ovl u_{\va}}{2}\Big )\Big\},
\end{aligned}
\end{equation}
where $\delta_{22}=1$ and $\delta_{2N}=0$ for $N\neq2$.
We next calculate all terms of \eqref{2.10}.

We first claim that there exists a sequence $\{R_i\}$  such that
\begin{equation}\label{A2-1}
\begin{aligned}
  \lim\limits_{R_i\to\infty}\int_{\p B_{R_i}(0)}\Big(\l|(\nabla-i\va^2 x^\Om)u_{\va}\r|^2+|u_{\va}|^2+|u_{\va}|^4\Big)|\nu_{N}|dS=0,
\end{aligned}
\end{equation}
where $\nu_{N}$ denotes the  outward unit normal vector of $\p B_{R_i}(0)$ along the $x_N$-axis.
On the contrary, assume that there exists a constant $C>0$ such that for any sufficiently large $R>0$,
\begin{equation*}
  \int_{\p B_{R}(0)}\Big(\l|(\nabla-i\va^2 x^\Om)u_{\va}\r|^2+|u_{\va}|^2+|u_{\va}|^4 \Big)|\nu_{N}|dS\geq C.
\end{equation*}
We then deduce that for some $R_1>0$,
\begin{equation*}
\begin{aligned}
  &\int_0^{+\infty}\int_{\p B_{R}(0)}\Big(\l|(\nabla-i\va^2 x^\Om)u_{\va}\r|^2+|u_{\va}|^2+|u_{\va}|^4 \Big)dSdR\\
  =&\inte \Big(\l|(\nabla-i\va^2 x^\Om)u_{\va}\r|^2+|u_{\va}|^2+|u_{\va}|^4 \Big)dx\geq \int_{R_1}^{+\infty}CdR=+\infty,
\end{aligned}
\end{equation*}
which however contradicts the fact that
\begin{equation*}
  \inte \Big(\l|(\nabla-i\va^2 x^\Om)u_{\va}\r|^2+|u_{\va}|^2+|u_{\va}|^4 \Big)dx<\infty.
\end{equation*}
We then conclude that the claim \eqref{A2-1} holds true.

By direct calculations, we obtain from \eqref{magn Lap2} and \eqref{A2-1} that
\begin{equation*}
\begin{aligned}
&\operatorname{Re} \l\{\inte (-\Delta_{\va^2x^\Om})u_{\va}\Big(\frac{\partial \ovl{u}_{\va}}{\partial x_N}+\f{-i\va^2\Om\delta_{2N}x_{1}\ovl u_{\va}}{2}\Big)\r\}\\
%   =&\operatorname{Re}\l\{\inte -(\nabla-i\va^2 x^\Om)^2u_{\va}\l(\frac{\partial \ovl{u}_{\va}}{\partial x_N}+\f{-i\va^2\Om\delta_{2N}x_{1}\ovl u_{\va}}{2}\r)\r\}\\
%   =&Re\l\{\inte (\nabla-i\va^2 x^\Om)u_{\va}  (\nabla+i\va^2 x^\Om)\l(\frac{\partial \ovl{u}_{\va}}{\partial x_N}+\f{-i\va^2\Om\delta_{2N}x_{1}\ovl u_{\va}}{2}\r)\r\}\\
% %=&Re\inte (\nabla-i\va^2 x^\Om)u_{\va}\Big[\f{\p}{\p x_{N}}(\nabla+i\va^2 x^\Om)\ovl u_{\va}-i\va^2\f{\p x^\Om}{\p x_{N}}\ovl u_{\va}\\
% %&\quad+(\nabla+i\va^2 x^\Om)\f{-i\va^2\Om\delta_{2N}x_{1}\ovl u_{\va}}{2}\Big]\\
% %=&Re\inte \l(\fh\f{\p}{\p x_{N}}\l|(\nabla-i\va^2 x^\Om)u_{\va}\r|^2-\f{i\va^2\Om\delta_{2N}x_{1}}{2}\l|(\nabla-i\va^2x^\Om)u_{\va}\r|^2\r)\\
=& \fh\inte \f{\p}{\p x_{N}}\l|(\nabla-i\va^2 x^\Om)u_{\va}\r|^2\\
=&\fh\lim\limits_{R_i\to\infty}\int_{B_{R_i}(0)}\f{\p}{\p x_{N}}\l|(\nabla-i\va^2 x^\Om)u_{\va}\r|^2=0.
\end{aligned}
\end{equation*}
Moreover, we have
\begin{equation*}
\begin{aligned}
&\operatorname{Re} \l\{\inte \l(u_{\va}-|u_{\va}|^2u_{\va}\r)\Big(\frac{\partial \ovl{u}_{\va}}{\partial x_N}+\f{-i\va^2\Om\delta_{2N}x_{1}\ovl u_{\va}}{2}\Big )\r\}\\
% =&\fh\inte \f{\p |u_{\va}|^2}{\p x_{N}}-\frac14\inte\frac{\partial |{u_{\va}}|^4}{\partial x_N}\\
=&\fh\lim\limits_{R_i\to\infty}\int_{\p B_{R_i}(0)}|u_{\va}|^2 \nu_{N}dS-\f{1}{4}\lim\limits_{R_i\to\infty}\int_{\p B_{R_i}(0)}|u_{\va}|^4 \nu_{N}dS=0,
\end{aligned}
\end{equation*}
and
\begin{equation*}
\begin{aligned}
&\va^2\operatorname{Re}\l\{\inte \big(|\va x+x_\va|^2-1\big)^2u_{\va}\Big(\frac{\partial \ovl{u}_{\va}}{\partial x_N}+\f{-i\va^2\Om\delta_{2N}x_{1}\ovl u_{\va}}{2}\Big )\r\}\\
%=&\f{\va^2}{2}\inte V_{\Om}(\va x+x_{\va})\frac{\partial |u_{\va}|^2}{\partial x_N}\\
=&-\f{\va^2}{2}\inte\f{\p \big(|\va x+x_\va|^2-1\big)^2}{\p x_{N}}  |u_{\va}|^2.\\
\end{aligned}
\end{equation*}
In addition, since $u_{\va}(x)=Q(x)+v_{\va}(x)$, we derive  that
\begin{equation*}
\begin{aligned}
&a_{\va,x_\va}\operatorname{Re}\Big\{\inte (iQ) \Big(\frac{\partial \ovl{u}_{\va}}{\partial x_N}+\f{-i\va^2\Om\delta_{2N}x_{1}\ovl u_{\va}}{2}\Big )\Big\}\\
&+
b_{\va,x_\va}\operatorname{Re}\Big\{\inte \frac{\partial Q}{\partial x_N}\Big(\frac{\partial \ovl{u}_{\va}}{\partial x_N}+\f{-i\va^2\Om\delta_{2N}x_{1}\ovl u_{\va}}{2}\Big )\Big\}\\
% =&a_{\va,x_\va}\operatorname{Re}\l(\inte (iQ) \frac{\partial \ovl{u}_{\va}}{\partial x_N}\r)+
% b_{\va,x_\va}\operatorname{Re}\l(\inte \frac{\partial Q}{\partial x_N}\frac{\partial \ovl{u}_{\va}}{\partial x_N}\r)\\
% &\quad-\f{\va^2\Om\delta_{2N}}{2}\l\{a_{\va,x_{\va}}Re\l(\inte i(iQ)x_{1}\ovl{( Q+ v_{\va})}\r)+b_{\va,x_{\va}}Re\l(\inte i\f{\p Q}{\p x_{N}}x_{1}\ovl {(Q+v_{\va})}\r)\r\}\\
=&a_{\va,x_\va}\operatorname{Re}\Big(\inte (iQ) \frac{\partial \ovl{u}_{\va}}{\partial x_N}\Big)+
b_{\va,x_\va}\operatorname{Re}\Big(\inte \frac{\partial Q}{\partial x_N}\frac{\partial \ovl{u}_{\va}}{\partial x_N}\Big),
\end{aligned}
\end{equation*}
where the radial symmetry of $Q(x)$ and the axial symmetry of $v_{\va}(x)$ are also used.
We therefore obtain from above that
\begin{equation}\label{2.94}
\begin{aligned}
&-\f{\va^2}{2}\inte \frac{\partial \big(|\va x+x_\va|^2-1\big)^2}{\partial x_N}|u_{\va}|^2\\
=&a_{\va,x_\va}\operatorname{Re}\Big(\inte (iQ) \frac{\partial \ovl{u}_{\va}}{\partial x_N}\Big)+
b_{\va,x_\va}\operatorname{Re}\Big(\inte \frac{\partial Q}{\partial x_N}\frac{\partial \ovl{u}_{\va}}{\partial x_N}\Big).
\end{aligned}
\end{equation}

We now claim that there exists $x_\va$ in the neighborhood of $P_0$ such that the left hand side of \eqref{2.94} equals to $0$. Indeed, since $u_\va(x)=Q(x)+v_\va(x)$, direct calculations yield  that
\begin{equation}\label{2.11}
\begin{split}
g(p_\va):=&\inte \f{\p \big(|\va x+x_\va|^2-1\big)^2}{\p x_N}|u_{\va}|^2\\
=&\va\Big(\inte 4Q^2\Big)(p^2_\va-1)p_\va +o(\va|p_\va-1|)+o(\va^2)\quad\hbox{as}\,\ \va\to 0.
\end{split}
\end{equation}
We then derive from \eqref{2.11} that $g(t)$ is a small perturbation of
$$f(t):=\va\Big(\inte 4Q^2\Big)(t^2-1)t,$$
and $\deg (f,B_\delta(1),0)\neq 0$.  This implies  that there exists at least one $p_\va$  such that $g(p_\va)=0$, and hence the above claim is true.

Applying the above claim, we now obtain from \eqref{2.93} and \eqref{2.94} that for sufficiently small $\va>0$,
\begin{equation*}%\label{2.95}
\left\{\begin{array}{lll}
a_{\va,x_\va}\operatorname{Re}\Big( \inte (iQ)\ovl{iu}_{\va}\Big)+b_{\va,x_\va} \operatorname{Re}\Big(\inte \f{\p Q}{\p x_N}\ovl {iu}_\va\Big)=0, \\[3mm]
a_{\va,x_\va}\operatorname{Re}\Big(\inte (iQ) \frac{\partial \ovl{u}_{\va}}{\partial x_N}\Big)
+
b_{\va,x_\va}\operatorname{Re}\Big(\inte \frac{\partial Q}{\partial x_N}\frac{\partial \ovl{u}_{\va}}{\partial x_N}\Big)=0.
\end{array}\right.
\end{equation*}
Moreover, we obtain from Lemmas \ref{lem2.4} and \ref{lemA.1} that for sufficiently small $\va>0$,
\begin{equation*}%\label{2.96}
\left | \begin{matrix}
\operatorname{Re} \Big(\inte (iQ)\ovl{iu}_{\va}\Big)&\operatorname{Re}\Big(\inte \f{\p Q}{\p x_N}\ovl{iu}_\va\Big) \\
\operatorname{Re}\Big(\inte (iQ) \frac{\partial \ovl{u}_{\va}}{\partial x_N}\Big) &\operatorname{Re}\Big(\inte \frac{\partial Q}{\partial x_N}\frac{\partial \ovl{u}_{\va}}{\partial x_N}\Big) \\
\end{matrix} \right | \neq 0,
\end{equation*}
which hence implies that $a_{\va,x_\va}=b_{\va,x_\va}=0$ holds for sufficiently small $\va>0$.
This completes the proof of Lemma \ref{lem2.5}. \qed

One can conclude from Lemmas \ref{lem2.4} and \ref{lem2.5} that there
exists a solution \(u_\va(x)=Q(x)+v_\va(x)\) of \eqref{GP-4} as
\(\va\to0\), where \(v_\va\in F_\va\) defined by \eqref{AA2.10}
satisfies \(\|v_\va\|_\va\to0\) as \(\va\to0\).

\subsection{Proof of Theorem \ref{thm1}}
In this subsection, we shall complete the proof of  Theorem \ref{thm1}. We start with the following estimate of the error term $v_\va(x)\in F_\va$.

\begin{lem}\label{lem2.1}
	% Under the assumptions of Lemma \ref{lem2.5},
	Assume  $x_{\va}=(0,\cdots,0,p_{\va})\in\R^N$ satisfies $x_{\va}\to P_0 =(0,\cdots,0,1)$ as $\va\to0$, and  suppose  $u_\va(x)=Q(x)+v_\va(x)$ is a  solution of \eqref{GP-4} satisfying \eqref{02.14}, where $v_\va(x)\in F_{\va}$. Then
	the error term $v_\va(x)$ satisfies
	\begin{equation}\label{2.91}
	\|v_\va\|_\va\leq C\Big[\va^4+\va^2(|x_\va|^2-1)^2\Big] \ \ as\ \ \va\to0,
	\end{equation}
	where the constant $C>0$ is independent of $\va>0$.
\end{lem}

Since it is not difficult to prove Lemma \ref{lem2.1} by applying Lemmas \ref{lem2.2} and \ref{lemA.1}, we omit the detailed proof for simplicity.
Based on Lemma \ref{lem2.1}, one can further derive the  following more refined estimates of the point $x_\va\in\R^N$ in \eqref{GP-4}.

\begin{lem}\label{Nlem3.0}
	Under the assumptions of Lemma \ref{lem2.1}, the point $x_{\va}=(0,\cdots,0,p_{\va})\in\R^N$ satisfies
	% Assume $u_{\va}(x)$ defined by \eqref{u eps} satisfies \eqref{3.5} and \eqref{3.6}, where $w_{a}$ is a solution of the equation \eqref{GP-1}, and $x_{a}=(0,0,p_{a})$ is a global maximum point of $|w_{a}|$. Then we have
	\begin{equation}\label{2:12}
	p_{\va}-1=A\va^2+o(\va^2)\ \ \hbox{as}\ \ \va\to0,
	\end{equation}
	where  the constant  $A:=-\frac{(N+2)\inte |x|^2Q^2dx}{2Na^*}<0$.
\end{lem}

\noindent{\bf Proof.}
Similar to the proof of  \eqref{2.94}, multiplying \eqref{GP-4} by  $\frac{\partial \ovl{u}_{\va}}{\partial x_N}+\f{-i\va^2\Om\delta_{2N}x_{1}\ovl u_{\va}}{2}$, integrating over $\R^N$ and then taking its real part, we deduce %from the exponential decay \eqref{2:11}
that
\begin{equation}\label{2:13}
\inte\frac{\partial \big(|\va x+x_\va|^2-1\big)^2}{\partial x_{N}}|u_{\va}|^2dx=0,\ \ \hbox{if}\ \ \va>0\ \ \hbox{is small enough}.
\end{equation}

On the other hand,
by the comparison principle, one can actually obtain from \eqref{NL-1} that
\begin{equation}\label{2-1}
|v_\va|\leq Ce^{-\frac{2}{3}|x|}\ \ \hbox{uniformly in}\ \ \R^N\ \ \hbox{as}\ \ \va\to 0.
\end{equation}
Since $u_{\va}(x)=Q(x)+v_{\va}(x)$, applying Lemmas \ref{lem2.1} and \ref{lemA.1},  we obtain from \eqref{2-1} that
\begin{equation*}%\label{2:14}
\begin{aligned}
&\quad\inte\frac{\partial \big(|\va x+x_\va|^2-1\big)^2}{\partial x_{N}}|u_{\va}|^2dx\\
% &=\inte4\va(|\va x+x_{\va}|^2-1)(\va x_{N}+p_{\va})\big|Q+v_{\va}\big|^2dx\\
&=\inte4\va\big(\va^2|x|^2+2\va x_{N}p_{\va}+p_{\va}^2-1\big)(\va x_{N}+p_{\va})\l(Q^2+2\operatorname{Re}(Qv_{\va})+|v_{\va}|^2\r)dx\\
&=8\va(p_{\va}-1)\inte Q^2dx\big(1+o(1)\big)+4\va^3\inte |x|^2Q^2dx\\
&\quad+8\va^3\inte x_{N}^2Q^2dx
% +o(|p_{\va}-1|)
+o\big(\va^3\big)\ \ \hbox{as}\ \ \va\to0,
\end{aligned}
\end{equation*}
and hence
\begin{equation*}%\label{2:18}
\begin{aligned}
p_{\va}-1&=-\frac{4\inte |x|^2Q^2dx+8\inte x_{N}^2Q^2dx}{8\inte Q^2dx}\va^2+o(\va^2)\\
%&=-\frac{\inte |x|^2Q^2dx+2\inte x_{N}^2Q^2dx}{2\inte Q^2dx}\va^2+o(\va^2)\\
&=-\frac{(N+2)\inte |x|^2Q^2dx}{2Na^*}\va^2+o(\va^2)\ \ \hbox{as}\ \ \va\to0.
\end{aligned}
\end{equation*}
 We therefore conclude that Lemma \ref{Nlem3.0} holds true.
% This implies that\eqref{2:12} holds true, and the proof of Lemma \ref{Nlem3.0} is therefore complete.
\qed

To prove Theorem \ref{thm1}, we also need to establish the following  Pohozaev identity.

\begin{lem}\label{lem2.6} Suppose $u_\va(x)=Q(x)+v_\va(x)$ is a solution of  \eqref{GP-4} as $\va\to0$, where  $\|v_\va\|_\va\to 0$ as $\va\to 0$. Then $u_\va(x)$ satisfies the following Pohozaev identity:
	\begin{equation}\label{2.92}
	\begin{split}
	&\quad\inte \Big(-\frac{N-2}{2}|\nabla u_\va|^2+2\operatorname{Re}\Big[i\va^2 (x^\Om\cdot\nabla u_\va)(x\cdot\nabla \ovl u_\va)\Big]+\frac{N}{4}|u_\va|^4-\f{N}{2}|u_\va|^2\Big)\\
	&=\inte \Big(\f{N+2}{2}\va^4|x^\Om|^2+\f{N\va^2}{2}\big(|\va x+x_\va|^2-1\big)^2\\
	&\qquad\quad \ \ +2\va^3(|\va x+x_\va|^2-1)(\va |x|^2+x\cdot x_\va)\Big)|u_\va|^2.
	\end{split}
	\end{equation}
\end{lem}

Lemma \ref{lem2.6} can be proved by multiplying \eqref{GP-4} with $x\cdot\nabla\bar u_{\va}$, integrating over $\R^N$ and then taking its real part. We omit the detailed proof for simplicity. Employing above lemmas, we next finish the proof of Theorem \ref{thm1}.

\vskip 0.05truein

\noindent{\bf Proof of Theorem \ref{thm1}.}
Let $0\le \Omega <\infty$ be fixed.
Following Lemma \ref{lem2.5}, we obtain that
%{\color{red}there exists a sufficiently small constant $\va_0>0$ such that for any $\va\in(0,\va_0)$},
for sufficiently small $\va>0$,
there exists a solution $u_\va(x)=Q(x)+v_\va(x)$ of \eqref{GP-4}, where $v_\va\in F_\va$   satisfies \eqref{02.14}. Denote
	\begin{equation}\label{R2-1}
	\hat{u}_\va(x):=\f{1}{\va}u_\va\Big(\f{x-x_\va}{\va}\Big)e^{ix\cdot x_\va^\Om},
	\end{equation}
where $x_\va=(0,\cdots,0,p_\va)\in\R^N$ is as in \eqref{GP-4}, and $x_\va^\Om$ is defined by
\begin{equation*}
 x_\va^\Omega:=\left\{
\begin{aligned}
&\frac{\Om}{2}(-p_\va,0),&\   \hbox{if}\ \ N=2;\\
&\frac{\Om}{2}(0,0,0),&    \hbox{if}\ \ N=3.
\end{aligned}
\right.
\end{equation*}
Setting
\begin{equation}\label{2.97}
a_\va:=\inte |\hat{u}_\va|^2,\,\ w_\va:=\f{\hat{u}_\va}{\sqrt{a_\va}},\ \ \hbox{and}\ \ \mu_\va:=-\f{1}{\va^2},
\end{equation}
we then derive that $w_{\va}(x)$ satisfies the following equation
\begin{equation}\label{2.13}
\left\{
\begin{aligned}
&-\Delta_{x^\Om} w_\va+(|x|^2-1)^2w_\va-a_\va|w_\va|^2w_\va=\mu_{\va}w_\va\ \ \hbox{in}\ \ \R^N,\\
&\inte |w_\va|^2dx=1,
\end{aligned}
\right.
\end{equation}
where $-\Delta_{x^\Om}$ is defined as in \eqref{GP-2}.
This implies that when $a=a_\va$, where $\va>0$ is sufficiently small, $(w_\va,\mu_\va)$  is an axially symmetric solution of \eqref{GP-2}.

We claim that for $N=2$ (or $N=3$), there exists a small constant $\delta>0$ such that for any $a\in (a^*-\delta,a^*)$ ( or $a\in (0,\delta)$, respectively), there exists a sufficiently small constant $\vaa>0$
%$\va_a\in(0,\va_0)$
such that
\begin{equation}\label{2:24}
a=a_{\va_a}=\inte |\hat{u}_{\va_a}|^2,
%\ \ \hbox{and}\ \ \mu_{\va_a}=-\f{1}{\va_a^2}<0.
\end{equation}
which further  gives that $(w_{\vaa},\mu_{\vaa})$
defined in \eqref{2.97} is an axially symmetric solution of  \eqref{GP-2} as $a\to a_*(N)$. We next prove  the claim \eqref{2:24} for two different cases.

{\em Case 1: $N=2$.}
Note from Lemmas \ref{lem2.1} and \ref{Nlem3.0} that
\begin{equation}\label{2.98}
\|v_\va\|_\va\leq C\va^4\quad\hbox{and}\quad p_{\va}-1=A\va^2+o(\va^2)\ \ \hbox{as}\ \ \va\to0,
\end{equation}
where $A<0$ is as in \eqref{2:12}. Since $a_\va=\int_{\R^2} |\hat{u}_\va|^2=\int_{\R^2} |{u}_\va|^2$, we next use  Lemma \ref{lem2.6} and \eqref{2.98} to give the estimates of $\int_{\R^2} |{u}_\va|^2$ in terms of $\va>0$.

%Next, we use the Pohozaev identity \eqref{2.92}to prove the claim \eqref{2:24}.
%We now estimate both hand sides of  the Pohozaev identity \eqref{2.92} established  in Lemma \ref{lem2.6}.
Similar to \cite[Lemma A.2]{GLY}, we obtain from \eqref{GP-4} that
\begin{equation}\label{F2:11}
|u_{\va}(x)|\leq Ce^{-\frac{2}{3}|x|}\ \ \hbox{and}\ \ |\nabla u_{\va}(x)|\leq Ce^{-\frac{1}{2}|x|}
\ \ \text{uniformly in}\, \ \R^N\ \ \hbox{as}\ \ \va\to0.
\end{equation} Following Lemma \ref{lemA.1}%to the Pohozaev identity \eqref{2.92}
, we then derive from \eqref{2.82}, \eqref{2.98} and \eqref{F2:11}  that
\begin{equation}\label{2.99}
\operatorname{Re}\int_{\R^2} (i\va^2 x^\Om\cdot\nabla u_\va)(x\cdot\nabla \ovl u_\va)=\operatorname{Re}\int_{\R^2} (i\va^2 x^\Om\cdot\nabla v_\va)(x\cdot\nabla \ovl u_\va)=O(\va^6)\quad\hbox{as}\,\ \va\to 0.
\end{equation}
Moreover, we have
\begin{equation}\label{2.51}
\begin{split}
\quad\int_{\R^2}\Big(\f 12|u_\va|^4-|u_\va|^2\Big)&=\int_{\R^2}\f 12|Q+v_\va|^4-a_\va\\
&=\int_{\R^2}\Big(\f{1}{2}Q^4+2\operatorname{Re}(Q^3v_\va) \Big)-a_\va+O(\va^8)\ \ \hbox{as}\ \ \va\to0.
\end{split}
\end{equation}
It follows from \eqref{NL-1} and \eqref{2.98} that
\begin{equation}\label{2.521}
\begin{split}
&\quad\int_{\R^2} \operatorname{Re}(Q^3v_\va)\\
&=\int_{\R^2}\operatorname{Re}\big((-\Delta Q+Q)v_\va\big)=\int_{\R^2}\operatorname{Re}\big((-\Delta v_\va+v_\va)Q\big)\\
&=\int_{\R^2}\operatorname{Re}\Big(\big(-2i\va^2 x^\Om\cdot\nabla v_\va-\va^4|x^\Om|^2v_\va-\va^2\big(|\va x+x_\va|^2-1\big)^2v_\va+Q^2v_\va\\
&\quad+2\operatorname{Re}(Q^2v_\va)+l_\va(x)+R_\va(v_\va)\big)Q\Big)\\
&=3\int_{\R^2}\operatorname{Re}(Q^3v_\va)-\f{\va^4(\Om^2+8)}{4}\int_{\R^2}|x|^2Q^2dx+O(\va^6)\ \ \hbox{as}\ \ \va\to0,
\end{split}
\end{equation}
which yields that
\begin{equation}\label{2.52}
\begin{split}
&2\int_{\R^2} \operatorname{Re}(Q^3v_\va)
=\f{\va^4(\Om^2+8)}{4}\int_{\R^2}|x|^2Q^2dx+O(\va^6)\ \ \hbox{as}\ \ \va\to0.
\end{split}
\end{equation}
%we find that $\int_{\R^2} \operatorname{Re}(Q^3v_\va)=O(\va^6)$.
We thus derive from \eqref{N1:1} and  \eqref{2.99}--\eqref{2.52} that
\begin{equation}\label{2.53}
% \int_{\R^2}\l(\f 12|u_\va|^4-|u_\va|^2\r)
\hbox{LHS of \eqref{2.92}}=a^*-a_\va+\f{\va^4(\Om^2+8)}{4}\int_{\R^2}|x|^2Q^2dx+O(\va^6)\quad\hbox{as}\,\ \va\to 0.
\end{equation}
Similarly, we deduce from \eqref{2.98} and Lemma \ref{lemA.1} that
\begin{equation}\label{2.54}
\begin{split}
&\qquad\hbox{RHS of \eqref{2.92}}\\
&=\int_{\R^2} \Big[2\va^4|x^\Om|^2+\va^2\big(|\va x+x_\va|^2-1\big)^2+2\va^3(|\va x+x_\va|^2-1)(\va |x|^2+x\cdot x_\va)\Big]|u_\va|^2\\
% &=\int_{\R^2} \l(2\va^4|x^\Om|^2+4\va^4x_2^2+4\va^4x_2^2\r)Q^2+O(\va^6)\\
&=\f{\Om^2+8}{2}\va^4\int_{\R^2}|x|^2Q^2+O(\va^6)\quad\hbox{as}\,\ \va\to 0.
\end{split}
\end{equation}
We then conclude from \eqref{2.92},  \eqref{2.53} and \eqref{2.54} that
\begin{equation}\label{2.55}
a_\va=a^*-\f{\Om^2+8}{4}\va^4\int_{\R^2}|x|^2Q^2+O(\va^6)\quad\hbox{as}\,\ \va\to 0.
\end{equation}
% which also implies that
% \begin{equation}\label{2:25}
%   \f{1}{\va^2}\l(\f{a^*-a_{\va}}{B_{1}}\r)^\fh=1+O(\va^2)\ \ \hbox{as}\ \ \va\to0,
% \end{equation}
% where $B_{1}$ satisfies
% \begin{equation*}
%   B_{1}:=\f{\Om^2+8}{4}\int_{\R^2}|x|^2Q^2.
% \end{equation*}

Let
$\va_{0}>0$
be a small constant
such that \eqref{GP-4} has a solution $u_{\va_{0}}$ satisfying  \eqref{2.82} and \eqref{02.14}.
Setting  $\delta=a^*-a_{\va_0}$, where  $a_{\va_0}=\int_{\R^2}|\hat{u}_{\va_{0}}|^2$, we then obtain from \eqref{2.55} that for any $a\in(a^*-\delta,a^*)$, there exists a small constant $\va_1\in (0,\va_0)$ such that
\begin{equation*}
a_{\va_0}=a^*-\delta<a<a^*-\f{\Om^2+8}{2}\va_{1}^4\int_{\R^2}|x|^2Q^2<a_{\va_{1}}.
\end{equation*}
Since $a_{\va}$ depends on $\va>0$ continuously, we deduce from the mean value theorem that there exists $\hat{\va}\in(\va_{1},\va_{0})$ such that $a=a_{\hat{\va}}=\int_{\R^2}|\hat{u}_{\hat{\va}}|^2$.
%and $\mu_{\hat{\va}}=-\f{1}{\hat{\va}^2}<0$.
This implies that the claim \eqref{2:24} holds true, and  we are thus done for this case.

{\em Case 2: $N=3$.} In this case, one can derive from \eqref{N1:1} and \eqref{2.92} that
\begin{equation}\label{2:2.12}
a_\va=a^*\va+o(\va)\quad\hbox{as}\ \ \va\to0.
\end{equation}
Let
$\va_{0}>0$
be a small constant
such that \eqref{GP-4} has a solution $u_{\va_{0}}$ satisfying  \eqref{2.82} and \eqref{02.14}. Setting $\delta=a_{\va_0}=\int_{\R^3}|\hat u_{\va_0}|^2dx$, we then deduce from \eqref{2:2.12} that for any $a\in (0,\delta)$,  there exists a small constant $\va_1\in (0,\va_0)$ such that
\begin{equation*}
  a_{\va_0}=\delta>a>2a^*\va_1>a_{\va_1}.
\end{equation*}
Similar to Case 1, applying the mean value theorem, one can further derive that  the claim \eqref{2:24} holds true for the case $N=3$.

%Similar to Case 1, applying (\ref{2:2.12}) and the mean value theorem, there exists a small constant $\delta>0$ such that for any $a\in(0,\delta)$, {\color{red}there exists a sufficiently small constant $\vaa>0$ such that $a=a_{\vaa}=\int_{\R^3}|\hat u_{\vaa}|^2$. We then conclude that the claim \eqref{2:24} is also true for the case $N=3$.}
%there exists a suitable constant {\color{red}$\mu_a<0$} such that equation \eqref{GP-2} has at least one axially symmetric solution $w_a(x)$.

To complete the proof of Theorem  \ref{thm1}, the rest is to prove that the axially symmetric solution $w_a(x)$ is a normalized concentrating solution of \eqref{GP-2} concentrating at $P_0=(0,\cdots,0,1)\in\R^N$ as $a\to a_*(N)$.
One can note from \eqref{2.82}, \eqref{02.14}, \eqref{R2-1}, \eqref{2.97}, \eqref{2.55} and \eqref{2:2.12} that  $w_a(x)$  satisfies
\begin{equation}\label{Asigpek}
u_a(x):=\sqrt {\f{a}{-\mu_a}} w_a\Big(\f{x}{\sqrt{-\mu_a}} +x_a\Big)e^{-\f{i}{\sqrt{-\mu_a}} x\cdot x_a^\Om}=Q(x)+v_a(x)\ \ \hbox{as}\ \ a\to a_*(N),
\end{equation}
where
\begin{equation}\label{A-sigpek}
  x_a=(0,\cdots,0,p_a)\to (0,\cdots,0,1)\in\R^N\ \ \hbox{and}\ \ \mu_a=-\f{1}{\vaa^2}\to-\infty\ \ \hbox{as}\ \  a\to a_*(N),
\end{equation}
and $v_a(x)\in F_\va$ is a lower order term in the sense $\|v_a\|_{\vaa}=o(1)$ as $a\to a_*(N)$.
If $x_a\in\R^N$ in \eqref{Asigpek} is a global maximum point of $|w_a|$, then one can directly note from above  that $w_a(x)$ is a normalized concentrating solution of \eqref{GP-2} in view of Definition \ref{def2}, and Theorem \ref{thm1} is therefore proved.

However, if $x_a\in\R^N$ in \eqref{Asigpek} is not a global maximum point of $|w_a|$, then we address the proof of Theorem \ref{thm1} as follows.
Following \eqref{Asigpek} and \eqref{A-sigpek}, we first deduce from Lemma \ref{lemA.1} that
\begin{equation*}
  u_a(x)\to Q(x)\ \ \hbox{strongly in}\ \ H^1(\R^N,\C)\ \ \hbox{as}\ \ a\to a_*(N).
\end{equation*}
Using the standard elliptic regularity estimates \cite{HL},
the uniqueness of the global maximum point of $|w_a|$ as $a\to a_*(N)$  can be then obtained by a similar argument of \cite[Proposition 3.3]{GLY}.
Since $w_a(x)$ is axially symmetric, we deduce that the unique global maximum point $\hat x_a\in\R^N$ of $|w_a|$ still lies on the $x_N$-axis.
Therefore, one may assume that  $\hat x_a=(0,\cdots,0,\hat p_a)$.
We claim  that the function
 \begin{equation}\label{F4}
 \begin{aligned}
     \widetilde u_a(x):&=\sqrt{\f{a}{-\mu_a}} w_a\l( \f{x}{\sqrt{-\mu_a}}+\hat x_a\r)e^{-\f{i}{\sqrt{-\mu_a}}  x\cdot \hat x_a^{\Om}}\\
     &=\l[Q\l(x+\sqrt{-\mu_a}(\hat x_a-x_{a})\r)+v_a\l(x+\sqrt{-\mu_a}(\hat x_a-x_{a})\r)\r] e^{-\f{i x\cdot(\hat{x}_a-x_a)^\Om}{\sqrt{-\mu_a}}}
\end{aligned}
 \end{equation}
satisfies
\begin{equation}\label{AAsigpek}
  \widetilde u_a(x)=Q(x)+\widetilde v_a(x)\ \ \hbox{as}\ \ a\to a_*(N),
\end{equation}
where $\widetilde v_a(x)\in H_\va$ satisfies $\|\widetilde v_a\|_{\vaa}=o(1)$ as $a\to a_*(N)$, and the unique global  maximum point $\hat x_a=(0,\cdots,0,\hat p_a)$ satisfies $\hat x_a\to P_0=(0,\cdots,0,1)$ as $a\to a_*(N)$.

In order to prove the above claim (\ref{AAsigpek}), similar to the argument of \cite[Lemma A.2]{GLY},  one can derive from \eqref{NL-1} that
\begin{equation}\label{R3:22}
  \big\|e^{\f{2}{3}|x|}v_a\big\|_{L^\infty(\R^N)}\leq C(-\mu_a)^{-2},\ \ \big\|e^{\f{1}{2}|x|}\nabla v_a\big\|_{L^\infty(\R^N)}\leq C(-\mu_a)^{-2}\ \ \hbox{as}\ \ a\to a_*(N),
\end{equation}
where we have used Lemmas \ref{lem2.1},   \ref{Nlem3.0} and \ref{lemA.1}.
Because the origin is the unique global maximum point of $|\widetilde u_a(x)|$, we obtain from \eqref{F4} that
$\sqrt{-\mu_a}|\hat{x}_a-x_a|=o(1)$ as $a\to a_*(N)$. Moreover, we have
\begin{equation}\label{F2}
\begin{aligned}
    % 0=\f{\p{\widetilde u_\va(0)}}{\p x_j}=\f{\p Q\l(x+\f{\hat x_\va-x_{\va}}{\va}\r)}{\p x_j}\Big|_{x=0}+\f{\p v_\va \l(x+\f{\hat x_\va-x_{\va}}{\va}\r)}{\p x_j}\Big|_{x=0},
    0=&\f{\p{|\widetilde u_a(x)|^2}}{\p x_j}\Big|_{x=0}=\f{\p Q^2\big(x+\sqrt{-\mu_a}(\hat x_a-x_{a})\big)}{\p x_j}\Big|_{x=0}\\
    &+\f{2\p \l\{Q\big(x+\sqrt{-\mu_a}(\hat x_a-x_{a})\big)\operatorname{Re} \l[v_a \big(x+\sqrt{-\mu_a}(\hat x_a-x_{a})\big)\r] \r\}}{\p x_j}\Big|_{x=0}
    \\
    &+\f{\p \l|v_a \big(x+\sqrt{-\mu_a}(\hat x_a-x_{a})\big)\r|^2}{\p x_j}\Big|_{x=0},
\end{aligned}
\end{equation}
where $j=1,2, \cdots, N$.
% Similar to the argument of proving \cite[Lemma A.1]{GLY}, one can derive from Lemmas \ref{lem2.1} and \ref{Nlem3.0} that
We then derive from  \eqref{R3:22} and \eqref{F2} that
% \begin{equation*}
%     C\l|\f{\hat x_\va-x_{\va}}{\va}\r|=O(\va^4)\ \ \hbox{as}\ \ \va\to0,
% \end{equation*}
% which implies that
\begin{equation}\label{F3}
    \l|\sqrt{-\mu_a}(\hat x_a-x_{a})\r|=O\l((-\mu_a)^{-2}\r)\ \ \hbox{as}\ \ a\to a_*(N).
\end{equation}
Combining \eqref{F4} with \eqref{F3} yields that
\begin{equation}\label{F5}
\begin{aligned}
    \widetilde u_a(x):&=\sqrt{\f{a}{-\mu_a}} w_a\l( \f{x}{\sqrt{-\mu_a}}+\hat x_a\r)e^{-\f{i}{\sqrt{-\mu_a}}  x\cdot \hat x_a^{\Om}}\\
&=Q(x)+v_a(x)+O\l((-\mu_a)^{-2}\r):=Q(x)+\widetilde v_a(x)\ \ \hbox{as}\ \ a\to a_*(N).
\end{aligned}
\end{equation}
Following \eqref{Asigpek}, \eqref{F3} and \eqref{F5}, one can further obtain that the above claim \eqref{AAsigpek} holds true.

Since the claim \eqref{AAsigpek} holds true, one can derive from Definition \ref{def2} that $w_a(x)$ is a normalized concentrating solution of \eqref{GP-2} concentrating at $P_0=(0,\cdots,0,1)$ as $a\to a_*(N)$. This completes the proof of Theorem  \ref{thm1}.\qed

\section{Refined Spike Profiles}
Under the assumptions that $0\leq\Omega<\infty$ is fixed,  and  $a>0$ satisfies (\ref{thm1:1}), Theorem \ref{thm1} gives the existence of normalized concentrating solutions for \eqref{GP-2} as $a\to a_*(N)$, which  concentrate at $P_0=(0,\cdots,0,1)$ and are axially symmetric. The main purpose of this section is to derive the refined spike profiles
of normalized concentrating solutions \(w_a\in\mathcal H\) for
\eqref{GP-2} concentrating at \(P_0=(0,\cdots,0,1)\) as \(a\to a_*(N)\).
When \(N=3\), we assume naturally that the maximum point \(x_a\) of \(|w_a|\) lies on the
\(x_3\)-axis.

%Note that the equation
%\eqref{GP-2} is invariant under the rotational transformation. Without loss of generality, we may assume that the global maximum point $x_a\in \R^N$ of $|w_a|$ lies in the $x_N$-axis, i.e., $x_a=(0,\cdots,0,p_a)$.
Denote
\begin{equation}\label{3.3}
\vaa:=\f{1}{\sqrt{-\mu_a}}>0 ,
\end{equation}
and
\begin{equation}\label{3.4}
u_a(x):=\sqrt{a}\vaa w_a(\vaa x+x_a)e^{-i\sigma_{a}-i\vaa x\cdot x_a^\Om},
\end{equation}
where
%$x_{a}=(0,\cdots,0,p_{a})\in\R^N$ is the unique global maximum point of $|w_{a}(x)|$, and
$\mu_a<0$ is as in Definition \ref{def2}, $x_a\in\R^N$ is a global maximum point of $|w_{a}|$, and $\sigma_{a}\in[0,2\pi)$ is chosen properly (similar to \cite[Lemma 3.2]{GLY}) such that
\begin{equation}\label{3:23}
\operatorname{Re}\Big(\inte u_{a}(iQ)\Big)=0.
\end{equation}
It follows from Definition \ref{def2} that
\begin{equation}\label{3.5}
\vaa\to 0\ \ \hbox{and}\ \ u_a(x)=Q(x)+v_{a}(x)\ \ \hbox{as}\ \ a\to a_*(N).
\end{equation}
Here and below  we consider $v_a(x)$ as a lower order term in the sense $\|v_a\|_{\vaa}=o(1)$ as $a\to a_*(N)$, where the norm $\|\cdot\|_{\vaa}$ is
defined  by the inner product \eqref{2.7}. Following \eqref{3.5}, by a similar argument of
\cite[Proposition 3.3]{GLY}, one can get the uniqueness of the global
maximum point \(x_a\) of \(|w_a|\) as \(a\to a_*(N)\). When
\(N=2\), using the rotational invariance of \eqref{GP-2}, without loss of generality we may assume
 that
\[
x_a=(0,p_a),\ \, p_a\to1   \ \ \hbox{as}\ \ a\to a_*(N).
\]
When \(N=3\), the property $(\ref{thmK:1})$ cannot be obtained by an arbitrary
rotation, and hence it is   assumed additionally  in Theorem~\ref{cor3}. Thus, we have for $N=2,3$,
\[
x_a=(0,\cdots,0,p_a),\ \, p_a\to1   \ \ \hbox{as}\ \ a\to a_*(N).
\]

It follows from \eqref{GP-2} and \eqref{3.4} that $u_a(x)$ satisfies
\begin{equation}\label{GP-5}
\big(-\Delta_{\vaa^2x^{\Omega}}\big) u_{a}+\Big[1+\vaa^2\big(|\vaa x+x_a|^2-1\big)^2\Big]u_{a}-|u_{a}|^2u_{a}=0\ \ \hbox{in}\ \ \R^N.
\end{equation}
We then derive
 from \eqref{A1}, \eqref{3.5} and \eqref{GP-5} that $v_{a}(x)$ satisfies
\begin{equation}\label{NL-2}
\begin{split}
&\quad\Big[-\Delta_{\vaa^2x^{\Omega}}+1+\vaa^2\big(|\vaa x+x_a|^2-1\big)^2\Big] v_{a}-Q^2v_{a}-2\operatorname{Re}(Qv_{a})Q\\
&=l_a(x)+R_a(v_{a})\ \,\ \hbox{in}\, \ \R^N,
\end{split}
\end{equation}
where
\begin{equation}\label{3.7}
l_a(x):=-\Big[\vaa^4|x^\Om|^2+\vaa^2\big(|\vaa x+x_a|^2-1\big)^2\Big]Q(x)\ \ \hbox{in}\ \ \R^N,
\end{equation}
and
\begin{equation}\label{3.8}
R_a(v_{a}):=(2\operatorname{Re}(v_{a})v_{a}+|v_{a}|^2)Q+|v_{a}|^2v_{a}\ \ \hbox{in}\ \ \R^N.
\end{equation}
We next follow \eqref{NL-2} to analyze the lower order term $v_a(x)$ as $a\to a_*(N)$.

\subsection{Refined spike profiles in terms of $\mu_a$}
Using above notations, in this subsection we shall derive the refined spike profiles of $w_a(x)$ in terms of $\mu_a<0$. We start with the following lemma on the refined estimates of $v_a(x)$ in \eqref{3.5} and the unique global maximum point $x_a\in\R^N$ of $|w_a|$ as $a\to a_*(N)$.
\begin{lem}\label{lem3.1}
Assume that $0\leq\Omega<\infty$ is fixed, and $a>0$ satisfies (\ref{thm1:1}). Suppose $w_a(x)\in\H$ is a normalized concentrating solution of \eqref{GP-2}  concentrating at $P_0=(0,\cdots,0,1)$ as $a\to a_*(N)$, and $x_a=(0,\cdots,0,p_a)\in\R^N$ is the unique global  maximum point of $|w_a|$ as $a\to a_*(N)$. Then we have
\begin{enumerate}
\item $v_a(x)$ defined in \eqref{3.5} satisfies
\begin{equation}\label{3:22}
|v_a|\leq C\vaa^4e^{-\f{2}{3}|x|},\ \ |\nabla v_a|\leq C\vaa^4e^{-\fh|x|}\quad \hbox{uniformly in}\, \ \R^N\,\ \hbox{as}\, \ a\to a_*(N),
\end{equation}
where $\vaa>0$ defined by \eqref{3.3} satisfies \eqref{3.5}.
\item $p_a\in\R$ satisfies
\begin{equation}\label{3:28}
p_{a}-1=A\vaa^2-\f{A^2}{2}\vaa^4+o(\vaa^4)\ \ \hbox{as}\ \ a\to a_*(N),
\end{equation}
where
\begin{equation}\label{R3}
A:=-\frac{(N+2)\inte |x|^2Q^2dx}{2Na^*}<0.
\end{equation}
\end{enumerate}
		
\end{lem}
\noindent{\bf Proof.}
1. To prove \eqref{3:22}, we first claim that there exists a constant $C>0$ such that
\begin{equation}\label{A3.6}
|v_a|\leq C\delta_a\vaa^2\ \ \hbox{uniformly in}\ \ \R^N\ \ \hbox{as}\ \ a\to a_*(N),
\end{equation}
where $\delta_a>0$ satisfies $\delta_a=o(1)$ as $a\to a_*(N)$,
and
\begin{equation}\label{3.9}
p_{a}-1=A\vaa^2+o(\vaa^2)\ \ \hbox{as}\ \ a\to a_*(N),
\end{equation}
where the constant $A<0$ is defined by \eqref{R3}.

In fact, on the contrary, assume that up to a subsequence if necessary, there exists a constant $C_0>0$ such that
$\frac{\|v_a\|_{L^\infty(\R^N)}}{\vaa^2}\geq C_0$ as $a\to a_*(N)$.
Denote $\mathcal{V}_{a}(x):=\frac{v_a(x)}{\|v_a\|_{L^\infty(\R^N)}}$, so that $\|\mathcal{V}_{a}\|_{L^\infty(\R^N)}=1$.
It then follows from \eqref{NL-2} that $\mathcal{V}_{a}(x)$ satisfies
	%We also deduce from \eqref{3:14} that
\begin{equation}\label{A3.9}
L_{1a}\mathcal{V}_{a}=\frac{N_{1a}}{\|v_a\|_{L^\infty(\R^N)}}\ \ \hbox{in}\ \ \R^N,
\end{equation}
where
\begin{equation*}%\label{3:14}
\begin{aligned}
N_{1a}(x):=l_{a}(x)+R_{a}(v_{a})-2i\vaa^2 x^\Om\cdot\nabla v_{a}\ \ \hbox{in}\ \ \R^N,
\end{aligned}
\end{equation*}
and the operator $L_{1a}$ is defined by
	\begin{equation}\label{3:54}
	L_{1a}:=-\Delta+\vaa^4|x^\Om|^2+1+\vaa^2\l(|\vaa x+x_{a}|^2-1\r)^2-Q^2-2\operatorname{Re}(Q\cdot)Q\,\ \ \hbox{in}\ \ \R^N.
	\end{equation}
By the comparison principle, one can  derive from \eqref{1:3} and \eqref{NL-2}--\eqref{3.8} that
	\begin{equation}\label{A3.7}
	|v_a|\leq C\delta_a e^{-\f{2}{3}|x|}\ \ \hbox{and}\ \ |\nabla v_a|\leq C\delta_a e^{-\fh|x|}\quad \hbox{uniformly in}\ \ \R^N\ \ \hbox{as}\ \ a\to a_*(N),
	\end{equation}
where $\delta_a>0$ satisfies $\delta_a=o(1)$ as $a\to a_*(N)$.
 We then deduce from \eqref{3.7} and \eqref{3.8}  that
	\begin{equation}\label{A3.8}
	\frac{|N_{1a}(x)|}{\|v_a\|_{L^\infty(\R^N)}}\leq C\Big(\f{\vaa^2}{\|v_a\|_{L^\infty(\R^N)}}+\|v_a\|_{L^\infty(\R^N)}\Big)e^{-\frac{1}{10}|x|}\ \ \hbox{uniformly in}\
	\  \R^N
	\end{equation}
as $a\to a_*(N)$, where the constant $C>0$ is independent of $a>0$. Applying the standard elliptic regularity theory, it yields from \eqref{A3.9} that $\|\mathcal{V}_{a}\|_{C^{2,\alpha}_{loc}(\R^N)}\leq C$ for some $\al\in(0,1)$ as $a\to a_*(N)$, where $C>0$ is independent of $a>0$. Therefore, taking a subsequence if necessary, we may assume that $\mathcal{V}_{a}(x)\to \mathcal{V}_0(x)$ in $C^{2}_{loc}(\R^N,\C)$ as $a\to a_*(N)$, where $\mathcal{V}_0(x)$ satisfies the following equation
	\[
	L\mathcal{V}_{0}=\left(-\Delta+1\right) \mathcal{V}_{0}-Q^2\mathcal{V}_{0}-2\operatorname{Re}(Q\mathcal{V}_{0})Q=0\, \ \ \hbox{in}  \ \R^N.
	\]
By the non-degeneracy of $Q$ in \eqref{nond}, we deduce that
	\begin{equation}\label{A3.10}
	\mathcal{V}_{0}=c_0(iQ)+\sum_{j=1}^Nc_j\f{\p Q}{\p x_j}\ \  \hbox{in}\, \ \R^N,
	\end{equation}
	where $c_j\in\R$  are some constants for $j=0, 1,\cdots,N$.

The orthogonality condition \eqref{3:23} gives \(c_0=0\). Moreover,
since the origin is the unique maximum point of \(|u_a|\), we  obtain  from \eqref{3.5} that when $a$ is close sufficiently  to $ a_*(N)$,
\begin{equation*}%\label{3.87}
\begin{split}
0=(\nabla |u_a|^2)(0)=\Big(\nabla\big(2\operatorname{Re}(v_aQ)+|v_a|^2\big)\Big)(0),
\end{split}
\end{equation*}
which yields that
\begin{equation*}
\operatorname{Re}(\nabla \mathcal{V}_{a}(0))=\f{\operatorname{Re}(\nabla v_a(0))}{\|v_a\|_{L^\infty(\R^N)}}=-\f{\operatorname{Re}\big(\nabla v_a(0)\ovl{v}_a(0)\big)}{Q(0)\|v_a\|_{L^\infty(\R^N)}}\to0\ \ \hbox{as}\ \ a\to a_*(N).
\end{equation*}
Thus
\[
\operatorname{Re}\nabla \mathcal V_0(0)=0.
\]
Using the representation of \(\mathcal V_0\) and the fact that \(c_0=0\),
we obtain
\[
0=\operatorname{Re}\nabla \mathcal V_0(0)
=
\sum_{j=1}^N c_j\nabla\partial_{x_j}Q(0)
=
D^2Q(0)(c_1,\ldots,c_N)^T .
\]
Since \(Q\) is radial and has a non-degenerate critical point at the
origin, \(D^2Q(0)\) is invertible. Hence,
\(c_1=\cdots=c_N=0\), and  \(\mathcal V_0\equiv0\).

On the other hand, assume $y_{a}\in\R^N$ is a global maximum point of $|\mathcal V_{a}(x)|$ such that $$|\mathcal V_{a}(y_{a})|
	=\max\limits_{x\in\R^N}\frac{|v_a(x)|}{\|v_a\|_{L^\infty(\R^N)}}=1.$$
By the comparison principle, one  then can deduce from \eqref{A3.9} and \eqref{A3.8} that $|y_{a}|\leq C$ uniformly in $a>0$. We thus conclude from above that up to a subsequence if necessary,
	$1=|\mathcal V_{a}(y_{a})|\to |\mathcal V_{0}(y_{0})|$ for some $y_{0}\in\R^N$ as $a\to a_*(N)$,
which however contradicts the fact that $\mathcal V_{0}(x)\equiv 0$ in $\R^N$. We then derive that \eqref{A3.6} holds true.

By the comparison principle,  one can obtain from \eqref{A3.6} that
	\begin{equation}\label{A3.17}
	|v_a|\leq C\vaa^2 e^{-\f23|x|}\ \ \hbox{and}\ \ |\nabla v_a|\leq C\vaa^2 e^{-\f12|x|}\ \ \hbox{uniformly in}\ \ \R^N\ \ \hbox{as}\ \ a\to a_*(N).
	\end{equation}
Applying \eqref{A3.17}, the similar argument of Lemma \ref{Nlem3.0} yields that \eqref{3.9} holds true.	
Based on  \eqref{3.9} and \eqref{A3.17}, the same argument of \eqref{A3.17} gives that \eqref{3:22} holds true. %This completes the proof of Lemma \ref{lem3.1} (1).

2. Following \eqref{3:22} and \eqref{3.9}, the similar argument of  Lemma \ref{Nlem3.0} yields that \eqref{3:28} is true. We therefore conclude that Lemma \ref{lem3.1} holds true.  \qed

We next use Lemma \ref{lem3.1} to derive the following refined spike profiles of  $w_{a}(x)$ in terms of $\vaa>0$.

\begin{lem}\label{lem3.3}
	Under the assumptions of Lemma \ref{lem3.1}, let $w_a(x)\in\H$ be  the normalized concentrating solution   of \eqref{GP-2} concentrating at $P_0=(0,\cdots,0,1)$ as $a\to a_*(N)$. Then $u_a(x)$ defined by \eqref{3.4} satisfies
	\begin{equation}\label{3:8}
	\begin{aligned}
	u_{a}(x)&=Q(x)+v_a(x)\\
	&=Q(x)+\vaa^4\psi_{1}(x)+\vaa^5\psi_{2}(x)+\vaa^6\psi_{3}(x)+o(\vaa^6)\ \ \hbox{in}\ \ \R^N\ \ \hbox{as}\ \ a\to a_*(N),
	\end{aligned}
	\end{equation}
	where $\psi_{j}(x)\in C^2(\R^N,\C)\cap L^\infty(\R^N,\C)$ solves uniquely
	\begin{equation}\label{3:9}
	\operatorname{Re} (\nabla \psi_{j}(0))=0,\ \operatorname{Re}\Big(\inte \psi_j(iQ)\Big)=0,\quad L\psi_{j}(x)=f_j(x)
	\,\ \mbox{in}\,\ \R^N,\ j=1,2,3,
	\end{equation}
	and $f_j(x)$ satisfies
	\begin{equation}\label{3:35}\arraycolsep=1.5pt
	f_j(x)=\left\{\begin{array}{lll}
	&-\big(|x^\Om|^2+4x_{N}^2\big)Q(x),\qquad \quad \ &\mbox{if}\ \ j=1;\\[3mm]
	&-\big(4|x|^2+8A\big)Q(x)x_{N}, \,\ &\mbox{if}\ \ j=2;\\[2mm]
	&  -\big(|x|^4+4A|x|^2+8Ax_{N}^2+4A^2\big)Q(x) -2ix^\Om\cdot\nabla\psi_1 , \,\ &\mbox{if}\ \ j=3;
	\end{array}\right.\end{equation}
	where $x=(x_{1}, \cdots, x_{N})\in\R^N$, and $x^\Om$ is as in \eqref{1.1}. Here $L$ and $A<0$ are defined by  \eqref{lin op} and \eqref{R3}, respectively.
	%where $\psi_{i}(x)\in C^2(\R^3)\cap L^\infty(\R^3)$ solves uniquely
	%\begin{equation}\label{3:9}
	%  \nabla\psi_{1}(0)=0,\ \ L\psi_{1}(x)=-\big(|x^\Om|^2+4x_{3}^2\big)Q(x)\ \ \hbox{in}\ \ \R^3.
	%\end{equation}
\end{lem}
\noindent{\bf Proof.}
Denote
\begin{equation}\label{3:10}
\mathcal{V}_{1a}(x):=v_{a}(x)-\vaa^4\psi_{1}(x)\ \ \hbox{in}\ \ \R^N.
\end{equation}
By the comparison principle, we derive from \eqref{1:3}, \eqref{3:9} and \eqref{3:35} that for $j=1, 2,3,$
\begin{equation}\label{3:11}
|\psi_{j}(x)|\leq Ce^{-\delta|x|}\ \ \hbox{uniformly in}\ \ \R^N,\ \ \frac{4}{5}<\delta<1.
\end{equation}
Using Lemma \ref{lem3.1}, the similar arguments of proving \eqref{A3.6} and   \cite[(2.42)]{Guo} give that
\begin{equation}\label{3:42}
\l|\mathcal{V}_{1a}(x)\r|,\ \l|\nabla \mathcal{V}_{1a}(x)\r|\leq C\delta_{a}\vaa^4e^{-\f{1}{11}|x|}\ \ \hbox{uniformly in}\ \ \R^N\ \ \hbox{as}\ \  a\to a_*(N),
\end{equation}
where $\delta_{a}>0$ satisfies $\delta_{a}=o(1)$ as $a\to a_*(N)$.

We next prove that \eqref{3:8} holds true.
In fact, denote
\begin{equation}\label{3:37}
\mathcal{V}_{2a}(x):=v_{a}(x)-\vaa^4\psi_{1}(x)-\vaa^5\psi_{2}(x)-\vaa^{6}\psi_{3}(x)\,\ \ \hbox{in}\, \ \R^N,
\end{equation}
where $\psi_{j}(x)$ is a solution of \eqref{3:9} for $j=1, 2, 3$. The uniqueness of $\psi_{j}(x)$ follows from \eqref{3:9} and the property \eqref{nond}.
It then follows from \eqref{NL-2} and \eqref{3:37} that
\begin{equation*}%\label{3:38}
\begin{aligned}
L_{1a}\mathcal{V}_{2a}&=l_{a}(x)+R_{a}(v_{a})-2i\vaa^2 x^\Om\cdot\nabla v_{a}-L\big(\vaa^4\psi_{1}(x)+\vaa^5\psi_{2}(x)+\vaa^6\psi_{3}(x)\big)\\
&\quad-(L_{1a}-L)\big(\vaa^4\psi_{1}(x)+\vaa^5\psi_{2}(x)+\vaa^6\psi_{3}(x)\big):=N_{2a}(x)\ \ \hbox{in}\ \ \R^N,
\end{aligned}
\end{equation*}
where the operator $L_{1a}$ is defined by \eqref{3:54}.
Applying Lemma \ref{lem3.1},  one can obtain from  \eqref{3.7}, \eqref{3.8}, \eqref{3:9}, \eqref{3:35}, \eqref{3:11} and \eqref{3:42} that
\begin{equation*}%\label{3:15}
\frac{\l|N_{2a}(x)\r|}{\vaa^6}\leq C\delta_{a}e^{-\frac{1}{12}|x|}\ \ \hbox{uniformly in}\ \ \R^N\ \ \hbox{as}\ \ a\to a_*(N),
\end{equation*}
where $\delta_{a}>0$ satisfies $\delta_{a}=o(1)$ as $a\to a_*(N)$. The same argument of proving \eqref{A3.6} gives that
\begin{equation*}
\mathcal{V}_{2a}=o(\vaa^6)\ \ \hbox{as}\ \ a\to a_*(N),
\end{equation*}
together with \eqref{3:37}, which yields that Lemma \ref{lem3.3} holds true.
\qed
\vskip 0.05truein

\vskip 0.05truein

\subsection{Refined spike profiles in terms of $\big|a_*(N)-a\big|$}
Following the estimates of the previous subsection, the main purpose of this subsection is  to derive the refined spike profiles of $w_a(x)$ in terms of  $\big|a_*(N)-a\big|$. Towards this aim,
we first give  the  refined estimates of $\vaa>0$ and $\mu_a<0$ for the case $N=2$, which implies from (\ref{A:1.16}) and (\ref{A:1.17}) that
\begin{equation}\label{A:R3:2}
a\nearrow a^*\ \ \Longleftrightarrow\ \ a\to a_{*}(2)=a^*=\int_{\R^2} Q^2dx.
\end{equation}

\begin{lem}\label{Nlem3.5}
	Under the assumptions of Lemma \ref{lem3.1} with $N=2$, let $w_a(x)\in\H$ be the normalized concentrating solution   of \eqref{GP-2} concentrating at $P_0=(0,1)$ as $a\nearrow a^*$, and suppose $\va _a>0$ and  $u_a(x)$ are defined by \eqref{3.3} and (\ref{3.4}), respectively. Then we have
	\begin{equation}\label{R3:2}
	\begin{aligned}
	\ap:=\Big(\f{a^*-a}{B_{1}}\Big)^\f{1}{4},\ \  \vaa=&\ap-\f{B_{2}}{4B_{1}}\ap^3
	+\f{9B_{2}^2-8B_{1}B_{3}}{32B_{1}^2}\ap^5
	+o\l(\ap^5\r)\ \ \hbox{as}\ \ a\nearrow a^*,
	\end{aligned}
	\end{equation}
	and
	\begin{equation}\label{R3.3}
	\begin{aligned}
	-\mu_{a}\ap^2=\Big(\f{\ap}{\vaa}\Big)^{2}=1+\f{B_{2}}{2B_{1}}\ap^2
	+\f{4B_{1}B_{3}-3B_{2}^2}{8B_{1}^2}\ap^4+o\l(\ap^4\r)\ \ \hbox{as}\ \ a\nearrow a^*,
	\end{aligned}
	\end{equation}
	where
	$B_{1}:=\f{\Om^2+8}{4}\int_{\R^2}|x|^2Q^2$, $B_{2}:=2\int_{\R^2}(|x|^4+4A|x|^2)Q^2$, and $B_{3}:=3\int_{\R^2}(|x^\Om|^2+4x_{2}^2)Q\psi_{1}$. Here the constant $A<0$ is as in \eqref{R3}, and $\psi_1(x)$ is defined by \eqref{3:9}.
\end{lem}
\noindent{\bf Proof.}
Similar to Lemma \ref{lem2.6},  one can check that $u_a$ satisfies
\begin{equation}\label{3:39}
\begin{split}
&\quad\int_{\R^2} \Big(\frac{1}{2}|u_a|^4-|u_a|^2\Big)\\
&=\int_{\R^2} \Big[2\vaa^4|x^\Om|^2+2\vaa^3(|\vaa x+x_a|^2-1)(\vaa |x|^2+x\cdot x_a)\\
&\quad\qquad+\vaa^2\big(|\vaa x+x_a|^2-1\big)^2\Big]|u_a|^2
-\operatorname{Re}\Big(2i\vaa^2 \int_{\R^2} (x^\Om\cdot\nabla u_a)(x\cdot\nabla \ovl u_a)\Big).
\end{split}
\end{equation}
 By  \eqref{N1:1}, \eqref{3.4}, \eqref{3.5} and \eqref{3:22}, direct  calculations yield  that
\begin{equation}\label{3:46}
\begin{aligned}
&\hbox{LHS of \eqref{3:39}}\ \ \\
% =&\fh\int_{\R^2}\Big[Q^4+4Re(Q^3v_{a})+\l(2Q^2|v_{a}|^2+4(Re(Qv_{a}))^2\r)\Big]-a+O(\vaa^{12})\\
=&a^*-a+\int_{\R^2}\Big[2\operatorname{Re}(Q^3v_{a})+Q^2|v_{a}|^2+2(\operatorname{Re}(Qv_{a}))^2\Big]+O(\vaa^{12})\ \ \hbox{as}\ \ a\nearrow a^*.\\
\end{aligned}
\end{equation}
We next estimate the right hand side of \eqref{3:46}.

Similar to \eqref{2.521}, we obtain from \eqref{A1} and  \eqref{NL-2} that
\begin{equation}\label{3:40}
\begin{aligned}
\operatorname{Re}\Big(\int_{\R^2}Q^3v_{a}\Big)
% =&Re\l(\int_{\R^2}(-\Delta Q+Q)v_{a}\r)
% =Re\l(\int_{\R^2}(-\Delta v_{a}+v_{a})Q\r)\\
=&\operatorname{Re}\Big(\int_{\R^2}\Big[-2i\vaa^2 (x^\Om\cdot\nabla v_{a})-\vaa^2\big(|\vaa x+x_a|^2-1\big)^2v_{a}\\
&\qquad\qquad-\vaa^4|x^\Om|^2v_{a}+3Q^2v_{a}+l_{a}(x)+R_{a}(v_{a})\Big]Q\Big).\\
\end{aligned}
\end{equation}
 By  the definition of $x^\Om$ in \eqref{1.1}, note  that
\begin{equation*}%\label{3:41}
\begin{aligned}
\operatorname{Re}\Big(\int_{\R^2}-2i\vaa^2 (x^\Om\cdot\nabla v_{a})Q\Big)
= \operatorname{Re}\Big(\int_{\R^2}2i\vaa^2 (x^\Om\cdot\nabla Q)v_{a}\Big)=0.
\end{aligned}
\end{equation*}
It then follows from \eqref{3.7}, \eqref{3.8}
%\eqref{3.9}
and \eqref{3:40} that
\begin{equation}\label{3:50}
\begin{aligned}
&2\operatorname{Re}\Big(\int_{\R^2}Q^3v_{a}\Big)\\
=&\operatorname{Re}\Big(\int_{\R^2}\Big\{\vaa^4|x^\Om|^2v_{a}
+\Big[\vaa^4|x^\Om|^2+\vaa^2\big(|\vaa x+x_a|^2-1\big)^2\Big]Q\\
&\qquad+\vaa^2\big(|\vaa x+x_a|^2-1\big)^2v_{a}-\l(2\operatorname{Re}(v_{a})v_{a}+|v_{a}|^2\r)Q-|v_{a}|^2v_{a}\Big\}Q\Big).
%+o(\vaa^8)\ \ \hbox{as}\ \ a\nearrow a^*.\\
\end{aligned}
\end{equation}
%Using Lemma \ref{lem3.3},
Using \eqref{3:22}, we obtain from \eqref{3:46} and \eqref{3:50} that
\begin{equation}\label{3:52}
\begin{aligned}
&\hbox{LHS of \eqref{3:39}}\\
% =&a^*-a
% +Re\Big(\int_{\R^2}\Big\{\vaa^4|x^\Om|^2v_{a}+\vaa^2V_{\Om}(\vaa x+x_{a})v_{a}
% +\l(\vaa^4|x^\Om|^2+\vaa^2V_{\Om}(\vaa x+x_{a})\r)Q\\
% &-\l(2Re(v_{a})v_{a}+|v_{a}|^2\r)Q\Big\}Q\Big)
% +\int_{\R^2}\Big[Q^2|v_{a}|^2+2(Re(Qv_{a}))^2\Big]+o(\vaa^{8})\\
=&a^*-a
+\operatorname{Re}\Big(\int_{\R^2}\Big\{\vaa^4|x^\Om|^2v_{a}+\vaa^2\big(|\vaa x+x_a|^2-1\big)^2v_{a}\\
&\quad+\Big[\vaa^4|x^\Om|^2+\vaa^2\big(|\vaa x+x_a|^2-1\big)^2\Big]Q\Big\}Q\Big)
+o(\vaa^{8})\ \ \hbox{as}\ \ a\nearrow a^*.\\
\end{aligned}
\end{equation}

On the other hand, we obtain from Lemma \ref{lem3.3} that
\begin{equation}\label{3:51}
\begin{aligned}
&\operatorname{Re}\Big(2i\vaa^2\int_{\R^2}(x^\Om\cdot\nabla u_a)(x\cdot\nabla \ovl u_a)\Big)\\
&=\operatorname{Re}\Big(2i\vaa^2\int_{\R^2}\Big[x^\Om\cdot\nabla \l(Q+\vaa^4\psi_{1}+\vaa^5\psi_{2}+\vaa^6\psi_{3}+o(\vaa^6)\r)\Big]\\
&\quad\qquad\cdot\Big[x\cdot\nabla \l(Q+\vaa^4\psi_{1}+\vaa^5\psi_{2}+\vaa^6\ovl{\psi}_{3}+o(\vaa^6)\r)\Big]\Big)\\
&=o(\vaa^{8})\ \ \hbox{as}\ \ a\nearrow a^*,
\end{aligned}
\end{equation}
where $\psi_{j}$ are defined by \eqref{3:9} for $j=1, 2, 3$, and we have used the fact that
\begin{equation*}%\label{3:49}
\int_{\R^2}(x^\Om\cdot\nabla\psi_{j})(x\cdot\nabla Q)=-\int_{\R^2}\psi_{j}x^\Om\cdot\nabla(x\cdot\nabla Q)=0.
\end{equation*}
Combining \eqref{3:39}, \eqref{3:52} and \eqref{3:51}, we derive from Lemmas \ref{lem3.1} and \ref{lem3.3}
that
\begin{equation}\label{3:53}
\begin{aligned}
a^*-a=&\int_{\R^2} \Big[2\vaa^4|x^\Om|^2+\vaa^2\big(|\vaa x+x_a|^2-1\big)^2\\
&\qquad+2\vaa^3(|\vaa x+x_a|^2-1)(\vaa |x|^2+x\cdot x_a)\Big]|u_a|^2\\
&-\operatorname{Re}\Big(\int_{\R^2}\Big\{\vaa^4|x^\Om|^2v_{a}+\vaa^2\big(|\vaa x+x_a|^2-1\big)^2v_{a}\\
&\quad\qquad+\Big[\vaa^4|x^\Om|^2+\vaa^2\big(|\vaa x+x_a|^2-1\big)^2\Big]Q\Big\}Q\Big)+o(\vaa^{8})\\
=&\operatorname{Re}\Big(\int_{\R^2}\Big\{\vaa^4|x^\Om|^2Q^2+3\vaa^4|x^\Om|^2v_{a}Q+\vaa^2\big(|\vaa x+x_a|^2-1\big)^2Qv_{a}\\
&\qquad+2\vaa^3(|\vaa x+x_{a}|^2-1)(\vaa |x|^2+x_{2}p_{a})(Q^2+2Qv_{a})\Big\}\Big)+o(\vaa^{8})\\
% =&\f{\Om^2}{4}\vaa^4\int_{\R^2}|x|^2Q^2
% +3\vaa^8\int_{\R^2}|x^\Om|^2Q\psi_{1}+4\vaa^8\int_{\R^2}x_{2}^2Q\psi_{1}\\
% &+2\vaa^3Re\Big(\int_{\R^2}\big[\vaa^3|x|^4+\l(\vaa^2|x|^2+2\vaa^2|x|^2+(p_{a}^2-1)\r)p_{a}x_{2}\\
% &+2\vaa p_{a}^2x_{2}^2+\vaa(p_{a}^2-1)|x|^2\big](Q^2+2Qv_{a})\Big)+o(\vaa^{8})\\
% =&\f{\Om^2}{4}\vaa^4\int_{\R^2}|x|^2Q^2
% +3\vaa^8\int_{\R^2}|x^\Om|^2Q\psi_{1}+4\vaa^8\int_{\R^2}x_{2}^2Q\psi_{1}\\
% &+2\vaa^3Re\Big(\int_{\R^2}\big[\vaa^3|x|^4+\l(\vaa^2|x|^2+2\vaa^2|x|^2+(p_{a}^2-1)\r)p_{a}x_{2}\\
% &+2\vaa p_{a}^2x_{2}^2+\vaa(p_{a}^2-1)|x|^2\big](Q^2+2Qv_{a})\Big)+o(\vaa^{8})\\
=&\f{\Om^2+8}{4}\vaa^4\int_{\R^2}|x|^2Q^2
+2\vaa^6\int_{\R^2}(|x|^4+4A|x|^2)Q^2\\
&+3\vaa^8\int_{\R^2}(|x^\Om|^2+4x_{2}^2)Q\psi_{1}
+o(\vaa^8)\\
:=&B_{1}\vaa^4+B_{2}\vaa^6+B_{3}\vaa^8+o(\vaa^8)\ \ \hbox{as}\ \ a\nearrow a^*,
\end{aligned}
\end{equation}
% Moreover, we have
% \begin{equation*}
% \begin{aligned}
%   &2\vaa^3Re\l(\int_{\R^2}(|\vaa x+x_{a}|^2-1)(\vaa |x|^2+x_{2}p_{a})(Q^2+2Qv_{a})\r)\\
%   =&2\vaa^3Re\Big(\int_{\R^2}\big[\vaa^3|x|^4+\l(\vaa^2|x|^2+2\vaa^2|x|^2+(p_{a}^2-1)\r)p_{a}x_{2}\\
%   &+2\vaa p_{a}^2x_{2}^2+\vaa(p_{a}^2-1)|x|^2\big](Q^2+2Qv_{a})\Big)\\
%   =&2\vaa^3\int_{\R^2}\Big\{\vaa^3|x|^4+2\vaa[1+2A\vaa^2+o(\vaa^4)]x_{2}^2
%   +\vaa[2A\vaa^2+o(\vaa^4)]|x|^2\Big\}Q^2\\
%   &+8\vaa^8\int_{\R^2}x_{2}^2 Q\psi_{1}+o(\vaa^8)\\
%   % =&4\vaa^4\int_{\R^2}x_{2}^2Q^2+2\vaa^6\l(\int_{\R^2}|x|^4Q^2+4A\int_{\R^2}x_{2}^2Q^2
%   % +2A\int_{\R^2}|x|^2Q^2\r)\\
%   % &+8\vaa^8\int_{\R^2}x_{2}^2 Q\psi_{1}+o(\vaa^8)\\
% =&2\vaa^4\int_{\R^2}|x|^2Q^2+2\vaa^6\int_{\R^2}(|x|^4+4A|x|^2)Q^2+8\vaa^8\int_{\R^2}x_{2}^2 Q\psi_{1}\\
%   &+o(\vaa^8)\ \ \hbox{as}\ \ a\nearrow a^*.\\
% \end{aligned}
% \end{equation*}
% We then obtain from \eqref{3:53} that
% \begin{equation*}%\label{3-1}
% \begin{aligned}
%   a^*-a=&\f{\Om^2+8}{4}\vaa^4\int_{\R^2}|x|^2Q^2
%   +2\vaa^6\int_{\R^2}(|x|^4+4A|x|^2)Q^2\\
%   &+3\vaa^8\int_{\R^2}(|x^\Om|^2+4x_{2}^2)Q\psi_{1}
%   +o(\vaa^8)\\
%   :=&B_{1}\vaa^4+B_{2}\vaa^6+B_{3}\vaa^8+o(\vaa^8)\ \ \hbox{as}\ \ a\nearrow a^*,
% \end{aligned}
% \end{equation*}
which gives  that
\begin{equation}\label{3-3}
\f{a^*-a}{B_{1}\vaa^4}=1+\f{B_{2}}{B_{1}}\vaa^2+\f{B_{3}}{B_{1}}\vaa^4+o(\vaa^4)\quad  \hbox{as}\ \ a\nearrow a^*.
\end{equation}
We next employ \eqref{3-3} to prove that both \eqref{R3:2} and \eqref{R3.3} hold true.

Denoting $\ap:=\Big(\f{a^*-a}{B_{1}}\Big)^\f{1}{4}$, we deduce from \eqref{3-3} that
\begin{equation}\label{3-2}
\vaa=\ap+o\l(\ap\r)\ \ \hbox{as}\ \ a\nearrow a^*.
\end{equation}
Substituting \eqref{3-2} into \eqref{3-3}, it yields that
\begin{equation}\label{3-4}
\Big(\f{\ap}{\vaa}\Big)^4=1+\f{B_{2}}{B_{1}}\ap^2+o\l(\ap^2\r)\ \ \hbox{as}\ \ a\nearrow a^*,
\end{equation}
which implies that
\begin{equation*}
\vaa=\ap-\f{B_{2}}{4B_{1}}\ap^3
+o\l(\ap^3\r)\ \ \hbox{as}\ \ a\nearrow a^*.
\end{equation*}
Similar to \eqref{3-4}, we have
\begin{equation}\label{3-5}
\begin{aligned}
\Big(\f{\ap}{\vaa}\Big)^4
% &=1+\f{B_{2}}{B_{1}}\l(\ap^2
% -\f{B_{2}}{2B_{1}}\ap^4\r)
% +\f{B_{3}}{B_{1}}\ap^4+o\l(\ap^4\r)\\
&=1+\f{B_{2}}{B_{1}}\ap^2+\f{2B_{1}B_{3}-B_{2}^2}{2B_{1}^2}\ap^4
+o\l(\ap^4\r)\ \ \hbox{as}\ \ a\nearrow a^*,
\end{aligned}
\end{equation}
from which we deduce that
\begin{equation*}
\begin{aligned}
\vaa=&\ap-\f{B_{2}}{4B_{1}}\ap^3
+\f{9B_{2}^2-8B_{1}B_{3}}{32B_{1}^2}\ap^5
+o\l(\ap^5\r)\ \ \hbox{as}\ \ a\nearrow a^*,
\end{aligned}
\end{equation*}
and hence \eqref{R3:2} holds true.
Following the definition of $\vaa$ in \eqref{3.3}, we also derive from \eqref{3-5} that
\begin{equation*}
-\mu_{a}\ap^2=\Big(\f{\ap}{\vaa}\Big)^{2}=1+\f{B_{2}}{2B_{1}}\ap^2
+\f{4B_{1}B_{3}-3B_{2}^2}{8B_{1}^2}\ap^4+o(\ap^4)\ \,\ \hbox{as}\ \ a\nearrow a^*,
\end{equation*}
which gives \eqref{R3.3}.
This completes the proof of Lemma \ref{Nlem3.5}.
\qed

\vskip 0.05truein
One can derive from Lemmas \ref{lem3.3} and \ref{Nlem3.5} that for $N=2$,
\begin{equation}\label{3-7}
\begin{aligned}
u_{a}(x):&=\sqrt{a}\vaa w_a(\vaa x+x_a)e^{-i\sigma_{a}-i\vaa x\cdot x_a^\Om}\\
&=Q(x)+\ap^4\psi_{1}(x)+o(\ap^4)\ \ \hbox{as}\ \ a\nearrow a^*,
\end{aligned}
\end{equation}
where $\sigma_{a}\in[0,2\pi)$ is chosen properly such that \eqref{3:23} holds true.
Moreover, we have the following  refined spike profile of the normalized concentrating solution $w_a(x)$ in terms of $a^*-a$.

\begin{lem}\label{lem3.6}
Under the assumptions of Lemma \ref{lem3.1} with $N=2$, let $w_a(x)\in\H$ be the normalized concentrating solution  of \eqref{GP-2} concentrating at $P_0=(0,1)$ as $a\nearrow a^*$. Then we have
\begin{equation*}%\label{R3.5}
\begin{aligned}
\widetilde{u}_{a}(x):&=\sqrt{a}\ap w_{a}(\ap x+x_{a})e^{-i\widetilde\sigma_{a}-i\ap x\cdot x_{a}^\Om}\\
&=Q(x)+\ap^2C_{1}(x)+\ap^4C_{2}(x)+o(\ap^4)\ \ \hbox{as}\ \ a\nearrow a^*,
	%&=Q(x)+\f{B_{2}}{4B_{1}}\ap^2\big(Q(x)+x\cdot\nabla Q\big)+\ap^4\Big[\f{8B_{1}B_{3}-7B_{2}^2}{32B_{1}^2}\big(Q(x)+x\cdot\nabla Q\big)\\
	%&\quad+\f{B_{2}^2}{16B_{1}^2}(x\cdot\nabla Q)+\f{B_{2}^2}{32B_{1}^2}x^T\cdot\nabla^2 Q\cdot x+\psi_{1}(x)\Big]+i\f{B_{2}\Om}{8B_{1}}\ap^3x_{1}Q\\
	%&\quad+o(\ap^4)\ \ \hbox{as}\ \ a\nearrow a^*,
\end{aligned}
\end{equation*}
where $\ap:=\Big(\f{a^*-a}{B_{1}}\Big)^\f{1}{4}>0$, $x_a$ denotes the unique global maximum point of $|w_a|$ as $a\nearrow a^*$, and  the constant phase  $\widetilde\sigma_{a}\in[0,2\pi)$ is chosen properly such that
\begin{equation}\label{R3.6}
\operatorname{Re}\Big(\int_{\R^2}  \widetilde{u}_{a}(iQ)\Big)=0.
\end{equation}
Here the function $C_{1}(x)$ and $C_{2}(x)$ satisfy
\begin{equation}\label{R3-1}
C_{1}(x):=\f{B_{2}}{4B_{1}}\l(Q+x\cdot\nabla Q\r),
\end{equation}
and
\begin{equation}\label{R3-2}
\begin{aligned}
C_{2}(x):=&\f{8B_{1}B_{3}-7B_{2}^2}{32B_{1}^2}\big(Q+x\cdot\nabla Q\big)
+\f{B_{2}^2}{32B_{1}^2}x^T(\nabla^2 Q) x\\
&+\f{B_{2}^2}{16B_{1}^2}\big(x\cdot\nabla Q\big)+\psi_{1}(x),
\end{aligned}
\end{equation}
where $\psi_1(x)$ is defined by \eqref{3:9}, and the constants  $B_j$ are as in Lemma \ref{Nlem3.5} for $j=1, 2, 3$.
\end{lem}

%similar to \cite[Lemma 3.2]{GLY},

\noindent{\bf Proof.}
Note from \eqref{R3.3} that
\begin{equation}\label{R3.7}
\f{\ap}{\vaa}=1+\f{B_{2}}{4B_{1}}\ap^2+\f{8B_{1}B_{3}-7B_{2}^2}{32B_{1}^2}\ap^4+o(\ap^4)\ \ \hbox{as}\ \ a\nearrow a^*.
\end{equation}
Using Taylor's expansion, we then derive from \eqref{3-7} and \eqref{R3.7} that
\begin{equation}\label{R3.8}
\begin{aligned}
&\sqrt{a}\ap w_{a}(\ap x+x_{a})e^{-i\sigma_{a}-i\ap x\cdot x_a^\Om}\\
=&\f{\ap}{\vaa}u_{a}\Big(\f{\ap}{\vaa} x\Big)
% =&\f{\ap}{\vaa}\l(Q\big(\f{\ap}{\vaa} x\big)+\ap^4\psi_{1}\big(\f{\ap}{\vaa} x\big)+o(\ap^4)\r)\\
% =&\l(1+\f{B_{2}}{4B_{1}}\ap^2+\f{8B_{1}B_{3}-7B_{2}^2}{32B_{1}^2}\ap^4+o(\ap^4)\r)\\
% &\quad\cdot\Big[Q(x)+\l(\f{\ap}{\vaa}-1\r)(x\cdot\nabla Q)+\l(\f{\ap}{\vaa}-1\r)^2\f{x^T(\nabla^2 Q) x}{2}
% +\ap^4\psi_{1}(x)+o(\ap^4)\Big]\\
=Q(x)+\ap^2C_{1}(x)+\ap^4C_{2}(x)+o(\ap^4)\ \ \hbox{as}\ \ a\nearrow a^*,
%=&Q(x)+\f{B_{2}}{4B_{1}}\ap^2\l(Q(x)+x\cdot\nabla Q\r)
%+\ap^4\Big[\f{8B_{1}B_{3}-7B_{2}^2}{32B_{1}^2}\big(Q(x)+x\cdot\nabla Q\big)\\
%&\quad+\f{B_{2}^2}{32B_{1}^2}x^T\nabla^2 Q x+\f{B_{2}^2}{16B_{1}^2}\big(x\cdot\nabla Q\big)+\psi_{1}%(x)\Big]
% +o(\ap^4)\ \ \hbo{as}\ \ a\nearrow a^*.
\end{aligned}
\end{equation}
where $C_{1}(x)$ and $C_{2}(x)$ are defined by \eqref{R3-1} and \eqref{R3-2}, respectively.

It follows from \eqref{R3.6} and
\eqref{R3.8} that
\begin{equation*}%\label{R3.9}
\begin{aligned}
0=&\operatorname{Re}\Big(\int_{\R^2}  \widetilde{u}_{a}(iQ)\Big)\\
% =&Re\l(\int_{\R^2} \sqrt{a}\ap w_{a}(\ap x+x_{a})e^{-i\widetilde{\sigma}_{a}-i\ap x\cdot x_{a}^\Om}(iQ)\r)\\
=&\operatorname{Re}\Big(\int_{\R^2}\sqrt{a}\ap w_{a}(\ap x+x_{a})e^{-i\sigma_{a}-i\ap x\cdot x_{a}^\Om}e^{-i(\widetilde{\sigma}_{a}-\sigma_{a})}
(iQ)\Big)\\
=&\operatorname{Re}\Big(\int_{\R^2}\Big(Q(x)+\ap^2C_{1}(x)+\ap^4C_{2}(x)\Big)
e^{-i(\widetilde{\sigma}_{a}-\sigma_{a})}
(iQ)\Big)
+o(\ap^4)\\
=&\sin(\widetilde{\sigma}_{a}-\sigma_{a})\int_{\R^2}\Big(Q(x)+\ap^2C_{1}(x)+\ap^4C_{2}(x)\Big)Q
+o(\ap^4)\ \ \hbox{as}\ \ a\nearrow a^*,
\end{aligned}
\end{equation*}
which implies that
\begin{equation}\label{R3.10}
|\widetilde{\sigma}_{a}-\sigma_{a}|=o(\ap^4)\ \ \hbox{as}\ \ a\nearrow a^*.
\end{equation}
Following  \eqref{R3.8} and \eqref{R3.10}, we have
\begin{equation*}
\begin{aligned}
\widetilde{u}_{a}(x):&=\sqrt{a}\ap w_{a}(\ap x+x_{a})e^{-i\widetilde\sigma_{a}-i\ap x\cdot x_{a}^\Om}\\
&=\sqrt{a}\ap w_{a}(\ap x+x_{a})e^{-i\sigma_{a}-i\ap x\cdot x_{a}^\Om}e^{-i(\widetilde{\sigma}_{a}-\sigma_{a})}\\
&=\Big(Q(x)+\ap^2C_{1}(x)+\ap^4C_{2}(x)\Big)e^{-i(\widetilde{\sigma}_{a}-\sigma_{a})}+o(\ap^4)\\
% &=\Big(Q(x)+\ap^2C_{1}(x)+\ap^4C_{2}(x)\Big)\Big[1+i(\vaa-\ap)x\cdot x_{a}^\Om\Big]+o(\ap^4)\\
&=Q(x)+\ap^2C_{1}(x)+\ap^4C_{2}(x)
+o(\ap^4)\ \ \hbox{as}\ \ a\nearrow a^*.
%&=\Big[Q(x)+\f{B_{2}}{4B_{1}}\ap^2\big(Q+x\cdot\nabla Q\big)+\ap^4\Big(\f{8B_{1}B_{3}-7B_{2}^2}{32B_{1}^2}\big(Q+x\cdot\nabla Q\big)\\
%&\quad+\f{B_{2}^2}{16B_{1}^2}(x\cdot\nabla Q)+\f{B_{2}^2}{32B_{1}^2}x^T\cdot\nabla^2 Q\cdot x+\psi_{1}(x)\Big)\Big]e^{-i(\widetilde{\sigma}_{a}-\sigma_{a})-i(\ap-\vaa)(x\cdot x_{a}^\Om)}\\
%&\quad+o(\ap^4)\\
%&=\Big[Q(x)+\f{B_{2}}{4B_{1}}\ap^2\big(Q+x\cdot\nabla Q\big)+\ap^4\Big(\f{8B_{1}B_{3}-7B_{2}^2}{32B_{1}^2}\big(Q+x\cdot\nabla Q\big)\\
%&\quad+\f{B_{2}^2}{16B_{1}^2}(x\cdot\nabla Q)+\f{B_{2}^2}{32B_{1}^2}x^T\cdot\nabla^2 Q\cdot x+\psi_{1}(x)\Big)\Big]\Big[1+i(\vaa-\ap)x\cdot x_{a}^\Om\Big]\\
%&\quad+o(\ap^4)\\
%&=Q(x)+\f{B_{2}}{4B_{1}}\ap^2\big(Q+x\cdot\nabla Q\big)+\ap^4\Big[\f{8B_{1}B_{3}-7B_{2}^2}{32B_{1}^2}\big(Q(x)+x\cdot\nabla Q\big)\\
%&\quad+\f{B_{2}^2}{16B_{1}^2}(x\cdot\nabla Q)+\f{B_{2}^2}{32B_{1}^2}x^T\cdot\nabla^2 Q\cdot x+\psi_{1}(x)\Big]+i\f{B_{2}\Om}{8B_{1}}\ap^3x_{1}Q\\
%&\quad+o(\ap^4)\ \ \hbox{as}\ \ a\nearrow a^*,
\end{aligned}
\end{equation*}
This completes the proof of Lemma \ref{lem3.6}.  \qed

Similar to Lemmas \ref{Nlem3.5} and \ref{lem3.6}, one can prove the following lemma with $N=3$, whose proof is omitted for simplicity.
% we next  establish the refined estimates of $\vaa$ and $\mu_a$ for $N=3$, based on which we give the refined spike profile of $w_a(x)$ in terms of $a$.

\begin{lem}\label{lem3.5}
Under the assumptions of Lemma \ref{lem3.1} with $N=3$, let $w_a(x)\in\H$ be a normalized concentrating solution   of \eqref{GP-2} concentrating at $P_0=(0,0,1)$ as $a\searrow0$. Then we have
	\begin{equation*}%\label{3-4}
	\vaa=\ap+o(\ap)\ \ \hbox{and}\ \
	-\mu_{a}\ap^2=\Big(\f{\ap}{\vaa}\Big)^2=1+o(1)\ \ \hbox{as}\ \ a\searrow0,
	\end{equation*}
where $\ap:=\f{a}{a^*}>0$. Moreover, we have
\begin{equation*}%\label{RR3-1}
	\begin{aligned}
	\widetilde{u}_{a}(x):&=\sqrt{a}\ap w_{a}(\ap x+x_{a})e^{-i\widetilde\sigma_{a}-i\ap x\cdot x_{a}^\Om}=Q(x)+o(1)\ \ \hbox{as}\ \ a\searrow0,
	%&=Q(x)+\f{B_{2}}{4B_{1}}\ap^2\big(Q(x)+x\cdot\nabla Q\big)+\ap^4\Big[\f{8B_{1}B_{3}-7B_{2}^2}{32B_{1}^2}\big(Q(x)+x\cdot\nabla Q\big)\\
	%&\quad+\f{B_{2}^2}{16B_{1}^2}(x\cdot\nabla Q)+\f{B_{2}^2}{32B_{1}^2}x^T\cdot\nabla^2 Q\cdot x+\psi_{1}(x)\Big]+i\f{B_{2}\Om}{8B_{1}}\ap^3x_{1}Q\\
	%&\quad+o(\ap^4)\ \ \hbox{as}\ \ a\nearrow a^*,
	\end{aligned}
\end{equation*}
where $\widetilde\sigma_{a}\in[0,2\pi)$ is chosen properly  such that
	\begin{equation*}%\label{RR3-2}
	\operatorname{Re}\Big(\int_{\R^3}  \widetilde{u}_{a}(iQ)\Big)=0.
\end{equation*}
\end{lem}

\section{Axial Symmetry  as $a\to a_*(N)$}
Following the analytical results of Section 3, the purpose of this section is to prove Theorem
\ref{cor3} and Proposition 1.3, for which we always assume that $0\leq\Omega<\infty$ is fixed, and $a>0$ satisfies \eqref{thm1:1}.
Suppose $w_a\in\H$ is a normalized concentrating solution of \eqref{GP-2}
concentrating at $P_0=(0,\cdots,0,1)$ as $a\to a_*(N)$,  where $N=2,3$.   When $N=3$, the uniqueness result is proved under the additional assumption that the maximum point of  $|w_a|$ lies on the $x_3$-axis.

We shall first prove the local uniqueness of Theorem
\ref{cor3}. Towards this aim, we shall prove separately the following Proposition \ref{thm2} for $N=2$ and Proposition \ref{prop:N3-axial-maximum-unique} for $N=3$.

\begin{prop}\label{thm2}
Let $N=2$ and $0\leq\Omega<\infty$ be fixed. For $j=1,2$, suppose
$w_{j,a}\in\H$ is a normalized concentrating solution of \eqref{GP-2}
concentrating at $P_0=(0,1)$ as $a\nearrow a^*$, where $\mu_{j,a}<0$ and
$x_{j,a}=(0,p_{j,a})\in\R^2$ is the unique global maximum point of
$|w_{j,a}|$. Then for all sufficiently small \(a^*-a>0\),
\[
\mu_{1,a}=\mu_{2,a}
\quad\hbox{and}\quad
w_{1,a}(x)\equiv w_{2,a}(\mathcal R x)e^{i\theta}\,\ \hbox{in}\,\ \R^2,
\]
where $\theta\in[0,2\pi)$ is a suitable constant phase and
$\mathcal R:\R^2\mapsto\R^2$ is a suitable rotational transformation.
\end{prop}

 To prove Proposition \ref{thm2}, we argue by contradiction. On the
contrary, suppose there exists a sequence $\{a_n\}$ satisfying \(a_n\nearrow a^*\) as $n\to\infty$ such that
\begin{equation}\label{RR4.1}
	w_{1,a_n}(x)\not\equiv w_{2,a_n}(\mathcal{R}x)e^{i\theta}
\end{equation}
for any rotation $\mathcal{R}: \R^2\mapsto\R^2$ and any constant phase
$\theta\in[0,2\pi)$. For convenience, we always use the parameter \(a\) to denote  \(a_n\) in the following proof.
%```
%
%
%
%
%
%In order to  prove Proposition \ref{thm2},
%on the contrary,  suppose that
%\begin{equation}\label{RR4.1}
%	w_{1,a}(x)\not\equiv w_{2,a}(\mathcal{R}x)e^{i\theta}
%\end{equation}
%for any rotation $\mathcal{R}: \R^2\mapsto\R^2$ and any constant phase $\theta\in[0,2\pi)$.

Denote
\begin{equation}\label{4-4}
	\ap:=\Big(\f{a^*-a}{B_{1}}\Big)^{\f{1}{4}},
\end{equation}
where $B_{1}>0$ is as in Lemma \ref{Nlem3.5}. Following Lemmas \ref{lem3.1} and \ref{Nlem3.5}, the similar argument of  Lemma \ref{Nlem3.0} gives that
\begin{equation}\label{N04-1}
	|x_{1,a}-x_{2,a}|=|p_{1,a}-p_{2,a}|=o(\ap^5)\ \ \ \hbox{as}\ \ a\nearrow a^*.
\end{equation}
Define
\begin{equation}\label{N4.1}
	\widetilde u_{j,a}(x):=\sqrt{a}\ap w_{j,a}\l(\ap x+x_{1,a}\r)e^{-i\widetilde\sigma_{j,a}-i\ap x\cdot x_{1,a}^\Om},\ \ j=1,2,
\end{equation}
where $ \widetilde\sigma_{j,a}\in[0,2\pi)$ is chosen properly (similar to \cite[Lemma 3.2]{GLY}) such that
\begin{equation}\label{N4.2}
	\operatorname{Re}\Big(\int_{\R^2} \widetilde u_{j,a}( iQ)\Big)=0,\ \ j=1,2.
\end{equation}
Using Lemma \ref{lem3.6}, we then derive from \eqref{N04-1} that
\begin{equation}\label{NNR4.1}
	\big|\widetilde u_{1,a}(x)-\widetilde u_{2,a}(x)\big|=o(\ap^4)\ \ \  \hbox{as}\ \ a\nearrow a^*.
\end{equation}

Denote
\begin{equation}\label{4.1}
	u_{1,a}(x):=\widetilde u_{1,a}(x),
	% =
	% \sqrt{a}\ap w_{1,a}\l(\ap x+x_{1,a}\r)e^{-i\sigma_{1,a}-i\ap x\cdot x_{1,a}^\Om},
\end{equation}
and
\begin{equation}\label{4:1}
	\begin{aligned}
		u_{2,a}(x):
		=  \sqrt{a}\ap w_{2,a}\big(\mathcal{R}_0(\ap x+x_{1,a})\big)e^{-i\sigma_{2,a}-i\ap x\cdot x_{1,a}^\Om},\\
	\end{aligned}
\end{equation}
where  $\mathcal{R}_0:\R^2\to\R^2$  is a suitable rotation  such that $\mathcal{R}_0(x_{1,a})$ satisfies
\begin{equation}\label{4:3}
	\l|w_{2,a}(\mathcal{R}_0(x_{1,a}))\r|=\max\limits_{|x|=|x_{1,a}|}|w_{2,a}(x)|,
\end{equation} and similar to \cite[Lemma 3.2]{GLY}, $\sigma_{2,a}\in[0,2\pi)$ is chosen properly such that
\begin{equation}\label{4.2}
	\operatorname{Re}\Big(\int_{\R^2} u_{2,a}( iQ)\Big)=0.
\end{equation}
One can then deduce from \eqref{N4.1} and \eqref{4.1}--\eqref{4:3} that
\begin{equation}\label{4:5}
\frac{\partial |u_{j,a}(x)|}{\partial x_{1}}\Big|_{x=0}=0,
\ \,j=1,2.
\end{equation}

Following above estimates, we deduce from Lemma \ref{lem3.6} and Lemma \ref{lemA.3} below that
\begin{equation}\label{R4:1}
	|u_{1,a}(x)-u_{2,a}(x)|=o(\ap^4)\ \ \ \hbox{as}\ \ a\nearrow a^*.
\end{equation}
% where $C_i(x)$ is as in  Lemma \ref{lem3.6} for $i=1, 2$.
Moreover, it follows from Lemma \ref{Nlem3.5} that
\begin{equation}\label{R3.2}
	\l|(\mu_{1,a}-\mu_{2,a})\ap^2\r|=o(\ap^4)\ \ \ \hbox{as}\ \ a\nearrow a^*.
\end{equation}
Since $w_{1,a}(x)$ and $w_{2,a}(x)$ are two different normalized concentrating solutions of \eqref{GP-2}, we obtain from \eqref{N4.1}, \eqref{4.1} and \eqref{4:1} that $u_{j,a}(x)$ satisfies
\begin{equation}\label{4.3}
	\Big[-\Delta_{\ap^2x^\Om}-\ap^2\mu_{j,a}+\ap^2\l(|\ap x+x_{1,a}|^2-1\r)^2\Big]u_{j,a}-|u_{j,a}|^2u_{j,a}=0\ \,\ \hbox{in}\ \ \R^2,\ \ j=1, 2.
\end{equation}
By the comparison principle (cf. \cite[Lemma A.2]{GLY}), we then derive
that for $j=1, 2$,
\begin{equation}\label{2:11}
	|u_{j,a}(x)|\leq Ce^{-\frac{2}{3}|x|}\ \ \hbox{and}\ \ |\nabla u_{j,a}(x)|\leq Ce^{-\frac{1}{2}|x|}
	\ \ \text{uniformly in}\ \ \R^2\ \ \hbox{as}\ \ a\nearrow a^*.
\end{equation}
Based on \eqref{R4:1}, \eqref{R3.2} and \eqref{2:11}, we next address the proof of Proposition \ref{thm2}.

\vskip 0.05truein

\noindent{\bf Proof of Proposition \ref{thm2}.}
Setting
\begin{equation}\label{4.4}
	\eta_a(x):=\f{u_{1,a}(x)-u_{2,a}(x)}
	{\|u_{1,a}-u_{2,a}\|_{L^\infty(\R^2)}}\quad \hbox{in}\, \ \R^2,
\end{equation}
we shall carry out the proof through the following five steps.

\vskip 0.05truein

{\em Step 1. Up to a subsequence if necessary, \(\eta_a(x)\) satisfies
\begin{equation}\label{RR4-1}
		\eta_a\to\eta_0\quad \hbox{in}\ \  C^1_{loc}(\R^2,\C)
\ \ \hbox{as}\ \  a\nearrow a^*,
\end{equation}
where \(\eta_0(x)\) satisfies
\begin{equation}\label{4.10}
\eta_0=b_0(iQ)+c_0(Q+x\cdot \nabla Q)
		+c_1\f{\p Q}{\p x_1}+c_2\f{\p Q}{\p x_2}
		\ \ \hbox{in}\ \, \R^2
\end{equation}
 for some  \(b_0,c_0,c_1,c_2\in\R\). }

We first deduce from \eqref{4.3} that \(\eta_a(x)\) satisfies the following equation
\begin{equation}\label{4.5}
\begin{split}
&\Big[-\Delta_{\ap^2x^\Om}-\ap^2\mu_{1,a}
+\ap^2\big(|\ap x+x_{1,a}|^2-1\big)^2-|u_{1,a}|^2\Big]\eta_a\\
		&\quad
		-\f{\ap^2(\mu_{1,a}-\mu_{2,a})}
		{\|u_{1,a}-u_{2,a}\|_{L^\infty(\R^2)}}u_{2,a}
		-(u_{1,a}\bar{\eta}_a+\eta_a\bar{u}_{2,a})u_{2,a}=0
		\quad \hbox{in}\ \ \R^2.
\end{split}
\end{equation}
For \(j=1,2\), multiplying \eqref{4.3} by \(\bar u_{j,a}\),
integrating over \(\R^2\), and taking the real part, we obtain that
\begin{equation*}
	\begin{aligned}
		\ap^{2}\mu_{j,a}a
		=&\operatorname{Re}\int_{\R^2}
		\Big\{
		\big(-\Delta_{\ap^2 x^\Om}u_{j,a}\big)\bar{u}_{j,a}
		+\ap^2\big(|\ap x+x_{1,a}|^2-1\big)^2|u_{j,a}|^2
		-|u_{j,a}|^4
		\Big\}dx .
	\end{aligned}
\end{equation*}
Since \(-\Delta_{\ap^2x^\Om}\) is self-adjoint, we obtain from \eqref{4.3} that
\begin{equation}\label{4.8}
	\begin{split}
		&\quad
		\f{\ap^{2}(\mu_{1,a}-\mu_{2,a})a}
		{\|u_{1,a}-u_{2,a}\|_{L^\infty(\R^2)}}\\
		&=\operatorname{Re}\Big(\int_{\R^2}
		\Big[\Big(-\Delta_{\ap^2x^\Om}
		+\ap^2\big(|\ap x+x_{1,a}|^2-1\big)^2\Big)(u_{1,a}+u_{2,a})\\
		&\quad\qquad
		-(|u_{1,a}|^2+|u_{2,a}|^2)(u_{1,a}+u_{2,a})
		\Big]\bar{\eta}_a\Big)\\
		&=\operatorname{Re}\Big(\int_{\R^2}
		\Big[
		(\ap^2\mu_{1,a}+|u_{1,a}|^2)u_{1,a}
		+(\ap^2\mu_{2,a}+|u_{2,a}|^2)u_{2,a}\\
		&\quad\qquad
		-(|u_{1,a}|^2+|u_{2,a}|^2)(u_{1,a}+u_{2,a})
		\Big]\bar{\eta}_a\Big)\\
		&=\operatorname{Re}\Big(\int_{\R^2}
		\Big[
		\ap^2(\mu_{2,a}-\mu_{1,a})u_{2,a}
		-(|u_{1,a}|^2u_{2,a}+|u_{2,a}|^2u_{1,a})
		\Big]\bar{\eta}_a\Big),
	\end{split}
\end{equation}
where the last identity follows from the fact
\begin{align*}
	\operatorname{Re}\Big(\int_{\R^2}
	\ap^2\mu_{1,a}(u_{1,a}+u_{2,a})\bar{\eta}_a\Big)=\operatorname{Re}\Big(\int_{\R^2}
	\ap^2\mu_{1,a}
	\f{|u_{1,a}|^2-|u_{2,a}|^2}
	{\|u_{1,a}-u_{2,a}\|_{L^\infty(\R^2)}}\Big)
	=0.
\end{align*}

Applying \eqref{4-4}, \eqref{R3.2} and \eqref{2:11}, we deduce from
\eqref{4.8} that
\begin{equation}\label{4.8-mu-bound}
	\f{|\ap^{2}(\mu_{1,a}-\mu_{2,a})|}
	{\|u_{1,a}-u_{2,a}\|_{L^\infty(\R^2)}}
	\leq C
	\quad \hbox{uniformly as}\ \ a\nearrow a^*.
\end{equation}
By the comparison principle, we further obtain from \eqref{4.5} that
\begin{equation}\label{R4.4}
	|\eta_{a}(x)|\leq Ce^{-\f{2}{3}|x|},\ \
	|\nabla \eta_{a}(x)|\leq Ce^{-\f{1}{2}|x|}
	\quad \hbox{uniformly in}\ \ \R^2\ \ \hbox{as}\ \ a\nearrow a^*.
\end{equation}
Therefore, the standard elliptic regularity theory yields from
\eqref{4.5} that
\[
\|\eta_a\|_{C^{1,\alpha}_{loc}(\R^2)}\leq C\quad \hbox{uniformly as}\ \ a\nearrow a^*
\]
for some $\al\in(0,1)$. Taking a
subsequence if necessary, we may assume that
\[
\eta_a\to\eta_0\quad \hbox{in}\ \ C^1_{loc}(\R^2,\C)
\ \  \hbox{as}\ \ a\nearrow a^*.
\]
Moreover, using Lemmas \ref{Nlem3.5} and \ref{lem3.6}, we derive from
\eqref{4.5} and \eqref{4.8} that \(\eta_0\) solves the limiting equation
\[
L\eta_0=(-\Delta+1-Q^2)\eta_0-2\operatorname{Re}(Q\eta_0)Q
=-\f{1}{a^*}
\Big(\operatorname{Re}\int_{\R^2}2Q^3\bar{\eta}_0\Big)Q
\quad \hbox{in}\ \ \R^2.
\]
Since
\[
L(Q+x\cdot\nabla Q)=-2Q\quad \hbox{in }\ \R^2,
\]
the non-degeneracy \eqref{nond} of \(L\) yields  that \eqref{4.10} holds true. This proves
Step 1.

\vskip 0.05truein

{\em Step 2. The constants \(b_0\) and \(c_2\) defined by \eqref{4.10}
	satisfy \(b_0=c_2=0\).}

Following \eqref{4.1} and \eqref{4.2}, one can obtain  from \eqref{N4.2} that
\[
\operatorname{Re}\Big(\int_{\R^2}\eta_a(iQ)\,dx\Big)=0.
\]
It then follows from \eqref{RR4-1} that
\begin{equation}\label{4:12}
	\operatorname{Re}\Big(\int_{\R^2}\eta_0(iQ)\,dx\Big)=0.
\end{equation}
Substituting \eqref{4.10} into \eqref{4:12} gives \(b_0=0\).
By the argument of proving \eqref{2.94}, multiplying \eqref{4.3}
with
\[
\frac{\partial \bar u_{j,a}}{\partial x_2}
-\frac{i\ap^2\Om x_1\bar u_{j,a}}{2},
\]
integrating over \(\R^2\), and taking the real part, we obtain that
\begin{equation}\label{4:13}
	0=\int_{\R^2}
	\f{\p \big(|\ap x+x_{1,a}|^2-1\big)^2}{\p x_2}
	|u_{j,a}|^2\,dx,
	\ \ j=1,2.
\end{equation}
By Lemmas \ref{lem3.1}, \ref{lem3.6} and \ref{lemA.3}, we derive from
\eqref{2:11}, \eqref{RR4-1}, \eqref{4.10} and \eqref{4:13} that
\begin{equation*}
	\begin{aligned}
		0
		&=\operatorname{Re}\Big(\int_{\R^2}
		\f{\p \big(|\ap x+x_{1,a}|^2-1\big)^2}{\p x_2}
		(u_{1,a}+u_{2,a})\bar{\eta}_a\Big)\\
		&=\operatorname{Re}\Big(\int_{\R^2}
		4\ap\big(|\ap x+x_{1,a}|^2-1\big)(\ap x_2+p_{1,a})
		(u_{1,a}+u_{2,a})\bar{\eta}_a\Big)\\
		&=\int_{\R^2}(8\ap^2x_2)2Q
		\Big(c_0(Q+x\cdot\nabla Q)
		+c_1\f{\p Q}{\p x_1}
		+c_2\f{\p Q}{\p x_2}\Big)+o(\ap^2)\\
		&=-8c_2\ap^2\int_{\R^2}Q^2+o(\ap^2)
		\quad \hbox{as }a\nearrow a^*.
	\end{aligned}
\end{equation*}
This implies \(c_2=0\), and Step 2 is hence proved.

\vskip 0.05truein

{\em Step 3. The constant \(c_0\) defined by \eqref{4.10} satisfies
	\(c_0=0\).}

Similar to Lemma \ref{lem2.6}, we obtain that for \(j=1,2\),
\begin{equation}\label{4:14}
	\begin{aligned}
		&\operatorname{Re}\Big(\int_{\R^2}2i\ap^2
		(x^\Om\cdot\nabla u_{j,a})(x\cdot\nabla \bar u_{j,a})\Big)\\
		=&-\ap^2\mu_{j,a}\int_{\R^2}|u_{j,a}|^2
		+\ap^2\int_{\R^2}
		\Big[
		2\ap^2|x^\Om|^2
		+\big(|\ap x+x_{1,a}|^2-1\big)^2\\
		&\quad\qquad
		+2\ap\big(|\ap x+x_{1,a}|^2-1\big)
		(\ap |x|^2+x\cdot x_{1,a})
		\Big]|u_{j,a}|^2
		-\frac12\int_{\R^2}|u_{j,a}|^4 .
	\end{aligned}
\end{equation}
Denote
\begin{equation}\label{4:16}
	\begin{aligned}
		B_{j,a}:=&-\ap^2\mu_{j,a}\int_{\R^2}|u_{j,a}|^2\\
		&+\ap^2\int_{\R^2}
		\Big[
		2\ap^2|x^\Om|^2
		+\big(|\ap x+x_{1,a}|^2-1\big)^2\\
		&\quad\qquad
		+2\ap\big(|\ap x+x_{1,a}|^2-1\big)
		(\ap |x|^2+x\cdot x_{1,a})
		\Big]|u_{j,a}|^2,
		\ \  j=1,2,
	\end{aligned}
\end{equation}
\begin{equation}\label{4:17}
	C_{j,a}:=-\frac12\int_{\R^2}|u_{j,a}|^4,
	\quad j=1,2,
\end{equation}
and
\begin{equation}\label{4:18}
	D_{j,a}:=\operatorname{Re}\Big(\int_{\R^2}2i\ap^2
	(x^\Om\cdot\nabla u_{j,a})(x\cdot\nabla \bar u_{j,a})\Big),
	\ \  j=1,2.
\end{equation}
It then yields from \eqref{4:14} that
\begin{equation}\label{4:19}
	\f{D_{1,a}-D_{2,a}}
	{\|u_{1,a}-u_{2,a}\|_{L^\infty(\R^2)}}
	=
	\f{(B_{1,a}-B_{2,a})+(C_{1,a}-C_{2,a})}
	{\|u_{1,a}-u_{2,a}\|_{L^\infty(\R^2)}}.
\end{equation}

We next estimate the terms of \eqref{4:19}. It first  follows from  \eqref{4:18} that
\begin{equation}\label{4:20}
	\begin{aligned}
		&\f {D_{1,a}-D_{2,a}}
		{\|u_{1,a}-u_{2,a}\|_{L^\infty(\R^2)}}\\
		=&\operatorname{Re}\Big(\int_{\R^2}2i\ap^2
		\Big[
		(x^\Om\cdot\nabla\eta_{a})(x\cdot\nabla\bar u_{1,a})
		+(x^\Om\cdot\nabla u_{2,a})(x\cdot\nabla\bar\eta_{a})
		\Big]\Big).
	\end{aligned}
\end{equation}
As for the term containing \(B_{j,a}\), we obtain from \eqref{4.1},
\eqref{4:1} and \eqref{4:16} that
\begin{equation}\label{4:21}
	\begin{aligned}
		&\f {B_{1,a}-B_{2,a}}
		{\|u_{1,a}-u_{2,a}\|_{L^\infty(\R^2)}}\\
		=&\ap^2\operatorname{Re}\Big(\int_{\R^2}
		\Big\{
		2\ap^2|x^\Om|^2
		+\big(|\ap x+x_{1,a}|^2-1\big)^2\\
		&\quad
		+2\ap\big(|\ap x+x_{1,a}|^2-1\big)
		(\ap |x|^2+x\cdot x_{1,a})
		\Big\}
		(\eta_a\bar u_{1,a}+u_{2,a}\bar\eta_a)\Big)\\
		&-\f{\ap^{2}(\mu_{1,a}-\mu_{2,a})a}
		{\|u_{1,a}-u_{2,a}\|_{L^\infty(\R^2)}} .
	\end{aligned}
\end{equation}
We now estimate the term containing \(C_{j,a}\). Note from \eqref{4.8} that
\begin{equation}\label{4:22}
	\begin{split}
		&\quad
		\f{\ap^{2}(\mu_{1,a}-\mu_{2,a})a}
		{\|u_{1,a}-u_{2,a}\|_{L^\infty(\R^2)}}
		+\frac12\operatorname{Re}\Big(\int_{\R^2}
		(|u_{1,a}|^2+|u_{2,a}|^2)(u_{1,a}+u_{2,a})\bar\eta_a\Big)\\
		&=\frac12\operatorname{Re}\Big(\int_{\R^2}
		(|u_{1,a}|^2-|u_{2,a}|^2)(u_{1,a}-u_{2,a})\bar\eta_a\Big)\\
		&\quad+\operatorname{Re}\Big(\int_{\R^2}
		\ap^2(\mu_{2,a}-\mu_{1,a})u_{2,a}\bar\eta_a\Big).
	\end{split}
\end{equation}
We then derive from \eqref{4:17} and \eqref{4:22} that
\begin{equation}\label{4:23}
	\begin{aligned}
		&\f {C_{1,a}-C_{2,a}}
		{\|u_{1,a}-u_{2,a}\|_{L^\infty(\R^2)}}\\
		=&-\frac12\operatorname{Re}\Big(\int_{\R^2}
		(|u_{1,a}|^2+|u_{2,a}|^2)(u_{1,a}+u_{2,a})\bar\eta_a\Big)\\
		=&\f{\ap^{2}(\mu_{1,a}-\mu_{2,a})a}
		{\|u_{1,a}-u_{2,a}\|_{L^\infty(\R^2)}}\\
		&-\frac12\operatorname{Re}\Big(\int_{\R^2}
		(|u_{1,a}|^2-|u_{2,a}|^2)(u_{1,a}-u_{2,a})\bar\eta_a\Big)\\
		&-\operatorname{Re}\Big(\int_{\R^2}
		\ap^2(\mu_{2,a}-\mu_{1,a})u_{2,a}\bar\eta_a\Big).
	\end{aligned}
\end{equation}
Moreover, following Lemmas \ref{lem3.6} and \ref{lemA.3}, we deduce from
\eqref{4:20} that
\begin{equation}\label{R4:3}
	\begin{aligned}
		&\f {D_{1,a}-D_{2,a}}
		{\|u_{1,a}-u_{2,a}\|_{L^\infty(\R^2)}}\\
		=&-\operatorname{Re}\Big(\int_{\R^2}2i\ap^2
		\Big(x^\Om\cdot\nabla
		\big[x\cdot\nabla(Q+\ap^2C_1+\ap^4C_2)\big]\Big)\eta_a\Big)\\
		&+\operatorname{Re}\Big(\int_{\R^2}2i\ap^2
		\Big(x^\Om\cdot\nabla(Q+\ap^2C_1+\ap^4C_2)\Big)
		(x\cdot\nabla\bar\eta_a)\Big)+o(\ap^4)\\
		=&o(\ap^4)\quad \hbox{as }\ a\nearrow a^*,
	\end{aligned}
\end{equation}
due to the radial symmetry of \(Q(x)=Q(|x|)\) and \(C_1(x)=C_1(|x|)\).
Applying Lemmas \ref{lem3.1}, \ref{Nlem3.5}, \ref{lem3.6} and
\ref{lemA.3}, we derive from \eqref{R4:1}, \eqref{R3.2}, \eqref{2:11},
\eqref{RR4-1}, \eqref{4.10}, \eqref{4:19}, \eqref{4:21},
\eqref{4:23} and \eqref{R4:3} that
\begin{equation*}
	\begin{aligned}
		o(\ap^4)
		=&\ap^2\operatorname{Re}\Big(\int_{\R^2}
		\Big\{
		2\ap^2|x^\Om|^2
		+\big(|\ap x+x_{1,a}|^2-1\big)^2\\
		&\quad
		+2\ap\big(|\ap x+x_{1,a}|^2-1\big)
		(\ap |x|^2+x\cdot x_{1,a})
		\Big\}
		(\eta_a\bar u_{1,a}+u_{2,a}\bar\eta_a)\Big)\\
		=&2\ap^4\int_{\R^2}
		\big(2|x^\Om|^2+8x_2^2\big)Q
		\Big[c_0(Q+x\cdot\nabla Q)+c_1\f{\p Q}{\p x_1}\Big]
		+o(\ap^4),
	\end{aligned}
\end{equation*}
where \(c_2=0\) and the radial symmetry of \(Q\) are used. The term
involving \(c_1\) vanishes by the oddness of $\frac{\partial Q}{\partial x_1}$.  Moreover,
since
\[
\int_{\R^2}|x|^2Q(Q+x\cdot\nabla Q)
=-\int_{\R^2}|x|^2Q^2\neq0,
\]
one can obtain that
\[
\int_{\R^2}\Big(2|x^\Om|^2+8x_2^2\Big)Q\big(Q+x\cdot\nabla Q\big)dx\neq0.
\]
Thus, \(c_0=0\) and Step 3 is proved.

\vskip 0.05truein

{\em Step 4. The constant \(c_1\) defined by \eqref{4.10} satisfies
	\(c_1=0\), and hence \(\eta_0\equiv0\) in \(\R^2\).}

Multiplying \eqref{4.4} by \(\overline{u_{1,a}+u_{2,a}}\), we derive that
\begin{equation*}
	\begin{aligned}
		\eta_a\overline{(u_{1,a}+u_{2,a})}
		=&\f{|u_{1,a}|^2-|u_{2,a}|^2
			+u_{1,a}\bar u_{2,a}-u_{2,a}\bar u_{1,a}}
		{\|u_{1,a}-u_{2,a}\|_{L^\infty(\R^2)}}\\
		=&\f{|u_{1,a}|^2-|u_{2,a}|^2}
		{\|u_{1,a}-u_{2,a}\|_{L^\infty(\R^2)}}
		+\eta_a\bar u_{2,a}-u_{2,a}\bar\eta_a
		\quad \hbox{in}\ \ \R^2.
	\end{aligned}
\end{equation*}
Following \eqref{4:5}, we have
\begin{equation}\label{4:8}
	\begin{aligned}
		&\partial_{1}\big[\eta_a\overline{(u_{1,a}+u_{2,a})}\big]\Big|_{x=0}
		=\partial_{1}\big[\eta_a\bar u_{2,a}-u_{2,a}\bar\eta_a\big]\Big|_{x=0}\\
		=&\Big[
		(\partial_1\eta_a)\bar u_{2,a}
		+\eta_a(\partial_1\bar u_{2,a})
		-(\partial_1\bar\eta_a)u_{2,a}
		-(\partial_1u_{2,a})\bar\eta_a
		\Big]\Big|_{x=0}.
	\end{aligned}
\end{equation}
Since \(b_0=0\) holds for \eqref{4.10}, we derive from \eqref{RR4-1} and
Lemma \ref{lemA.3} that
\begin{equation}\label{4:9}
	\hbox{RHS of \eqref{4:8}}\to0
	\quad \hbox{as}\ \ a\nearrow a^*.
\end{equation}
On the other hand,
\begin{equation}\label{4:10-step4}
	\begin{aligned}
		\hbox{LHS of \eqref{4:8}}
		=&\partial_{1}\big[\eta_a\overline{(u_{1,a}+u_{2,a})}\big]\Big|_{x=0}\\
		=&\Big[
		\partial_1\eta_a\overline{(u_{1,a}+u_{2,a})}
		+\eta_a\partial_1\overline{(u_{1,a}+u_{2,a})}
		\Big]\Big|_{x=0}\\
		&\to 2Q(0)\partial_1\eta_0(0)
		\quad \hbox{as}\ \ a\nearrow a^*.
	\end{aligned}
\end{equation}
We then derive from \eqref{4:9} and \eqref{4:10-step4} that
\begin{equation}\label{4:11}
	\partial_1\eta_0(0)=0.
\end{equation}
Substituting \eqref{4.10} into \eqref{4:11}, and using
\(b_0=c_0=c_2=0\), we obtain that
\[
c_1\f{\p^2 Q(0)}{\p x_1^2}=c_1Q''(0)=0.
\]
Since \(Q''(0)<0\), it follows that \(c_1=0\). Hence, \(\eta_0\equiv0\) in \(\R^2\) and Step 4 is proved.

\vskip 0.05truein

{\em Step 5. The assumption \eqref{RR4.1} is false. }

We derive from above steps that \(\eta_0\equiv0\) holds in \(\R^2\).
However, since
\[
\|\eta_a\|_{L^\infty(\R^2)}=1,
\]
 the exponential decay  of \eqref{R4.4} gives that the maximum of \(|\eta_a|\) cannot escape to
infinity. Thus,  a subsequence of \(\{\eta_a\}\) converges to a nonzero limit
in \(C^1_{loc}(\R^2,\C)\), which  however   contradicts  the fact that \(\eta_0\equiv0\) holds in \(\R^2\). Therefore,
the assumption \eqref{RR4.1} is false.

We therefore derive that for sufficiently small \(a^*-a>0\),
\[
w_{1,a}(x)\equiv w_{2,a}(\mathcal R x)e^{i\theta}
\ \  \hbox{in}\ \, \R^2
\]
holds for some suitable constant phase \(\theta\in[0,2\pi)\) and some rotation
\(\mathcal R:\R^2\mapsto\R^2\).  It then yields from
\eqref{GP-2} that
\[
\mu_{1,a}=\mu_{2,a}
\]
holds for sufficiently small \(a^*-a>0\).   Proposition
\ref{thm2} is therefore proved. \qed

\begin{prop}\label{prop:N3-axial-maximum-unique}
Let $N=3$ and $0\leq\Omega<\infty$ be fixed. For $j=1,2$, suppose
$w_{j,a}\in\H$ is a normalized concentrating solution of \eqref{GP-2}
concentrating at $P_0=(0,0,1)$ as $a\searrow0$, where
$\mu_{j,a}<0$, and the maximum point $x_{j,a}$ of $|w_{j,a}|$
satisfies
\[
x_{j,a}=p_{j,a}e_3=(0,0,p_{j,a}),\ \  p_{j,a}\to1
\ \ \hbox{as}\ \  a\searrow 0,\ \  j=1,2.
\]
Then for  sufficiently small $a>0$, there exists a constant
phase $\theta_a\in[0,2\pi)$ such that
\[
w_{1,a}(x)\equiv e^{i\theta_a}w_{2,a}(x)
\ \ \hbox{in}\ \ \R^3,
\]
and $\mu_{1,a}=\mu_{2,a}$.
\end{prop}

\noindent{\bf Proof.}
By contradiction, suppose  there exists a sequence $\{a_n\}$ satisfying
	\(a_n\searrow0\) as $n\to\infty$ such that
\begin{equation}\label{assum1}
w_{1,a_n}\not\equiv e^{i\theta}w_{2,a_n}
\ \  \hbox{for any}\ \ \theta\in[0,2\pi).
\end{equation}
For convenience, next we   still denote the sequence $\{w_{i,a_n}\}$ by $\{w_{i,a}\}$ for $i=1,2$.

Recall that for $N=3$,
\[
\alpha_{a}=\frac{a}{a^*}.
\]
Since $x_{j,a}=p_{j,a}e_3$, we have $x_{j,a}^{\Omega}=0$. Define
\begin{equation}\label{eq:N3-new-uja}
	u_{j,a}(x):=\sqrt {a}\,\alpha_a\,w_{j,a}(\alpha_a x+p_{j,a}e_3)e^{-i\sigma_{j,a}},
	\ \  j=1,2,
\end{equation}
where $\sigma_{j,a}\in[0,2\pi)$ is chosen so that
\begin{equation}\label{eq:N3-phase-condition}
	\operatorname{Re}\int_{\R^3}u_{j,a}(iQ)\,dx=0,
	\ \  j=1,2.
\end{equation}
Applying the concentration estimates of Section 3,  we obtain that
\begin{equation}\label{eq:N3-uja-Q}
u_{j,a}\to Q\quad \hbox{in}\ \ C^1_{\rm loc}(\R^3,\C)\cap L^\infty(\R^3,\C)\ \ \hbox{as}\ \  a\searrow0,
\ \  j=1,2,
\end{equation}
and $u_{j,a}$ satisfies the uniform exponential decay  \eqref{2:11}. Moreover,
\begin{equation}\label{eq:N3-mass-uja}
	\int_{\R^3}|u_{j,a}|^2\,dx=a^*,
	\ \  j=1,2.
\end{equation}
Set
\[
V_q(x):=\big(|\alpha_a x+qe_3|^2-1\big)^2,
\]
so that $u_{j,a}$ satisfies
\begin{equation}\label{eq:N3-uja-equation}
	\Big[-\Delta_{\alpha_a^2x^\Omega}-\alpha_a^2\mu_{j,a}
	+\alpha_a^2V_{p_{j,a}}(x)\Big]u_{j,a}-|u_{j,a}|^2u_{j,a}=0
	\ \  \hbox{in}\ \ \R^3,\ \  j=1,2.
\end{equation}

Denote
\[
\delta_a:=\|u_{1,a}-u_{2,a}\|_{L^\infty(\R^3)}\ge 0.
\]
We claim that \(\delta_a>0\). On the contrary, suppose \(\delta_a=0\), then \(u_{1,a}\equiv u_{2,a}\). This implies that
\[
w_{1,a}(y)=e^{i\theta}w_{2,a}(y+h_ae_3),
\ \  h_a=p_{2,a}-p_{1,a}.
\]
Compare the equations of $w_{1,a}$ and $w_{2,a}$ through (\ref{eq:N3-new-uja}) and $(\ref{eq:N3-uja-equation})$. By the fact \((y+h_ae_3)^\Omega=y^\Omega\),  it then yields that
\[
\left[(|y|^2-1)^2-(|y+h_ae_3|^2-1)^2\right]w_{1,a}
=(\mu_{1,a}-\mu_{2,a})w_{1,a}\ \  \hbox{in}\ \ \R^3.
\]
Since \(w_{1,a}\) is nonzero on a  non-empty  open set, the polynomial in  the
bracket of the above equation  must be a constant. This is impossible unless \(h_a=0\). Thus,
\(w_{1,a}\equiv e^{i\theta}w_{2,a}\) in $\R^3$,  which however contradicts \eqref{assum1}.
Therefore, the claim \(\delta_a>0\) holds true.

We next claim that there exists $C>0$, independent of $a>0$, such that
\begin{equation}\label{eq:N3-p-difference}
	|p_{1,a}-p_{2,a}|\leq C\alpha_a\delta_a\ \   \hbox{as}\ \ a\searrow0.
\end{equation}
Indeed, applying the \(x_3\)-direction Pohozaev identity  of
\eqref{eq:N3-uja-equation}, together with
\[
\partial_{x_3}x^\Omega=0,\ \  (x^\Omega)_3=0,
\]
we obtain that
\begin{equation}\label{eq:N3-x3-poho}
	\int_{\R^3}\partial_{x_3}V_{p_{j,a}}(x)|u_{j,a}|^2\,dx=0,
	\ \  j=1,2.
\end{equation}
Define
\[
F_q(x):=\big(|\alpha_a x+qe_3|^2-1\big)(\alpha_a x_3+q).
\]
Since
\[
\partial_{x_3}V_q(x)=4\alpha_a\big(|\alpha_a x+qe_3|^2-1\big)(\alpha_a x_3+q),
\]
one can derive from (\ref{eq:N3-x3-poho}) that
\[
\int_{\R^3}F_{p_{j,a}}(x)|u_{j,a}|^2\,dx=0,
\ \ j=1,2,
\]
which further gives that
\begin{equation}\label{eq:N3-subtract-F}
	0=\int_{\R^3}(F_{p_{1,a}}-F_{p_{2,a}})|u_{1,a}|^2\,dx
	+\int_{\R^3}F_{p_{2,a}}(|u_{1,a}|^2-|u_{2,a}|^2)\,dx .
\end{equation}
Note that
\[
F_{p_{1,a}}-F_{p_{2,a}}
=(p_{1,a}-p_{2,a})\int_0^1
\partial_qF_{p_{2,a}+t(p_{1,a}-p_{2,a})}\,dt,
\]
where
\[
\partial_qF_q(x)=2(\alpha_a x_3+q)^2+\big(|\alpha_a x+qe_3|^2-1\big).
\]
Since $p_{j,a}\to1$ as $a\searrow0$, we derive from \eqref{eq:N3-uja-Q} that
\begin{equation}\label{eq:N3-positive-coeff}
	\begin{split}
	 \int_{\R^3}\left(\int_0^1\partial_qF_{p_{2,a}+t(p_{1,a}-p_{2,a})}\,dt\right)
	|u_{1,a}|^2\,dx=2\int_{\R^3}Q^2\,dx+o(1)>0\ \ \hbox{as}\ \ a\searrow0.
	\end{split}
\end{equation}
On the other hand, the above definition of $F_q(x)$ gives that
\begin{equation}\label{eq:N3-F-expansion}
	F_q(x)=q(q^2-1)+\alpha_a(3q^2-1)x_3
	+\alpha_a^2q(|x|^2+2x_3^2)+\alpha_a^3|x|^2x_3.
\end{equation}
Note that the constant term $p_{2,a}(p_{2,a}^2-1)$ gives no contribution to the
second integral of \eqref{eq:N3-subtract-F},  due to  the mass identity
\eqref{eq:N3-mass-uja}. Thus,  the uniform exponential decay of
$u_{j,a}$ yields that
\begin{equation}\label{eq:N3-second-F-bound}
	\begin{aligned}
		\left|\int_{\R^3}F_{p_{2,a}}(|u_{1,a}|^2-|u_{2,a}|^2)\,dx\right|
		&\leq C\alpha_a\int_{\R^3}(1+|x|^3)(|u_{1,a}|+|u_{2,a}|)|u_{1,a}-u_{2,a}|\,dx \\
		&\leq C\alpha_a\delta_a.
	\end{aligned}
\end{equation}
We hence combine \eqref{eq:N3-subtract-F}--\eqref{eq:N3-second-F-bound} to obtain that the claim
\eqref{eq:N3-p-difference} holds true.

Define
\begin{equation}\label{eq:N3-eta-def}
	\eta_a(x):=\frac{u_{1,a}(x)-u_{2,a}(x)}{\delta_a}.
\end{equation}
By \eqref{eq:N3-p-difference}, we have
\begin{equation}\label{eq:N3-potential-error}
	\left|
	\frac{\alpha_a^2(V_{p_{1,a}}-V_{p_{2,a}})}{\delta_a}
	\right|
	\leq C\alpha_a^4(1+|x|^3)
	\quad\hbox{in } \R^3 .
\end{equation}
Repeating the analysis of \eqref{4.8} and \eqref{4.8-mu-bound}, which contains an additional
potential-difference term controlled by \eqref{eq:N3-potential-error}, we obtain that
\[
\left|\frac{\alpha_a^2(\mu_{1,a}-\mu_{2,a})}{\delta_a}\right|
\leq
C\delta_a
\left|\frac{\alpha_a^2(\mu_{1,a}-\mu_{2,a})}{\delta_a}\right|
+C.
\]
Since $\delta_a\to 0$ as $ a\searrow0$, this  further yields that
\begin{equation}\label{eq:N3-mu-bound}
	\left|\frac{\alpha_a^2(\mu_{1,a}-\mu_{2,a})}{\delta_a}\right|\le C
	\quad\hbox{uniformly as }\ a\searrow0.
\end{equation}

Subtracting the equations (\ref{eq:N3-uja-equation}) of $u_{1,a}$ and $u_{2,a}$, we  obtain that
\begin{equation}\label{eq:N3-eta-equation}
	\begin{aligned}
		&\Big[-\Delta_{\alpha_a^2x^\Omega}-\alpha_a^2\mu_{1,a}
		+\alpha_a^2V_{p_{1,a}}-|u_{1,a}|^2\Big]\eta_a \\
		& -\frac{\alpha_a^2(\mu_{1,a}-\mu_{2,a})}{\delta_a}u_{2,a}
		+\frac{\alpha_a^2(V_{p_{1,a}}-V_{p_{2,a}})}{\delta_a}u_{2,a} \\
		& -(u_{1,a}\overline{\eta}_a+\eta_a\overline{u}_{2,a})u_{2,a}=0
		\ \ \hbox{in}\ \ \R^3.
	\end{aligned}
\end{equation}
Similar to  Step 1
in the proof of Proposition \ref{thm2}, by  the comparison
principle and standard elliptic estimates, one can derive from  \eqref{eq:N3-potential-error}--(\ref{eq:N3-eta-equation}) that
\begin{equation}\label{eq:N3-eta-decay}
	|\eta_a(x)|\leq Ce^{-\frac23|x|},\ \
	|\nabla\eta_a(x)|\leq Ce^{-\frac12|x|}\ \ \hbox{uniformly in}\ \ \R^3,
\end{equation}
and up to a subsequence if necessary,
\begin{equation}\label{eta-con}
	\eta_a\to\eta_0 \ \ \hbox{in}\ \ C^1_{\rm loc}(\R^3,\C)\ \ \hbox{as}\ \ a\searrow0,
\end{equation}
where
\begin{equation}\label{eq:N3-eta0-form}
	\eta_0=b_0(iQ)+c_0(Q+x\cdot\nabla Q)+\sum_{k=1}^3c_k\frac{\partial Q}{\partial x_k}
	\ \ \hbox{in}\ \ \R^3
\end{equation}
holds for some $b_0,c_0,c_1,c_2,c_3\in\R$.

We now prove that $b_0=c_0=0$ holds for $(\ref{eq:N3-eta0-form})$. In fact, we deduce from
\eqref{eq:N3-phase-condition}, \eqref{eq:N3-eta-decay} and \eqref{eta-con} that
\[
\operatorname{Re}\int_{\R^3}\eta_0(iQ)\,dx=0,
\]
which implies that  $b_0=0$ in view of \eqref{eq:N3-eta0-form}. Note from \eqref{eq:N3-mass-uja} that
\[
0=\frac1{\delta_a}\int_{\R^3}(|u_{1,a}|^2-|u_{2,a}|^2)\,dx
=\operatorname{Re}\int_{\R^3}(u_{1,a}+u_{2,a})\overline{\eta}_a\,dx,
\]
which further  gives that
\[
0=2\int_{\R^3}Q\operatorname{Re}\eta_0\,dx.
\]
Since
\[
\int_{\R^3}Q(Q+x\cdot\nabla Q)\,dx
=-\frac12\int_{\R^3}Q^2\,dx\ne0,
\]
and
\[
\int_{\R^3}Q\frac{\partial Q}{\partial x_k}\,dx=0,
\ \  k=1,2,3,
\]
we obtain that $c_0=0$.

We next prove that $c_1=c_2=c_3=0$ holds for $(\ref{eq:N3-eta0-form})$, and hence  $\eta_0\equiv0$ in $\R^3$. Since $0$ is a maximum point of $|u_{j,a}|$, we have
\[
\nabla |u_{j,a}|^2(0)=0,
\ \  j=1,2,
\]
which gives that
\[
\operatorname{Re}\Big(
\nabla\eta_a(0)\overline{u}_{1,a}(0)
+\nabla u_{2,a}(0)\overline{\eta}_a(0)
\Big)=0.
\]
Noting that $u_{j,a}\to Q$ and  $\nabla u_{2,a}(0)\to0$ as $a\searrow0$,
we then deduce that
\[
\operatorname{Re}\nabla\eta_0(0)=0.
\]
Since $b_0=c_0=0$, it follows from \eqref{eq:N3-eta0-form} that
\[
\eta_0=\sum_{k=1}^3c_k\frac{\partial Q}{\partial x_k}.
\]
Thus, for $l=1,2,3$,
\[
0=\operatorname{Re}\partial_{x_l}\eta_0(0)
=\sum_{k=1}^3c_k\frac{\partial^2Q(0)}{\partial x_l\partial x_k}.
\]
Since $Q=Q(|x|)$, we have
\[
\frac{\partial^2Q(0)}{\partial x_l\partial x_k}=Q''(0)\delta_{lk},
\]
and hence $c_1=c_2=c_3=0$ in view of $Q''(0)<0$, which further implies that  $\eta_0\equiv0$ in $\R^3$.

It follows from  \eqref{eq:N3-eta-decay} and
$\|\eta_a\|_{L^\infty(\R^3)}=1$  that the maximum of $|\eta_a|$    in $\R^3$ is
attained in a fixed ball. Hence, up to a subsequence if necessary, $\eta_a$ converges to a
nonzero limit in $C^1_{\rm loc}(\R^3,\C)$, which however contradicts the fact that  $\eta_0\equiv0$ in $\R^3$.
Therefore, the assumption \eqref{assum1} is false, and there exists a constant phase
$\theta_a\in[0,2\pi)$ such that for all sufficiently small $a>0$,
\[
w_{1,a}\equiv e^{i\theta_a}w_{2,a}\ \ \hbox{in}\ \ \R^3,
\]
which  further implies from \eqref{GP-2} that  $\mu_{1,a}=\mu_{2,a}$ holds for all sufficiently small $a>0$.
This completes the proof of Proposition  \ref{prop:N3-axial-maximum-unique}. \qed

\subsection{Proofs of Theorem \ref{cor3} and Proposition \ref{cor4}}

Following Theorem \ref{thm1}, Proposition \ref{thm2}, and Proposition
\ref{prop:N3-axial-maximum-unique}, we are now ready to prove Theorem
\ref{cor3} and Proposition \ref{cor4}.

\vskip 0.05truein
\noindent{\bf Proof of Theorem \ref{cor3}.}
Let $w_a\in\H$ be a normalized concentrating solution of \eqref{GP-2}
concentrating at $P_0=(0,\cdots,0,1)$ as $a\to a_*(N)$, and suppose $w_{0,a}$ is
the axially symmetric concentrating solution constructed in Theorem
\ref{thm1}.

If $N=2$, then Proposition \ref{thm2} implies that for $a$ sufficiently close to
$a^*$, there exist a constant phase $\theta\in[0,2\pi)$ and a rotation
$\mathcal R: \R^2\to\R^2$ such that
\begin{equation}\label{04-1-N2}
	w_a(x)\equiv w_{0,a}(\mathcal R x)e^{i\theta}
	\ \  \hbox{in}\ \ \R^2.
\end{equation}
Since $w_{0,a}$ is axially symmetric, Theorem \ref{cor3} is proved.

If $N=3$, then we assume additionally that the maximum point of $|w_a|$ lies on the
$x_3$-axis. The solution $w_{0,a}$ constructed in Theorem \ref{thm1} is
axially symmetric and has its maximum point on the $x_3$-axis. It then yields from
Proposition \ref{prop:N3-axial-maximum-unique}   that for sufficiently small $a>0,$
\begin{equation}\label{04-1-N3}
	w_a(x)\equiv e^{i\theta_a}w_{0,a}(x)
	\ \  \hbox{in}\ \ \R^3
\end{equation}
holds for some $\theta_a\in[0,2\pi)$. Hence, $w_a$ is axially symmetric. This therefore
completes the proof of Theorem \ref{cor3}. \qed

\vskip 0.05truein
\noindent{\bf Proof of Proposition \ref{cor4}.}
For $N=3$,  let $w_a\in\H$ be a normalized concentrating solution of
\eqref{GP-2} concentrating at $P_0=(0,0,1)$ as $a\searrow0$, where we also assume that the maximum point of $|w_a|$ lies on the $x_3$-axis. It follows from
Theorem \ref{cor3} that  up to a constant phase, $w_a$ coincides with the
axially symmetric solution constructed in Theorem \ref{thm1}. Therefore, it
suffices to prove the assertion for this axially symmetric solution.

For the axially symmetric solution $w_a$, we have
\[
x^\Omega\cdot\nabla w_a\equiv0\ \ \hbox{in}\ \ \R^3,
\]
and hence the equation \eqref{GP-2} becomes
\begin{equation}\label{04-3}
	-\Delta w_a+(V(x)+1)w_a-a|w_a|^2w_a=\mu_a w_a
	\ \ \hbox{in}\ \ \R^3,
\end{equation}
where $V(x)$ is defined by \eqref{pot}. Let $x_a$ be a maximum point of
$|w_a|$, and set
\[
u_a(x)=\sqrt a\,\vaa\,w_a(\vaa x+x_a)e^{-i\sigma_a-i\vaa x\cdot x_a^\Omega}
=:R_a(x)+iI_a(x),
\ \  \vaa:=\frac1{\sqrt{-\mu_a}},
\]
where the constant phase $\sigma_a\in[0,2\pi)$ is chosen so that
\[
\operatorname{Re}\int_{\R^3}u_a(iQ)\,dx=0.
\]
It yields from Definition \ref{def2} that $R_a\to Q$ and $I_a\to0$ in $H^1(\R^3)$ as
$a\searrow0$. Taking the imaginary part of the equation for $u_a$, we obtain that
\begin{equation}\label{4.11}
	\mathcal L_a I_a=0\ \ \hbox{in }\ \R^3,
	\ \  \int_{\R^3}I_aQ\,dx=0,
\end{equation}
where
\[
\mathcal L_a=-\Delta+\vaa^2|x_a^\Omega|^2
+1+\vaa^2\big[V(\vaa x+x_a)+1\big]-|u_a|^2.
\]
Set
\begin{equation}\label{4.12}
	\mathcal L:=-\Delta+1-Q^2\quad \hbox{in }\ L^2(\R^3).
\end{equation}
The spectral property of \(\mathcal L\) gives that, under the
orthogonality condition in \eqref{4.11},
\[
(\mathcal L I_a,I_a)\geq c\|I_a\|_{L^2(\R^3)}^2
\]
holds for  sufficiently small $a>0$. Multiplying
\eqref{4.11} by $I_a$ and integrating over $\R^3$, we then get that  for  sufficiently small $a>0$,
\[
0=(\mathcal L_aI_a,I_a)
\geq (\mathcal LI_a,I_a)-o(1)\|I_a\|_{L^2(\R^3)}^2
\geq \frac c2\|I_a\|_{L^2(\R^3)}^2,
\]
which implies that $I_a\equiv0$ holds for  sufficiently small $a>0$. Thus, up to a constant phase, $w_a$ is real-valued
in $\R^3$ and hence $w_a$ is free of vortices  for  sufficiently small $a>0$. The proof of Proposition
\ref{cor4} is therefore complete. \qed

\section*{Statements and Declarations}

\noindent\textbf{Data availability statement.}
No datasets were generated or analyzed during the current study.

\medskip

%\noindent\textbf{Author contribution statement.}
%All authors contributed to the formulation of the problem, the development of the mathematical analysis, and the preparation of the manuscript. All authors read and approved the final manuscript.

%\medskip
%
%\noindent\textbf{Competing interests.}
%The authors declare that they have no competing interests.

\section*{Acknowledgments}
This paper is supported by the National Key R \& D Program of China
(2022YFA1006900). Y. Guo is partially supported by NSFC under Grants
12225106 and 11931012. Y. Luo is partially supported by NSFC under Grant
12571119. S. Peng is partially supported by the Key Project of NSFC under
Grant 11831009.

\appendix
\section{Appendix}
\renewcommand{\theequation}{A.\arabic{equation}}
\setcounter{equation}{0}
For the reader's convenience, in this appendix we shall give the detailed proofs of the following results used in this paper.

\begin{lem}\label{lemA.1}
	Assume that $0\leq\Omega<\infty$ is fixed, and $H_\va$ is the Hilbert space defined by \eqref{2.83}, where $\va>0$ is small. Then there exists a constant $C>0$, independent of $\va$, such that for any $u\in H_\va$,
	\begin{equation*}%\label{A1}
	\|u\|^2_{H^1(\R^N)}\leq C\|u\|_\va^2,
	\end{equation*}
	where
$\|u\|^2_\va=\langle u,u\rangle_\va$
is defined by \eqref{2.7},
and $x_{\va}\in\R^N$ satisfies \eqref{2.1}.
\end{lem}
\noindent{\bf Proof.}
We first note that for any fixed $0\leq\Om<+\infty$, there exists a constant $C_{1}>0$ such that
\begin{equation*}
|x^\Om|^2\leq C_{1}\Big[(|x|^2-1)^2+1\Big]\ \ \hbox{in}\ \ \R^N,
\end{equation*}
where $x^\Om$ is as in \eqref{1.1}.
This implies that for sufficiently small $\va>0$,
\begin{equation*}
|(\va x+x_{\va})^\Om|^2\leq C_{1}\Big[\l(|\va x+x_{\va}|^2-1\r)^2+1\Big]\ \ \hbox{in}\ \ \R^N.
\end{equation*}
This gives that there exists a constant $C_{2}>0$ such that for sufficiently small $\va>0$,
\begin{equation*}
\va^2|x^\Om|^2\leq C_{2}\Big[\l(|\va x+x_{\va}|^2-1\r)^2+1\Big]\ \ \hbox{in}\ \ \R^N,
\end{equation*}
due to the fact that $x_{\va}\to P_{0}=(0,\cdots,0,1)$ as $\va\to0$.
By Cauchy's inequality, we derive from \eqref{2.7} that for sufficiently small $\va>0$,
\begin{equation}\label{A.1}
\begin{aligned}
\|u\|_\va^2
% =&Re\l(\inte \l\{(\nabla u-i\va^2 x^\Om u)\ovl{(\nabla u-i\va^2 x^\Om u)}
% +\big(1+\va^2 V_{\Om}(\va x+x_{\va})\big)|u|^2\r\}\r)\\
=&\operatorname{Re}\Big(\inte\Big\{|\nabla u|^2-2\va^2 x^\Om(iu,\nabla u)\\
&\quad\quad\qquad+\Big[\va^4|x^\Om|^2+1+\va^2\l(|\va x+x_{\va}|^2-1\r)^2\Big]|u|^2\Big\}\Big)\\
\geq& \inte\Big\{\fh|\nabla u|^2-2\va^4|x^\Om|^2|u|^2+|u|^2\Big\}\\
\geq& \inte\Big\{\fh|\nabla u|^2-C_{3}\Big[\va^2\l(|\va x+x_{\va}|^2-1\r)^2+1\Big]|u|^2+|u|^2\Big\},\\
\end{aligned}
\end{equation}
where $C_{3}>0$ is a constant.

On the other hand, applying the diamagnetic inequality \eqref{1:4}, we obtain from \eqref{2.7} that
\begin{equation*}
\|u\|_\va^2\geq \inte \Big[1+\va^2\l(|\va x+x_{\va}|^2-1\r)^2\Big]|u|^2.
\end{equation*}
Together with \eqref{A.1}, we then deduce that for sufficiently small $\va>0$,
\begin{equation*}
\|u\|_\va^2\geq \f{1}{2(1+C_{3})}\inte \l(|\nabla u|^2+|u|^2\r)\geq C\|u\|_{H^1(\R^N)}^2.
\end{equation*}
This completes the proof of Lemma \ref{lemA.1}.\qed

\begin{lem}\label{lemA.3}
For \(N=2\),  let \(\mathcal{R}_0\) be the rotational
transformation defined in \eqref{4:3}. Then
\begin{equation}\label{A.4}
	u_{2,a}(x)
	=
	Q(x)+\ap^2C_1(x)+\ap^4C_2(x)+o(\ap^4)
	\ \ \hbox{as }\ a\nearrow a^* .
\end{equation}
\end{lem}

\noindent{\bf Proof.}
It follows from \eqref{4:3} that $\mathcal{R}_0: \R^2 \rightarrow \R^2$ is a suitable rotation  such that $\mathcal{R}_0(x_{1,a})$ satisfies
\begin{equation}\label{A4:3}
  \l|w_{2,a}(\mathcal{R}_0(x_{1,a}))\r|=\max\limits_{|x|=|x_{1,a}|}|w_{2,a}(x)|,
\end{equation}
where $w_{2,a}$ is as in \eqref{RR4.1}, and $x_{1,a}\in\R^2$ is the unique maximum point of $|w_{1,a}|$ as $a\to a_*(2)$.
We first claim that
\begin{equation}\label{NN4:1}
\l|\mathcal{R}_0 (x_{1,a})-x_{1,a}\r|=
o(\ap^5)|x_{1,a}|\ \ \hbox{as}\ \ a\nearrow a^*.
\end{equation}
% from which one can deduce that $u_{2,a}(x)$ also satisfies the refined spike profile \eqref{NNR4.1}.
Indeed,  on the contrary, assume that
\begin{equation}\label{NN4:2}
\l|\mathcal{R}_0 (x_{1,a})-x_{1,a}\r|\geq C\ap^5|x_{1,a}|\ \ \hbox{as}\ \ a\nearrow a^*,
\end{equation}
where $C>0$ is a constant. Since $x_{2,a}\in\R^2$ is the unique maximum point of $|w_{2,a}|$ as $a\to a_*(2)$, we derive from \eqref{N04-1} and \eqref{NN4:2} that
$$\big|w_{2,a}(x_{1,a})\big|>\big|w_{2,a}\big(\mathcal{R}_0 (x_{1,a})\big)\big|,$$
which is however a contradiction in view of \eqref{A4:3}. The claim \eqref{NN4:1} thus holds true.

Since  $\mathcal{R}_0$  is  a rotational transformation, we also have  for any  $x\in\R^2$,
\begin{equation}\label{NN4:3}
\l|\mathcal{R}_0 (x)-x\r|=o(\ap^5)|x|\ \ \hbox{as}\ \ a\nearrow a^*,
\end{equation}
Similar to the argument of \eqref{R3.10}, one then deduces from \eqref{NNR4.1}, \eqref{4:1}, \eqref{4.2} and \eqref{NN4:3} that
\begin{equation}\label{A:4}
|\sigma_{2,a}-\widetilde\sigma_{2,a}|=
o(\ap^4)\ \ \hbox{as}\ \ a\nearrow a^*,
\end{equation}
Following Lemma \ref{lem3.6}, we thus derive from \eqref{NNR4.1},
\eqref{NN4:3} and \eqref{A:4} that \eqref{A.4} holds true. This completes the proof of Lemma \ref{lemA.3}. \qed

\end{document}